\documentclass[11pt]{elsarticle}
\RequirePackage{snapshot}
\usepackage[utf8]{inputenc}
\usepackage[T1]{fontenc}
\usepackage{graphicx}
\usepackage{longtable}
\usepackage{float}
\usepackage{wrapfig}
\usepackage[normalem]{ulem}
\usepackage{amssymb}
\usepackage{amstext}
\usepackage{hyperref}
\usepackage{amsmath}
\usepackage{amsthm}
\usepackage{tabularx}
\usepackage{xfrac}
\usepackage{array}
\usepackage{appendix} 
\usepackage[capitalize]{cleveref}
\crefname{appsec}{Appendix}{Appendices}
\crefname{lemma}{Lemma}{Lemmata}
\Crefname{lemma}{Lemma}{Lemmata}
\crefname{thm}{theorem}{theorems}
\Crefname{thm}{Theorem}{Theorems}
\crefname{alg}{algorithm}{algorithms}
\Crefname{alg}{Algorithm}{Algorithms}
\crefname{table}{Table}{Tables}
\Crefname{table}{Table}{Tables}

\newcommand{\Bezier}{{B\'{e}zier} }

\usepackage{scalerel}
\usepackage{setspace} 

\DeclareMathOperator*{\assembly}{\scalerel*{\text{\sf \Huge A}}{\sum}}

\usepackage[margin=1in]{geometry}
\usepackage[table]{xcolor}

\usepackage[numbers]{natbib}
\usepackage{pdfcomment}
\usepackage{tikz}
\usepackage{arydshln}
\newenvironment{lbmatrix}[1]
  {\left[\array{@{}*{#1}{c}@{}}}
  {\endarray\right]}

\makeatletter
  \renewcommand*\env@matrix[1][*\c@MaxMatrixCols c]{%
    \hskip -\arraycolsep
    \let\@ifnextchar\new@ifnextchar
  \array{#1}}
\makeatother
\makeatletter
\renewcommand*\env@matrix[1][*\c@MaxMatrixCols c]{%
  \hskip -\arraycolsep
  \let\@ifnextchar\new@ifnextchar
  \array{#1}}
\makeatother
\newcommand{\tikzmark}[2]{
    \tikz[overlay,remember picture,baseline] 
    \node[anchor=base] (#1) {$#2$};
  }
  \usepackage{overpic}

\usepackage{pgfplots}

\usepackage{bm}
\usepackage{mathrsfs}
\usepackage{multirow}
\usepackage{nicefrac}
\usepackage[section]{placeins}
\usepackage{subcaption}
\biboptions{numbers,sort&compress} 
\begin{document}
\begin{frontmatter}

\title{Isogeometric B\'ezier dual mortaring: Refineable higher-order spline dual bases and weakly continuous geometry}

\author[byu]{Z. Zou\corref{cor}}
\ead{zhihuizou1988@gmail.com}
\cortext[cor]{Corresponding author}
\author[byu]{M. A. Scott}

\author[byu]{M. J. Borden}

\author[coreform]{D. C. Thomas}

\author[TUK]{W. Dornisch}

\author[Pavia]{E. Brivadis}

\address[byu]{Department of Civil and Environmental Engineering,
  Brigham Young University,
  Provo, Utah 84602, USA}

\address[coreform]{Coreform LLC, P.O. Box 970336, Orem, Utah 84097, USA} 

\address[TUK]{Lehrstuhl f{\"u}r Technische Mechanik, Technische Universit{\"a}t Kaiserslautern,
  Gottlieb-Daimler-Str., Kaiserslautern 67663, Germany}

\address[Pavia]{Istituto di Matematica Applicata e Tecnologie Informatiche del
  CNR, Via Ferrata 1, Pavia 27100, Italy}

\begin{abstract}
 In this paper we develop the isogeometric \Bezier dual mortar method. It is based on \Bezier extraction and projection and is applicable to any spline space which can be represented in \Bezier form (i.e., NURBS, T-splines, LR-splines, etc.). The approach weakly enforces the continuity of the solution at patch interfaces and the error can be adaptively controlled by leveraging the
refineability of the underlying dual spline basis without introducing any additional
degrees of freedom. We also develop weakly continuous geometry as a particular application of isogeometric \Bezier dual mortaring. Weakly continuous geometry is a geometry description where the weak continuity constraints are built into properly modified \Bezier extraction operators. As a result, multi-patch models can be processed in a solver directly without having to employ a mortaring solution strategy. We demonstrate the utility of the approach on several challenging benchmark problems.
\end{abstract}


\end{frontmatter}
\section{Introduction}
Isogeometric Analysis (IGA), first introduced by Hughes et al.~\cite{hughes2005isogeometric}, adopts the Computer-aided design (CAD) basis as the basis for analysis. This unifying paradigm has the potential to eliminate the costly geometry clean-up and mesh generation steps which encumber traditional simulation pipelines and improve simulation accuracy through a higher-order smooth basis~\cite{hughes2005isogeometric, Cottrell:2009rp, COTTRELL20074160, HuEvRe13}. To introduce additional flexibility into the isogeometric approach, weak coupling techniques are often used to sew together models composed of multiple patches. These approaches can accommodate patches with differing parameterizations and trimming~\cite{NME:NME4568,Breitenberger2015401,guo2015nitsche,ruess2014weak,Hesch2012104,NME:NME4918,Coox2017235,Nguyen2014,Seitz2016259,Dornisch2017449,Brivadis2015292}. However, if not done properly, these coupling techniques can negatively impact the accuracy and robustness of the analysis~\cite{NME:NME4568,Breitenberger2015401,Brivadis2015292}.

In this paper, we present a new local dual mortar method for the coupling of
nonconforming higher-order smooth meshes that is based on B\'{e}zier extraction
and projection~\cite{Thomas2015}. Since it is a biproduct of the B\'{e}zier
extraction and projection framework, it can be employed during the creation and
editing of geometry through properly modified extraction operators and is
applicable to any spline space which has a representation in B\'{e}zier form
(i.e., NURBS, T-splines, LR-splines, etc.). For this reason, we call the method
the isogeometric \Bezier dual mortar method. Since the dual basis can be refined
in a fashion which is similar to the corresponding spline basis, the error in
the method can be adaptively controlled without adding \textit{any} additional
degrees of freedom to the linear system. For matched parameterizations, knot insertion can be employed to generate a dual mortaring with optimal approximation properties. All numerical examples show that the proposed method works equally well for arbitrary pairings of the master and slave patches. We also develop weakly continuous geometry as an application of dual mortaring in the context of geometric design. A weakly continuous model is a model in which weak continuity constraints are embedded directly into the geometry description. In this way, multi-patch models can be processed in a solver directly, without having to employ a dual mortaring approach during the construction of the linear system.

\subsection{A review of weak coupling methodologies}
To provide some context and background for the method proposed in this paper, we provide a brief review of existing weak coupling methodologies that have been used in FEA and IGA. There is a vast literature on the subject so only those contributions most closely related to the proposed approach have been included in the overview.

The penalty method~\cite{babuvska1973finite, NME:NME4568, Breitenberger2015401} weakly imposes a coupling constraint by introducing a penalty term into the variational formulation. It is simple to implement and it does not introduce any additional degrees-of-freedom. The drawback is that, to get an accurate result, a problem and mesh dependent penalty parameter must be selected. This parameter, if not properly adjusted during mesh refinement, results in ill-conditioned linear systems~\cite{NME:NME4568}.

The Lagrange multiplier method employs a field of Lagrange multipliers to weakly enforce a coupling constraint. For structural mechanics problems, the field of Lagrange multipliers can be interpreted as the traction forces across an interface. In the context of mesh coupling, this method is also called the mortar method~\cite{Belgacem1999}. The additional Lagrange multiplier field leads to a saddle point variational formulation, which requires that the Lagrange multiplier space satisfy inf-sup stability and ideally have enough approximability to recover optimal convergence rates~\cite{brezzi:13,Belgacem1999}. In the context of IGA, the mortar method was first used to couple multiple non-uniform rational B-splines (NURBS) patches by Dornisch et al.~\cite{PAMM:PAMM201110095}, and then applied in nonlinear elasticity by Hesch and Betsch~\cite{Hesch2012104}. Brivadis et al.~\cite{Brivadis2015292} explored several choices for the Lagrange multiplier space theoretically and numerically.

The Nitsche method~\cite{Nitsche1971}, originally introduced for the weak treatment of Dirichlet boundary conditions, is a method that has a variational structure between the Lagrange multiplier and penalty methods. In this approach, the Lagrange multiplier in the variational formulation is replaced by the normal flux, and an extra penalty-like stabilization term is added to restore the coercivity of the bilinear form. This method has been applied to the coupling of non-conforming meshes in many areas, including IGA~\cite{Nguyen2014, NME:NME4568, guo2015nitsche, ruess2014weak}. Like the penalty method, the stabilization term contains a parameter which must be estimated~\cite{NME:NME4568}.

The approach proposed in~\cite{bernardi1993domain, bernardi1994new} embeds the coupling constraints into the finite element space
directly, thus leading to a positive definite nonconforming variational problem. 
Based on~\cite{bernardi1993domain, bernardi1994new}, Wohlmuth~\cite{wohlmuth2000mortar, wohlmuth2001mortar} then proposed a local dual
Lagrange multiplier space for linear triangular discretizations and called the resulting formulation a dual mortar method. In contrast to a standard Lagrange multiplier method,
in a dual mortar method the Lagrange multipliers can be eliminated easily leading to greater computational efficiency. In
addition, the compact support of the local dual basis preserves
the sparsity of the stiffness matrix. Unfortunately, it is not easy to construct a
local dual basis that possesses a high-order polynomial reproduction
property~\cite{Lamichhane2002,oswald2001polynomial}.

Dornisch et al.~\cite{NME:NME4918} developed a dual mortar method based on a global B-spline
dual basis, and derived a relation matrix which enabled a condensation of the
Lagrange multiplier degrees-of-freedom. A similar relation matrix is derived by
Coox et al.~\cite{Coox2017235} by inserting virtual knots on either side of an
interface. This method is mathematically identical to the global dual method
in~\cite{NME:NME4918} but is more efficient. However, it is limited to the case
where the neighbouring patches have the same degree and parameterization along
the interface, which is a very restrictive requirement. Seitz et al.~\cite{Seitz2016259}
proposed a local dual mortar method based on a NURBS basis for both patch
coupling and contact mechanics. In this case, the local dual basis
does not satisfy the polynomial reproduction property, so only reduced
convergence rates are obtained. Other types of local dual basis functions, such
as, the explicit de Boor-Fix dual basis~\cite{DEBOOR197319, de1975local,
  schumaker2007spline} and the approximate dual basis~\cite{CHUI2004141} are
explored in~\cite{Dornisch2017449}. The de Boor-Fix dual basis functions have the same support as the B-spline basis functions. However, the polynomial reproduction property does not hold, leading to significantly deteriorated convergence rates. The approximate dual basis fulfils the polynomial reproduction property but not biorthogonality. Therefore, the fully populated inverse matrix of the original mortar matrix must be approximated by a diagonal matrix to maintain the locality. This implies that the coupling constraints are not imposed exactly. Even though several numerical examples show that the approximate dual mortar method achieves convergence rates which are comparable to the global dual mortar method, a mathematical analysis of the effects of the approximation is still missing.

The outline of this paper is as follows. In Section \ref{sec:GeomFund}, we
briefly review fundamental concepts for splines and dual
  bases which are needed throughout the paper. Section~\ref{sec:problem} describes multi-patch domain decomposition and the model problem we
will use to define our method. Isogeometric \Bezier dual mortaring is then
described in Section~\ref{sec:method}. We then define weakly continuous geometry and its relationship to \Bezier dual mortaring in Section~\ref{sec:weak_geom}. Several challenging benchmark problems are solved in Section~\ref{sec:results} to illustrate the properties of the method. We then draw conclusions in Section~\ref{sec:conclude}.

\section{Spline fundamentals and dual bases}
\label{sec:GeomFund}
\subsection{B\'{e}zier, B-spline, and NURBS fundamentals}
The $i$th Bernstein polynomial of degree $p$ on $[\xi_1,\xi_2]$ can be defined as
\begin{align}
  B_{i,p}(\xi) = {{p}\choose{i-1}}\left (\frac{\xi_2-\xi}{\xi_2-\xi_1} \right )^{p-i+1}\left ( \frac{\xi-\xi_1}{\xi_2 - \xi_1} \right )^{i-1}, 
\end{align}
where ${{p}\choose{i-1}} = \frac{p!}{(i-1)!(p-i+1)!}$ is the binomial
coefficient. 

The set of Bernstein polynomials
$\mathbf{B}(\xi)=\{B_{i,p}(\xi)\}_{i=1}^{p+1}$ forms a basis for the space of
polynomials of degree $p$.  The Bernstein polynomials
$\tilde{\mathbf{B}}(\tilde{\xi})$ defined on $[\tilde{\xi_1},  \tilde{\xi_2}]$, can be related to the Bernstein basis $\mathbf{B}(\xi)$ defined on $[\xi_1, \xi_2]$ through the relation
\begin{align}\label{eq:splitBezier}
\tilde{\mathbf{B}}(\tilde{\xi})= (\mathbf{M})^{-\text{T}} \mathbf{B}(\xi)
\end{align}
with $\mathbf{M}$ the transformation matrix. A formula for the inverse of the transformation matrix $\mathbf{M}$ can be found in~\cite{farouki1990numerical} and is written as
\begin{align}\label{eq:Mmatrix}
  (M)_{jk}^{-1} = \sum_{l = \text{max}(1, j+k-p-1)}^{\text{min}(j,k)}B_{l,j-1}(\xi_2)B_{k - l  + 1, p - j + 1}(\xi_1), \quad \quad  1 \leq j, k \leq p+1.
\end{align}
A degree $p$ \Bezier curve in $\mathbb{R}^{d}$ can be written as 
 \begin{align}
   \mathbf{x}(\xi) = \sum_{i=1}^{p+1}\mathbf{P}_i \,\,B_{i,p}(\xi), \quad \quad  \; \xi \in [ \xi_1,\xi_2 ]
 \end{align}
 where $\mathbf{P}_i$ is called a control point.
A univariate B-spline basis is defined by a knot vector ${\varXi}= \left\{ \xi_1,\xi_2,\ldots,\xi_{n+p+1}\right\}$, which consists of a non-decreasing sequence of real numbers, $\xi_i \leq \xi_{i+1},i=1,\ldots,n+p+1$, where $p$ is the degree of the B-spline basis functions and $n$ is the number of basis functions. The $i$th B-spline basis function of degree $p$, denoted by $N_{i,p}(\xi)$, can be recursively defined by
\begin{gather*}
N_{i,0}(\xi) = 
 \begin{cases} 
      1 ,& \text{if} \quad \xi_i \leq \xi < \xi_{i+1} \\
      0 , &  \text{otherwise}
   \end{cases}\\[2ex]
  N_{i,p}(\xi) = \frac{\xi-\xi_i}{\xi_{i+p} - \xi_i}N_{i,p-1}(\xi) + \frac{\xi_{i+p+1} - \xi}{\xi_{i+p+1} - \xi_{i+1}}N_{i+1,p-1}(\xi).
\end{gather*}
A B-spline curve can be viewed as the smooth composition of multiple \Bezier curves. A B-spline curve of degree $p$ can be written as
\begin{align}
\mathbf{x}(\xi) = \sum_{i =1}^n \mathbf{P}_iN_{i,p}(\xi),
\quad  \quad \; \xi \in [ \xi_1,\xi_{n+p+1}] 
\end{align}
and a $p$th-degree NURBS curve can be written as
\begin{align}
\mathbf{x}(\xi) = \sum_{i = 1}^n\mathbf{P}_i w_iR_{i,p}(\xi),
\quad \quad \; \xi \in [ \xi_1,\xi_{n+p+1}]
\end{align}
where the NURBS basis function $R_{i,p}$ is defined by
\begin{align}
R_{i,p}(\xi) = \frac{N_{i,p}(\xi)}{W(\xi)}
\end{align}
where $N_{i,p}(\xi)$ is the $i$th $p$-degree B-spline basis functions,
\begin{align}
  \label{eq:weight}
  W(\xi) = \sum_{i=1}^{n}{w_iN_{i,p}(\xi)}
\end{align}
is the weighting function, and $w_i$ is the weight corresponding to $N_{i,p}(\xi)$. Since a NURBS curve is a rational polynomial it can be used to exactly represent conic sections. Higher dimensional analogs to these concepts can be created using tensor products or more advanced construction schemes like T-splines or hierarchical B-splines.
\subsection{\Bezier extraction}
The \Bezier extraction process~\cite{NME:NME2968,NME:NME3167} generates a linear operator, called the extraction operator, that maps a Bernstein basis onto a B-spline basis. In the context of one-dimensional B-splines, the extraction operator encodes the result of repeated knot insertion~\cite{PiegTil97} such that the multiplicity of all interior knots of a knot vector is $p+1$. At the element level, the resulting linear transformation, $\mathbf{C}^e$, is called the element extraction operator. This element-level operator is used to map a Bernstein basis $\mathbf{B}$ defined over an element $e$ onto a B-spline basis restricted to that same element. In other words, $\mathbf{N}^e = \mathbf{C}^e\mathbf{B}$. See~\cite{NME:NME2968,NME:NME3167} for additional details.
\subsection{Dual bases}
Suppose $\mathcal{B}_{p}$ is a $(p+1)$-dimensional linear space generated by a set of linearly independent functions $\left\{b_i\right\}_{i=1}^{p+1}$ of maximal degree $p$. Given an inner product $(\cdot, \cdot)\colon\mathcal{B}_p \times \mathcal{B}_p \mapsto \mathbb{R} $, the functions from the set
\begin{align}
  \bm{\lambda}_p := \{ \lambda_i\}_{i=1}^{p+1}
\end{align}
satisfying the following conditions
\begin{align}\label{dualconditions}
  \begin{cases} 
      \text{span}\, \bm{\lambda}_p  = \mathcal{B}_p, \\[1ex]
      (b_i, \lambda_i) = \delta_{ij}, \quad 1 \leq i, j \leq p+1
   \end{cases}
\end{align}
form the so-called dual basis corresponding to the basis
$\left\{b_i\right\}_{i=1}^{p+1}$ with respect to the inner product $(\cdot,
\cdot)$. The first condition in (\ref{dualconditions}) is called the
reproduction property of order $p$ and the second property is called the
biorthogonality property. With local dual mortar method, $(p-1)$-order reproduction
property of the Lagrange multiplier space is required to
guarantee the optimality of the finite element space of order $p$~\cite{oswald2001polynomial}. 

\section{Problem description}
\label{sec:problem}

\subsection{Domain decomposition}
Let $\Omega$ be a bounded domain decomposed into $K$ non-overlapping subdomains $\Omega^k$, i.e.,
\[
\overline{\Omega} = \bigcup_{k=1}^K \overline{\Omega}^k, \text{ and } \Omega^i \cap \Omega^j = \emptyset, i \neq j.
\]
We define the interfaces as the interior of the intersections of the boundaries,
i.e., $\Gamma^\ell = \partial \Omega^i \cap \partial \Omega^j$. On each $\Omega^k$
the solution space $\mathcal{S}^k$ is defined as 
\[\mathcal{S}^k=\{ \mathbf{u}^k \in H^1(\Omega^k),  \mathbf{u}^{k}{ |_{\partial \Omega \cap \partial \Omega^k}}=\mathbf{u}_0\}
\]
where $H^1(\Omega^k)$ are the standard Sobolev spaces and $\mathbf{u}_0$ are the
Dirichlet boundary conditions. The corresponding weighting function spaces, $\mathcal{V}^k$, are similarly defined with homogeneous boundary conditions on $\partial \Omega \cap \partial \Omega^k$.
The displacement solution space on $\Omega$ is then the broken Sobolev space
$\mathcal{S}$ defined as $\mathcal{S}= \prod_{k=1}^K \mathcal{S}^k$, along with
continuity conditions defined along the interfaces. To simplify the exposition
of the proposed mortaring technique, we employ a two-patch geometry with one
interface, i.e., $K=2$. The interface is denoted as $\Gamma = \partial
{\Omega}^{m} \cap \partial
{\Omega}^{s} $, where the superscripts $m$ and $s$ are used to denote the
master and slave patches, respectively.

\subsection{A linear elastic model problem}
\begin{figure}[h]
  \centering
\includegraphics[width=0.5\textwidth]{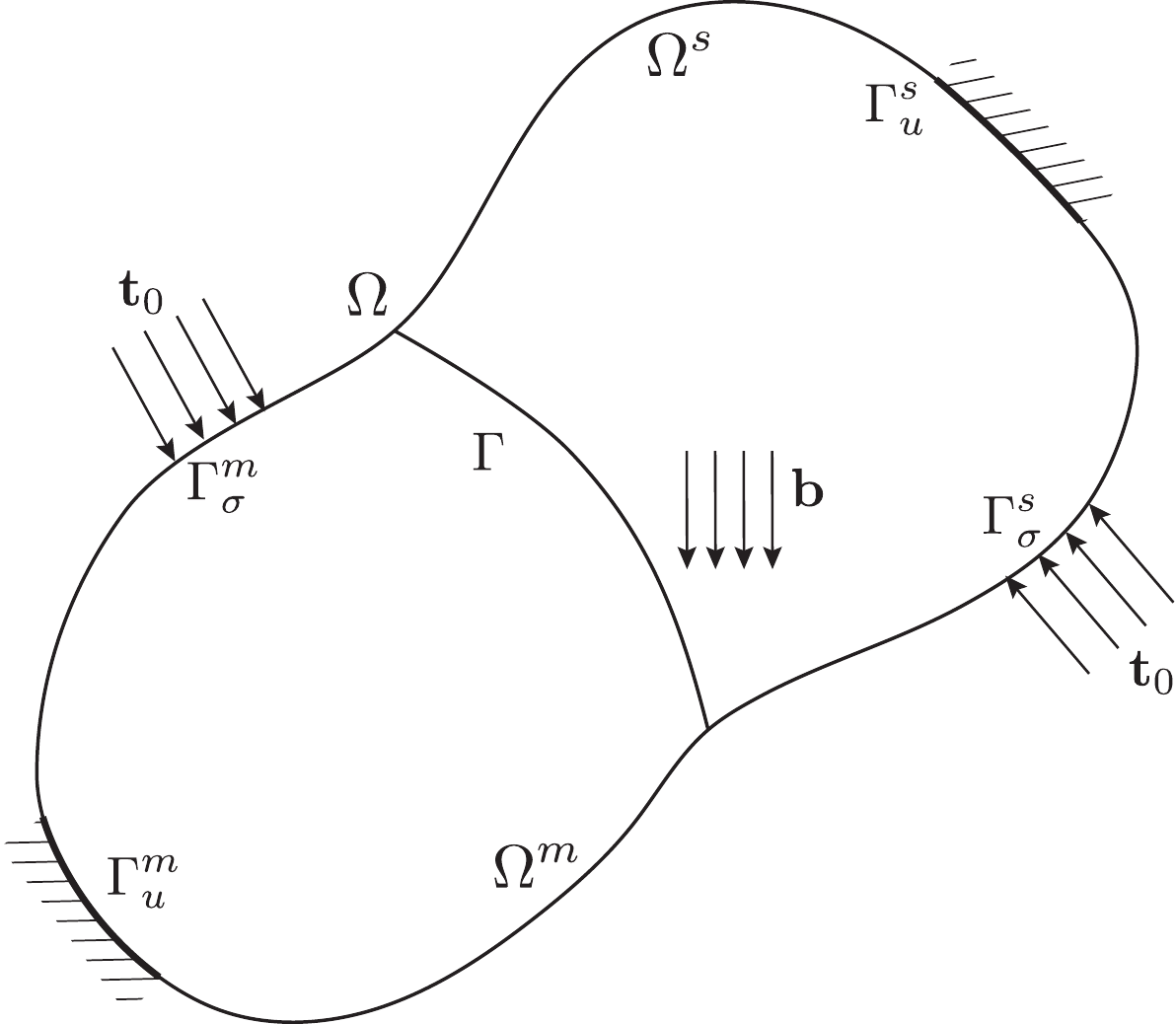}
\caption{A schematic of the linear elastic model problem.}
\label{fig:elasticitydescription}
\end{figure}

To ground our approach in a practical example, we consider the following linear elasticity problem:
\begin{subequations}\label{eq:strongform}
\begin{align}
 \text{div}~\bm{\sigma} + \mathbf{b} &= 0 \quad  \text{in} \: \Omega = \Omega^m \cap \Omega^s  \\
  \mathbf{u} &= \mathbf{u}_0 \quad  \text{on} \: \Gamma_u = \Gamma_u^m \cap \Gamma_u^s  \\
  \bm{\sigma}\cdot \mathbf{n} &= \mathbf{t}_0 \quad  \text{on} \: \Gamma_{\sigma} = \Gamma_\sigma^m \cap  \Gamma_\sigma^s\\
 \mathbf{u}^m &=\mathbf{u}^s \quad  \text{on} \: \Gamma \label{eq:strongconstraints}
\end{align}
\end{subequations}
where $\Gamma_u \cap \Gamma_{\sigma}=\emptyset $, $\Gamma_u \cap \Gamma=\emptyset $, $\Gamma_{\sigma} \cap \Gamma=\emptyset $, $\bm{\sigma}$ is the stress tensor, $\mathbf{b}$ is the body force, $\mathbf{u}_0$ and $\mathbf{t}_0$ are the prescribed Dirichlet and Neumann boundary conditions applied on $\Gamma_u$ and $\Gamma_{\sigma}$, respectively and  $\mathbf{n}$ is the unit outward normal vector on $\partial \Omega$, see Figure~\ref{fig:elasticitydescription}. The kinematic coupling condition $\mathbf{u}^m = \mathbf{u}^s$ is introduced along the interface $\Gamma$ where $\mathbf{u}^m$ and $\mathbf{u}^s$ are the master and slave interface displacements, respectively. 

 The total potential energy $\Pi$ of the system $\Omega$ is
 \begin{align}\label{eq:total_potential}
 \Pi(\mathbf{u}) = \Pi^m(\mathbf{u}^m) +  \Pi^s(\mathbf{u}^s) + \int_{\hat{\Gamma}} \mathbf{\Phi}\cdot( \mathbf{u}^m - \mathbf{u}^s) \, \rm{d}s
 \end{align}
 where $\Pi^m$ and $\Pi^s$ are the potential energy on $\Omega^m$ and $\Omega^s$, respectively, and $\mathbf{\Phi}$ is a Lagrange multiplier weakly enforcing the continuity constraint along the interface.  Invoking the stationarity of $\Pi$ with respect to $\mathbf{u}^m$, $\mathbf{u}^s$ and $\mathbf{\Phi}$, we obtain the weak formulation of (\ref{eq:strongform}) that reads as: find $\mathbf{u}^m \in \mathcal{S}^m$, $\mathbf{u}^s \in \mathcal{S}^s$ and  $\mathbf{\Phi} \in \mathcal{S}^{\ell}$ such that for all variations $\delta \mathbf{u}^m \in \mathcal{V}^m$, $\delta \mathbf{u}^s \in \mathcal{V}^s$ and  $\delta \mathbf{\Phi} \in \mathcal{V}^{\ell}$

\begin{subequations}\label{eq:weakform}
\begin{align}
 \delta \Pi^*(\mathbf{u}, \delta \mathbf{u}^m) &=\delta \Pi^m(\mathbf{u}^m,\delta \mathbf{u}^m) + \int_{\hat\Gamma} \mathbf{\Phi}\cdot \delta \mathbf{u}^m \, \rm{d}\Gamma = 0 \\
 \delta \Pi^*(\mathbf{u}, \delta \mathbf{u}^s) &=\delta \Pi^s(\mathbf{u}^s,\delta \mathbf{u}^s) - \int_{\hat\Gamma} \mathbf{\Phi}\cdot \delta \mathbf{u}^s \, \rm{d}\Gamma = 0 \\
 \delta \Pi^*(\mathbf{u}, \delta \mathbf{\Phi}) &= \int_{\hat\Gamma} \delta \mathbf{\Phi}\cdot (\mathbf{u}^m - \mathbf{u}^s) \, \rm{d}s = 0 \label{eq:weakformconstraint}
\end{align}
\end{subequations}
where $\mathcal{S}^m$, $\mathcal{S}^s$ and $\mathcal{S}^{\ell}$ are the
displacement solution approximation spaces on $\Omega^m$ and $\Omega^s$ and the
Lagrange multiplier space, respectively, and $\mathcal{V}^m$, $\mathcal{V}^s$
and $\mathcal{V}^{\ell}$ are the corresponding weighting function spaces. 
Note that in (\ref{eq:total_potential}) we define the interface energy on
the parametric domain of the slave interface, denoted by $\hat{\Gamma}$, which results in the interface
continuity condition (\ref{eq:weakformconstraint}). As will be shown subsequently, this will allow us to define a dual basis which is independent of geometry, an important simplification which improves the efficiency of the approach.
\section{Isogeometric \Bezier dual mortaring}
\label{sec:method}
We will choose the Lagrange multiplier spaces to be those spanned by a dual
spline basis defined over $\hat\Gamma$ which emanate from the B\'ezier extraction and projection framework. When a local dual basis is chosen for the Lagrange multiplier spaces the method is often called a dual mortar method.
A weighted dual basis for each element domain $\hat{\Gamma}^e$ is defined as
 \begin{align}\label{eq:elemdual}
  \bar{\mathbf{N}}^{e} &= \operatorname{diag}(\boldsymbol{\omega}^e) (\mathbf{R}^e)^{T} (\mathbf{G}_{B,B}^e)^{-1} \mathbf{B}^{e,s} \\
                       &= \mathbf{D}^e \mathbf{B}^{e,s}
 \end{align}
where $\mathbf{B}^{e,s}$ is the set of Bernstein polynomials defined on the
$e$th slave interface element and 
\begin{align}
\mathbf{G}_{B,B}^e &= \left [ \int_{\hat{\Gamma}^e}  B_i^{e,s}(\xi^s) B_j^{e,s}(\xi^s) \, ds \right ]
\end{align}
is the Gramian matrix for the Bernstein basis~\cite{Thomas2015,juttler1998},
$\mathbf{R}^e$ is the element reconstruction operator~\cite{Thomas2015}  and
$\mathbf{D}^e$ is called a dual element extraction operator. Note that $\mathbf{R}^e$ is restricted to the element boundary $\hat{\Gamma}^e$. We use the standard \Bezier projection weighting, i.e.,
\begin{align}
\omega_i^e = \frac{\int_{\hat{\Gamma}^{e}} N_i^{e,s} \, ds}{ \int_{\hat{\Gamma}^{I}} N^s_{I(i,e)} \, ds}
\end{align} 
where $\hat{\Gamma}^I$ is the domain of support for the interface basis function $N^s_I$ and $I(i,e)$ is a standard mapping from element nodal indexing to a global index $I$. While other weightings could be used this weighting has been shown to give particularly accurate results~\cite{Thomas2015}. We can easily show that the proposed dual basis satisfies the biorthogonality condition (\ref{dualconditions}) by noting that
\begin{align}
  \int_{\hat{\Gamma}^e } \bar{\mathbf{N}}^{e} (\mathbf{N}^{e,s})^T ds &= \operatorname{diag}(\boldsymbol{\omega}^e)
\end{align}
and
\begin{align}
\assembly_{e} \left[\int_{\hat{\Gamma}^e } \bar{\mathbf{N}}^{e} (\mathbf{N}^{e,s})^T \, ds \right] &= \mathbf{I}.
\end{align}
where $\assembly$ is the usual finite element assembly operator. In other words,
\begin{align}
 \int_{\hat{\Gamma}} \bar{N}_I N^s_J ds = \delta_{IJ}
 \end{align}
 as desired.

 Note that even though this dual basis does not possess a higher-order polynomial reproduction property, optimal higher-order rates can be easily recovered through a simple refinement step as described in Section~\ref{sec:refine}.

\subsection{Rational dual basis functions}
If rational basis functions are used, we define the dual basis as
\begin{align}
  \bar{R}_I &= W \bar{N}_i
\end{align}
where $W$ is the rational weight given in (\ref{eq:weight}). Now
\begin{align}
  \int_{\hat{\Gamma}} \bar{R}_I R^s_J \, ds &= \int_{\hat{\Gamma}} \bar{N}_I N^s_J \, ds = \delta_{IJ}.
\end{align}
\subsection{Discretization}
\begin{figure}[h]
  \centering
  \begin{overpic}[scale=0.8]{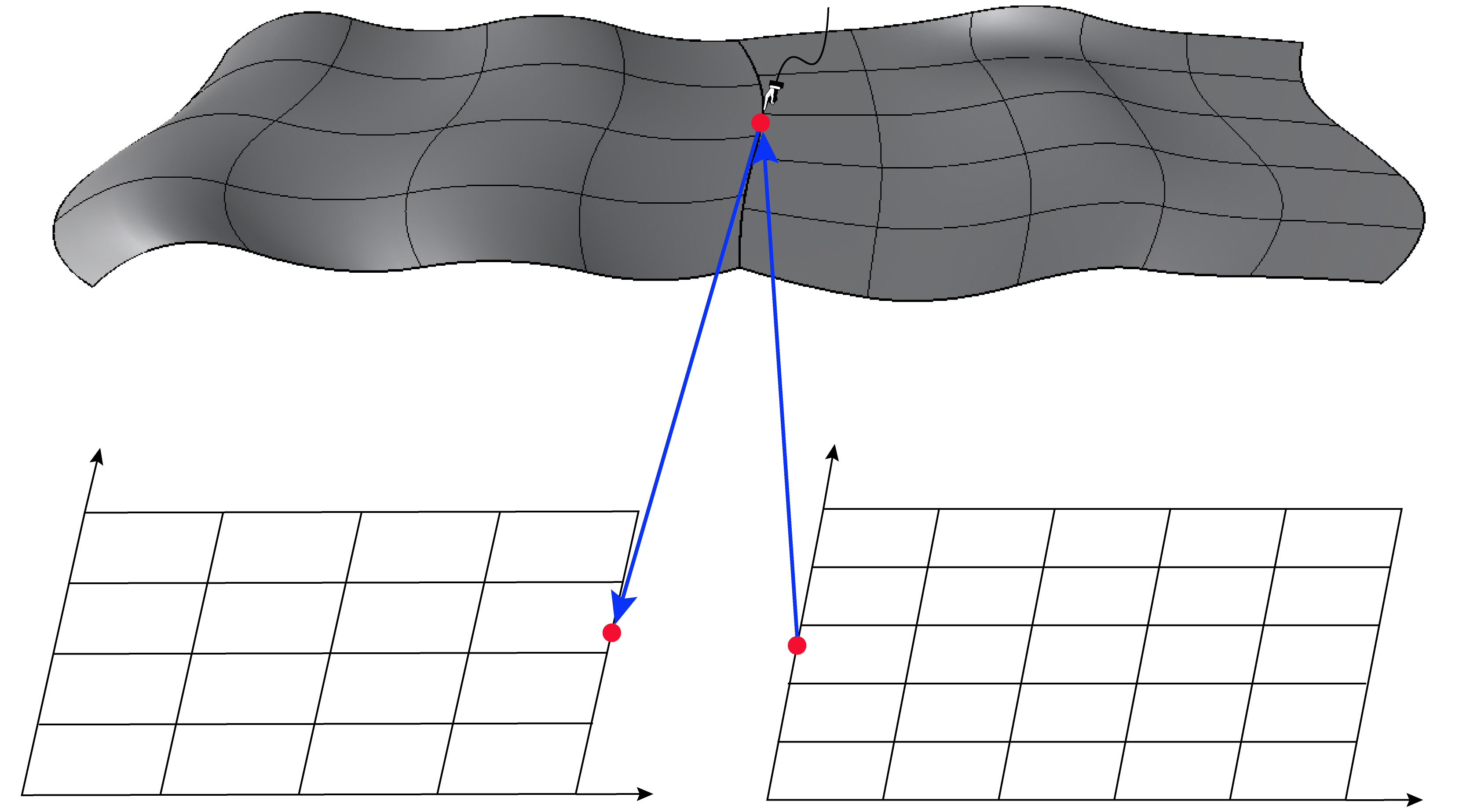}
    \put(3, 2){$\hat{\Omega}^m$}
    \put(54, 2){$\hat{\Omega}^s$}
    \put(9, 40){$\Omega^m$}
    \put(90, 41){$\Omega^s$}
    \put(43, 3){$\xi^m_1$}
    \put(2, 23){$\xi^m_2$}
    \put(95, 3){$\xi^s_1$}
    \put(58, 23){$\xi^s_2$}
    \put(56, 10){$\xi^s$}
    \put(33, 8){$\xi^m = \varphi(\xi^s)$}
    \put(54, 30){$\mathbf{x}^s$}
    \put(35, 30){$(\mathbf{x}^m)^{-1}$}
    \put(52, 56){$\mathbf{x}^m(\xi^m) = \mathbf{x}^s(\xi^s)$}
  \end{overpic}
  \caption{Slave and master geometric mappings, $\mathbf{x}^s$ and $\mathbf{x}^m$.}	
\label{fig:geometrymap}
\end{figure}
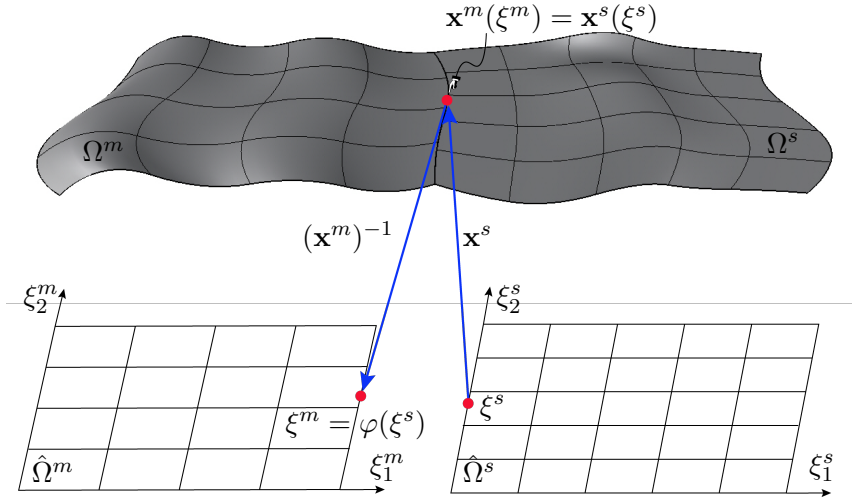

Over the slave interface we introduce the discretizations
\begin{align}
\mathbf{u}^{m} &=  \sum_{I} N^m_I \left( \varphi(\xi^s) \right) \, \mathbf{d}_{I}^{m} \\
\mathbf{u}^s &= \sum_{I} N^s_I (\xi^s) \, \mathbf{d}_I^{s} \\
\delta\mathbf{\Phi} &= \sum_{I} \bar{N}_I(\xi^s) \, \delta \bm{\Phi}_{I}
\end{align}
where $\xi^s \in \hat{\Gamma}^s$ is a parametric position on the slave interface and $\varphi : \hat{\Gamma}^s \rightarrow \hat{\Gamma}^m$ is a compositional mapping defined to be
\begin{align}
  \varphi &= (\mathbf{x}^m)^{-1} \circ \mathbf{x}^s
\end{align}
where $\mathbf{x}^s : \hat{\Omega}^s \rightarrow \Omega^s$ and $\mathbf{x}^m  :
\hat{\Omega}^m \rightarrow \Omega^m$ are the slave and master geometric
mappings, respectively, as shown in Figure~\ref{fig:geometrymap}. Note that we say the master and slave parameterizations are matched if the mapping $\varphi$ is
linear, otherwise, we say the master and slave parameterizations are mismatched. In the mismatched case $\varphi$ can be computed using the Newton-Raphson algorithm.

Discretizing (\ref{eq:weakformconstraint}) and leveraging the biorthogonality property of the dual basis results in
\begin{align}
  \label{discretedispcons}
 \mathbf{d}^{s} &= \left[  \int_{\hat{\Gamma}} \bar{N}_I (\xi^s)N_{J}^{m}\left( \varphi(\xi^s) \right)  ds \right] \mathbf{d}^{m} \nonumber \\
 &= \mathbf{G}_{\bar{N}, N^m} \mathbf{d}^{m}.
\end{align}

The matrix form of (\ref{eq:weakform}) can be written as
\begin{align} \label{eq:lag_sys_eqs}
  \begin{bmatrix}
    \mathbf{K}^m          &             0            & \mathbf{K}^{\ell m} \\
    0               &  \mathbf{K}^s            &  -\mathbf{K}^{\ell s}\\
    ( \mathbf{K}^{\ell m} )^{\text{T}} &  - ( \mathbf{K}^{\ell s} )^{\text{T}} &                       0\\ 
  \end{bmatrix} 
  \begin{bmatrix}
    \mathbf{d}^m \\
    \mathbf{d}^s \\
    \mathbf{d}^\ell
  \end{bmatrix} =
  \begin{bmatrix}
    \mathbf{f}^m \\
    \mathbf{f}^s \\
    0
  \end{bmatrix} 
\end{align} 
where $\mathbf{K}^m$ and $\mathbf{K}^s$ are standard patch-level
stiffness matrices, $\mathbf{f}^m$ and $\mathbf{f}^s$ are the
corresponding force vectors, and $\mathbf{K}^{\ell m}$ and
$\mathbf{K}^{\ell s}$ are stiffness matrices with all entries equal to zero except for those related to the $I$th Lagrange multiplier basis $\bar{N}_I$ and the $J$th master and slave interface basis $N_J^m$ and $N_J^s$, respectively. In other words, the components of $\mathbf{K}^{\ell m}$ and $\mathbf{K}^{\ell s}$ can be written as
\begin{align}
  K^{\ell m}_{\bar{N}_I, N^m_J} = \int_{\hat{\Gamma}} \bar{N}_I N_{J}^{m} ds
\end{align}
and
\begin{align}
  K^{\ell s}_{\bar{N}_I, N^s_J}=  \int_{\hat{\Gamma}} \bar{N}_I N_{J}^{s} ds.
\end{align}
The displacement vectors $\mathbf{d}^m$ and $\mathbf{d}^s$ can be split such that
\begin{align}
  \mathbf{d}^m =
  \begin{bmatrix}
    \mathbf{d}^m_d \\
    \mathbf{d}^m_c \\
  \end{bmatrix} \quad \quad 
  \mathbf{d}^s =
  \begin{bmatrix}
    \mathbf{d}^s_d \\
    \mathbf{d}^s_c \\
  \end{bmatrix}
\end{align}
where the subscript $d$ indicates the \textit{distinct}
degrees-of-freedom internal to each patch, and $c$ indicates the
degrees-of-freedom along the interface.
Choosing a local dual space as the Lagrange multiplier space allows us to
condense out the Lagrange multiplier coefficients 
$\mathbf{d}^\ell$ and the slave patch degrees-of-freedom $\mathbf{d}_c^s$
through (\ref{discretedispcons}) while preserving the sparsity of the global stiffness matrix. Condensing (\ref{eq:lag_sys_eqs}) results in
\begin{align} \label{eq:dua_sys_eqs}
  \begin{bmatrix}
    \mathbf{K}^m_{dd} &  \mathbf{K}^m_{dc}                                                                        &  0  \\
    \mathbf{K}^m_{cd} &  \mathbf{K}^m_{cc}  + ( \mathbf{G}_{\bar{N},N^m} )^\text{T} \mathbf{K}^s_{cc} \mathbf{G}_{\bar{N},N^m} & ( \mathbf{G}_{\bar{N},N^m})^\text{T} \mathbf{K}^s_{cd}\\
    0                & \mathbf{K}^s_{dc}  \mathbf{G}_{\bar{N},N^m}                                                       &                       \mathbf{K}^s_{dd} 
  \end{bmatrix} 
  \begin{bmatrix}
    \mathbf{d}^m_d \\
    \mathbf{d}^m_c \\
    \mathbf{d}^s_d
  \end{bmatrix} =
  \begin{bmatrix}
    \mathbf{f}^m_d \\
    \mathbf{f}^m_c \\
    \mathbf{f}^s_d \\
  \end{bmatrix}.
\end{align}

Note that if the interface energy term in (\ref{eq:total_potential}) is defined
on the physical domain instead of the parametric domain the dual basis must be
defined as
 \begin{align}
  \bar{\mathbf{N}}^{e} &= \frac{1}{|\mathbf{J}|}\operatorname{diag}(\boldsymbol{\omega}^e) (\mathbf{R}^e)^{T} (\mathbf{G}_{B,B}^e)^{-1} \mathbf{B}^{e,s} \\
                       &= \frac{1}{|\mathbf{J}|}\mathbf{D}^e \mathbf{B}^{e,s}
 \end{align}
where $\mathbf{J}$ is the Jacobian of the geometric mapping
$\mathbf{x}^s$.

\subsection{Refinement of the dual basis}
\label{sec:refine}
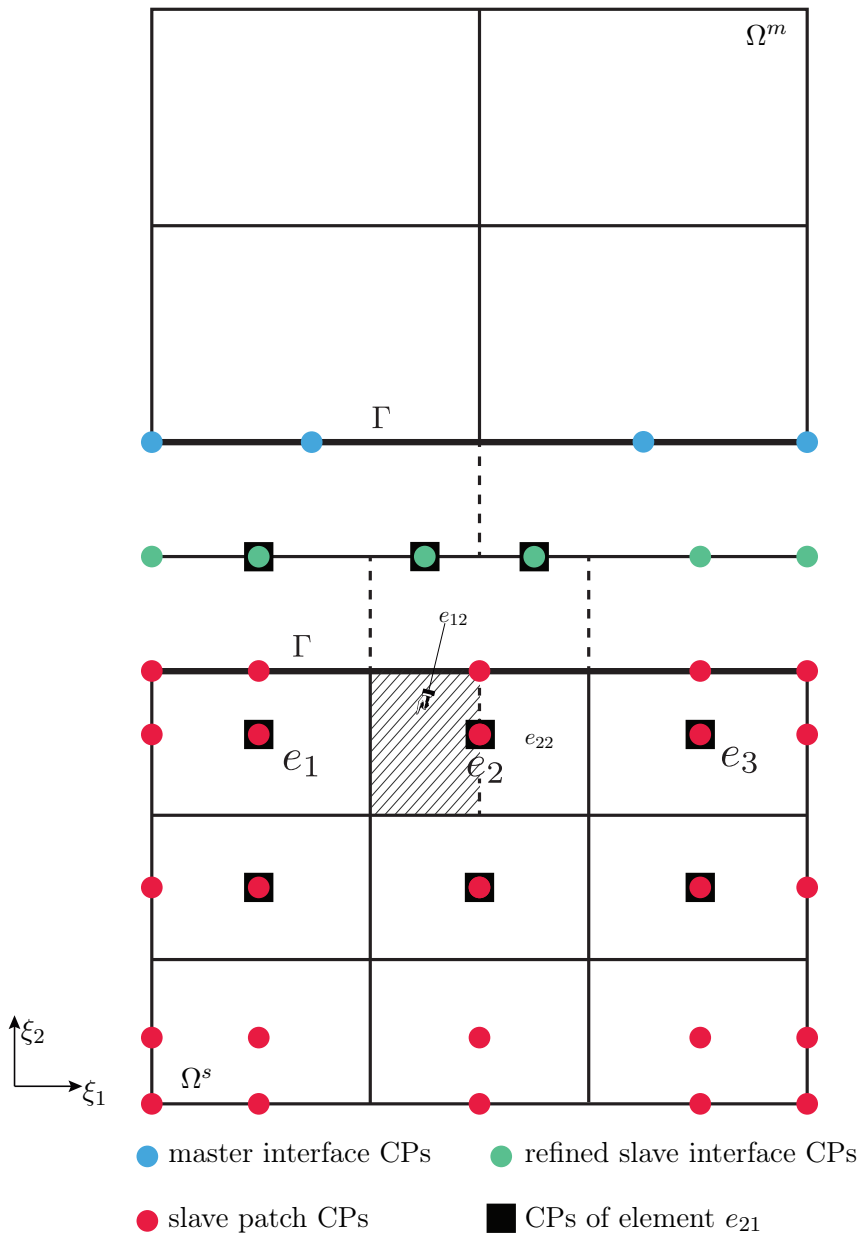
\begin{figure}[h]
  \centering
  \begin{overpic}[scale=1.2]{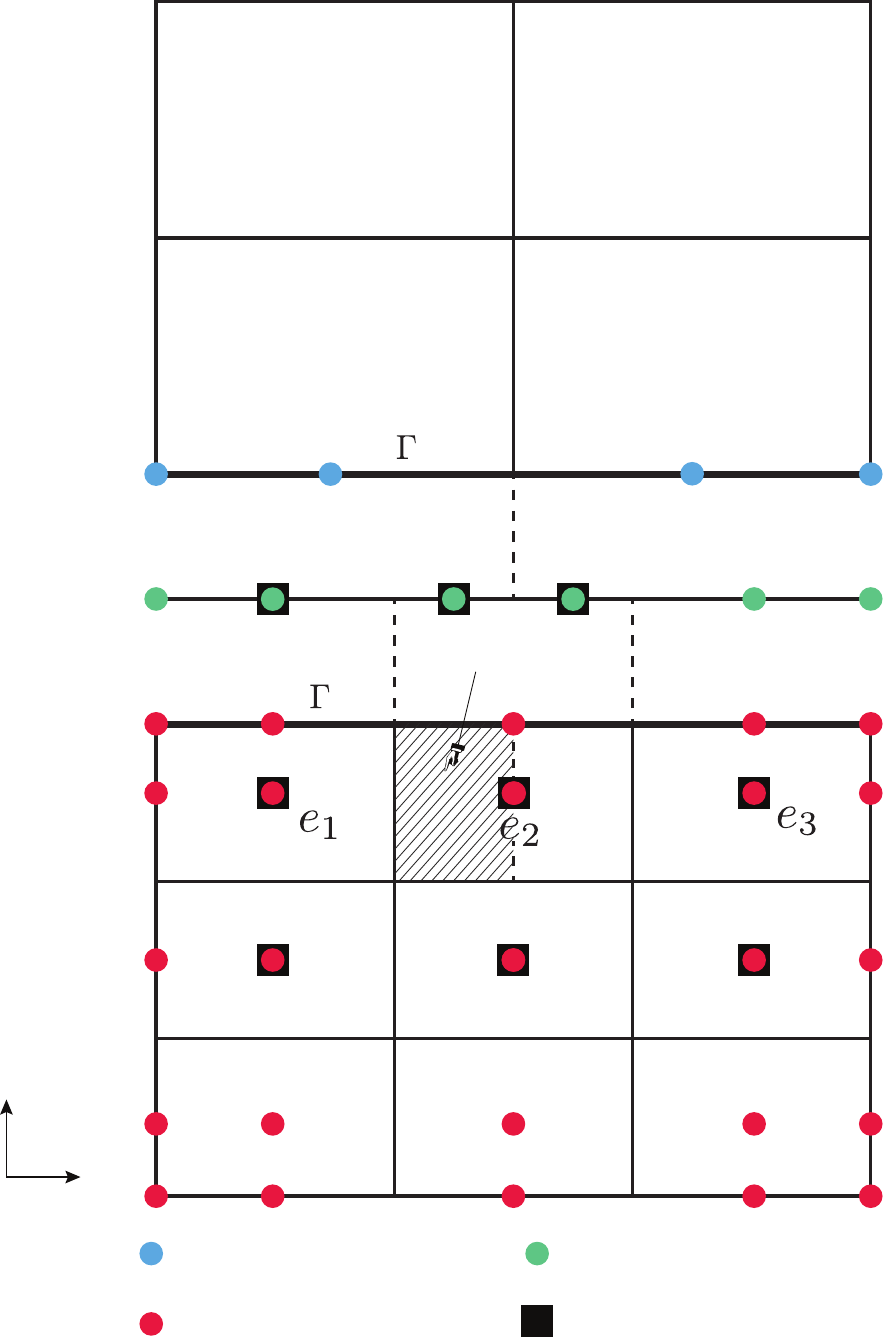}
    \put(14, 12){$\Omega^s$}
    \put(60, 97){$\Omega^m$}
   \put(35, 50){\scalebox{0.8}{$e_{12}$}}
   \put(42, 40){\scalebox{0.8}{$e_{22}$}}
   \put(13, 6){master interface CPs}
   \put(13, 0.7){slave patch CPs}
   \put(42, 6){refined slave interface CPs}
   \put(42, 0.7){CPs of element $e_{21}$}
    \put(6, 11){$\xi_1$}
    \put(1, 16){$\xi_2$}
  \end{overpic} 
  \caption{Refinement of a slave interface and corresponding control points (CPs).}	
\label{fig:schematics_of_spliting_element}
\end{figure}

If the master and slave parameterizations are matched, the underlying basis have
the same degrees, and the knots along the master interface are contained in the slave interface the interface constraint
(\ref{eq:weakformconstraint}) is imposed exactly. In this case, $\mathbf{G}_{\bar{N}, N^m}$ is then a standard spline refinement operator. In any case, the approximation can be improved without adding additional degrees-of-freedom to the global system by refining the slave interface and dual basis. We highlight that if the slave interface is refined, quadrature error accumulates if the new lines of reduced continuity in the slave interface are not accounted for in the element domains $\hat{\Omega}^{s,e}$ which touch the slave interface.
For example, in Figure~\ref{fig:schematics_of_spliting_element}, two quadratic,
linearly parameterized B-spline patches meet at a common interface. A
refinement is performed in which all knots in the master interface which are not
already present in the slave interface are added, i.e., the knot $\frac{1}{2}$
is inserted into the slave interface. To properly account for the new line of
reduced continuity in the slave interface element $e_2$, it is subdivided into
two elements $e_{21}$ and $e_{22}$ and quadrature is performed on both
subelements. The nodes whose basis functions are supported by element $e_{21}$ are depicted in Figure~\ref{fig:schematics_of_spliting_element}.  Note that this element subdivision is only for quadrature and can be performed at the \Bezier element level. The subdivision \textit{does not} add any additional degrees of freedom to the slave patch.

\section{Weakly continuous geometry}
\label{sec:weak_geom}
Since the weak continuity constraint is defined on the parametric domain of
the slave interface \Bezier dual
mortaring can also be viewed as an isogeometric design methodology for building multi-patch
geometry where the weak continuity constraint is built into the space spanned by
the geometric basis. In this case, weak geometric compatibility is preserved for
any choice of control points and the dual mortaring no longer needs to be incorporated into the finite element assembly algorithm. To build the weak continuity constraint into the element extraction operators we start by noticing that
\begin{align}\label{eq:global_basis_relation}
  \mathbf{N}^m &= (\mathbf{G}_{\bar{N}, N^m})^\text{T} \mathbf{N}^s
\end{align}
which can be localized to each element on the interface
\begin{align}\label{eq:elem_basis_relation}
  \mathbf{N}^{m,e} &= (\mathbf{G}^e_{\bar{N}, N^m})^\text{T} \mathbf{R}^e \mathbf{B}^{e,s} \\
  &= \tilde{\mathbf{R}}^e \mathbf{B}^{e,s}
\end{align}
where $\tilde{\mathbf{R}}^e$ is called a weakly continuous element extraction operator. Since
\begin{align}\label{eq:weakly_elem_extrac_op}
\tilde{\mathbf{R}}^e &= (\mathbf{G}^e_{\bar{N}, N^m})^\text{T} \mathbf{R}^e
\end{align}
it is clear that each row of $\tilde{\mathbf{R}}^e$ (which corresponds to a
master basis function) is formed by taking a linear combination of rows in
$\mathbf{R}^e$ (which correspond to slave basis functions) where the weighting
in the linear combination comes from the columns of $\mathbf{G}^e_{\bar{N}, N^m}$. Note that while only weak $C^0$ continuity is considered in this paper, this same framework can be utilized to build other types of constraints into a geometric representation.

\begin{figure}[h]
  \centering
  \begin{subfigure}{0.47\textwidth}
    \begin{overpic}[scale=0.4]{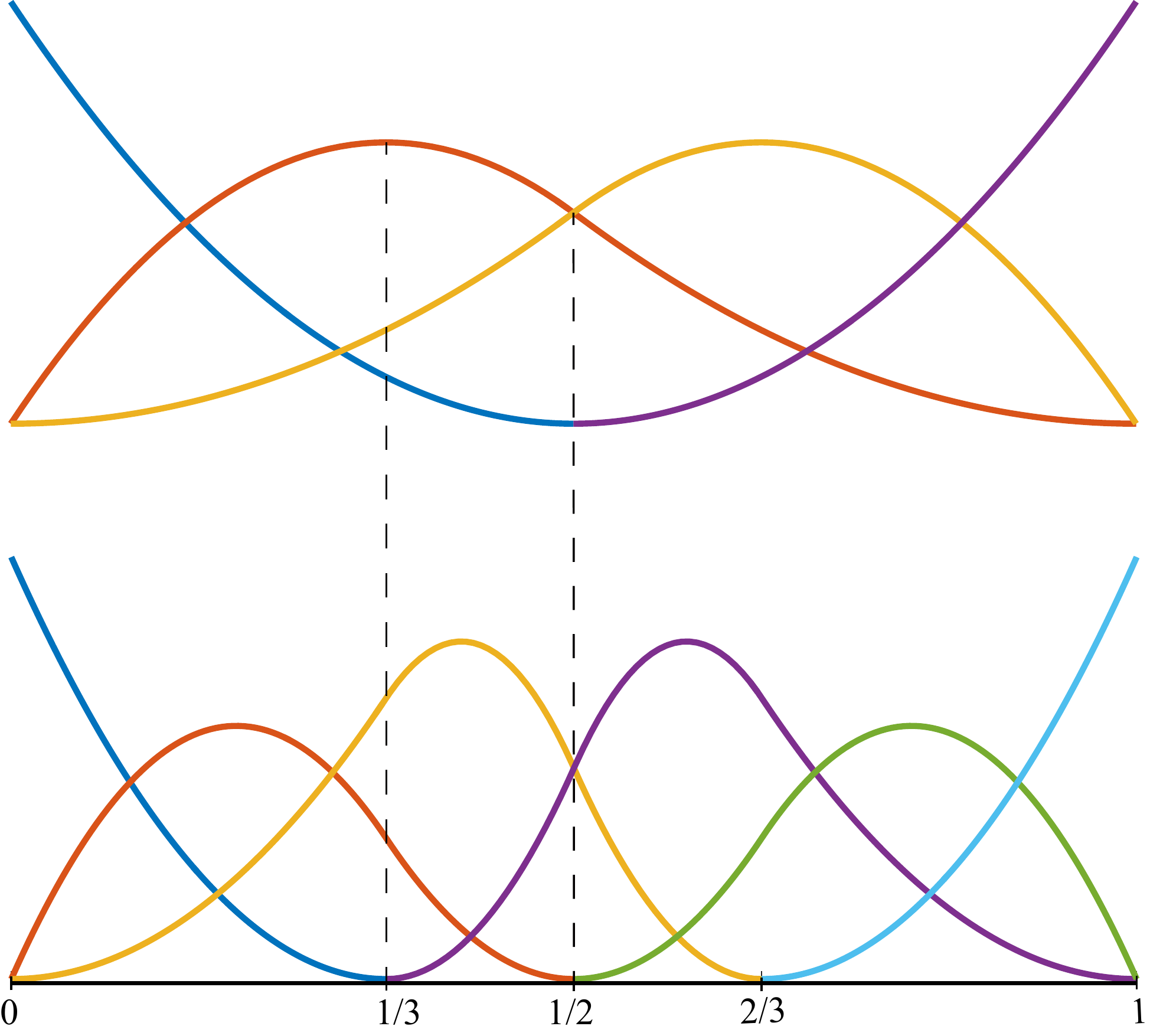}
      \put(5 , 85){\scalebox{0.8}{$N^m_1$}}
      \put(31, 81){\scalebox{0.8}{$N^m_2$}}
      \put(68, 81){\scalebox{0.8}{$N^m_3$}}
      \put(89, 85){\scalebox{0.8}{$N^m_3$}}
      \put(5 , 38){\scalebox{0.8}{$N^r_1$}}
      \put(18, 29){\scalebox{0.8}{$N^r_2$}}
      \put(38, 36){\scalebox{0.8}{$N^r_3$}}
      \put(58, 36){\scalebox{0.8}{$N^r_4$}}
      \put(78, 29){\scalebox{0.8}{$N^r_5$}}
      \put(89, 38){\scalebox{0.8}{$N^r_6$}}
    \end{overpic}
    \caption{Master interface basis (top) and refined
      slave interface basis (bottom).}
  \end{subfigure} \quad \quad
  \begin{subfigure}{0.47\textwidth}
    \begin{overpic}[scale=0.4]{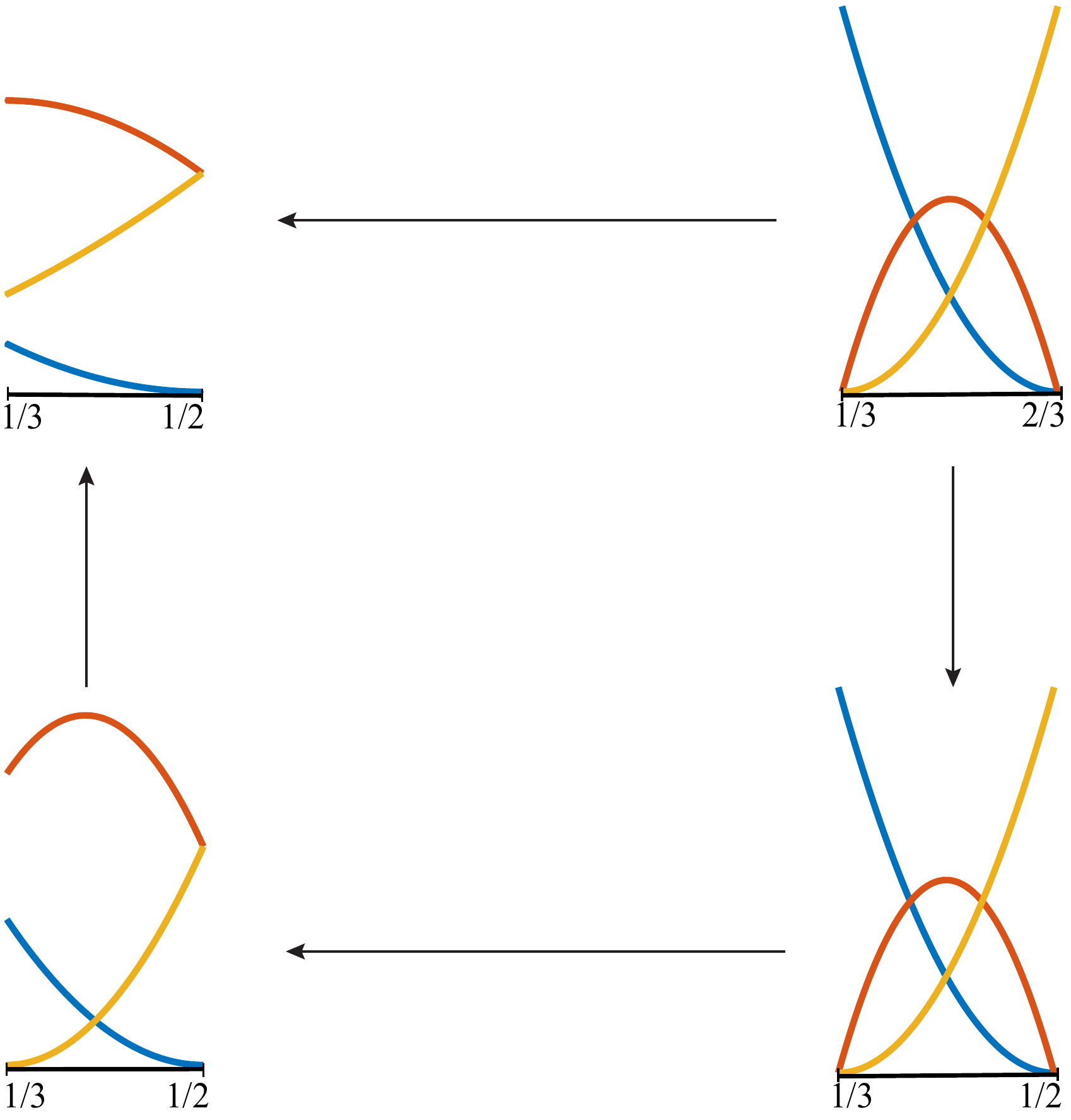}
      \put(12 , 35){\scalebox{0.8}{$\mathbf{N}^{e,r}$}}
      \put(12 , 92){\scalebox{0.8}{$\mathbf{N}^{e,m}$}}
      \put(82 , 32){\scalebox{0.8}{$\mathbf{B}^{e,r}$}}
      \put(82 , 93){\scalebox{0.8}{$\mathbf{B}^{e,s}$}}
      \put(10, 48){\scalebox{0.8}{$\mathbf{G}^e_{\bar{N}^r,N^m}$}}
      \put(46, 17){\scalebox{0.8}{$\mathbf{R}^{e,r}$}}
      \put(71, 48){\scalebox{0.8}{$\mathbf{M}^{-\text{T}}$}}
      \put(23, 84){\scalebox{0.8}{$\tilde{\mathbf{R}}^e = (\mathbf{G}^e_{\bar{N}^r,N^m})^{\text{T}}\mathbf{R}^{e,r}\mathbf{M}^{-\text{T}}$}}
    \end{overpic}
    \caption{Transformation of basis functions for element $e_{21}$.}
  \end{subfigure} %
  \caption{Construction of a refined interface extraction operator for element $e_{21}$ from Figure~\ref{fig:schematics_of_spliting_element}.}	
     \label{fig:basis_transf}
\end{figure}

\begin{figure}
   \centering
   \begin{subfigure}{0.48\textwidth}
    \includegraphics[width=\textwidth]{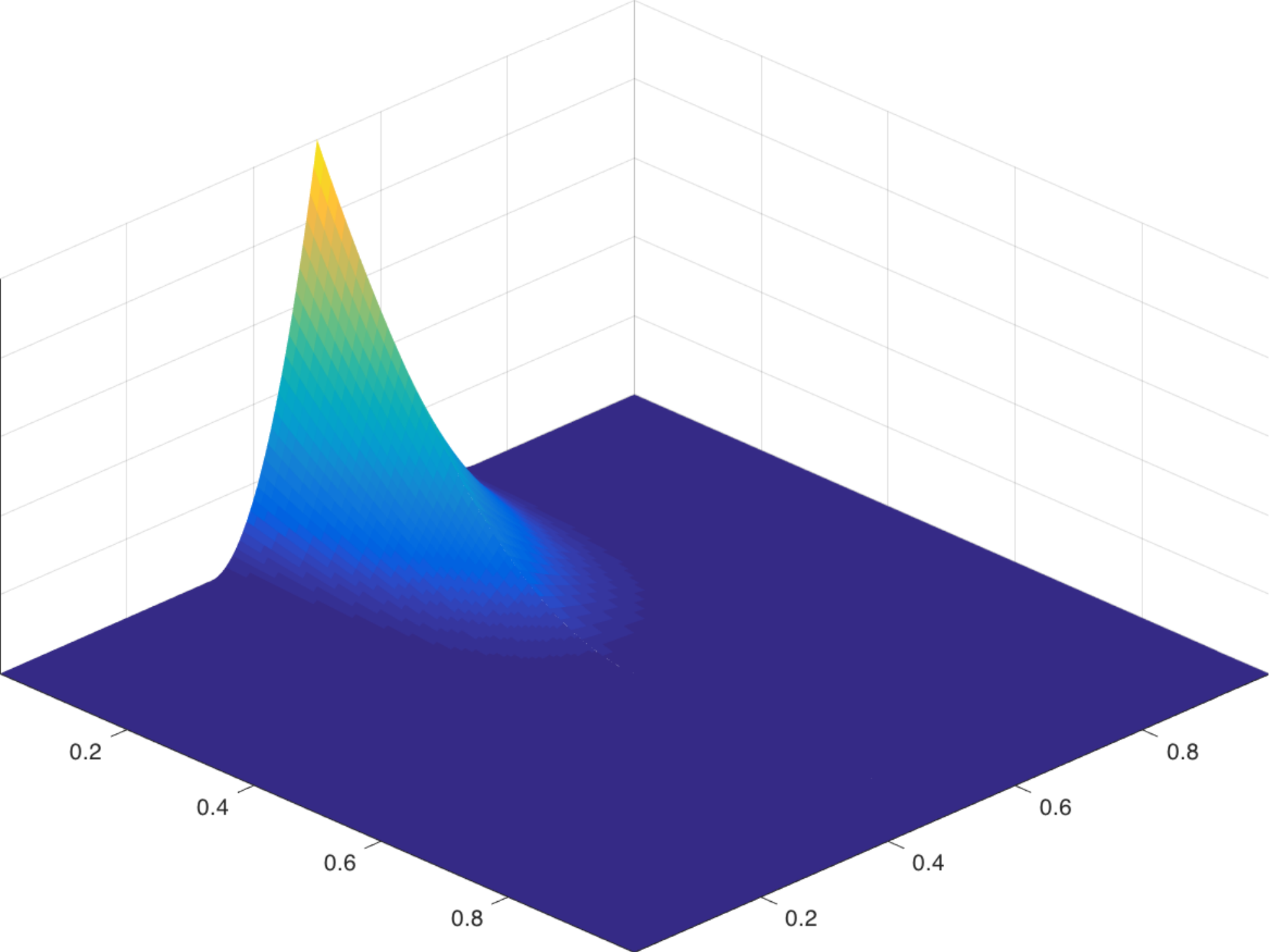} 
     \caption{First interface basis.}	
    \end{subfigure}
   \begin{subfigure}{0.48\textwidth}
    \includegraphics[width=\textwidth]{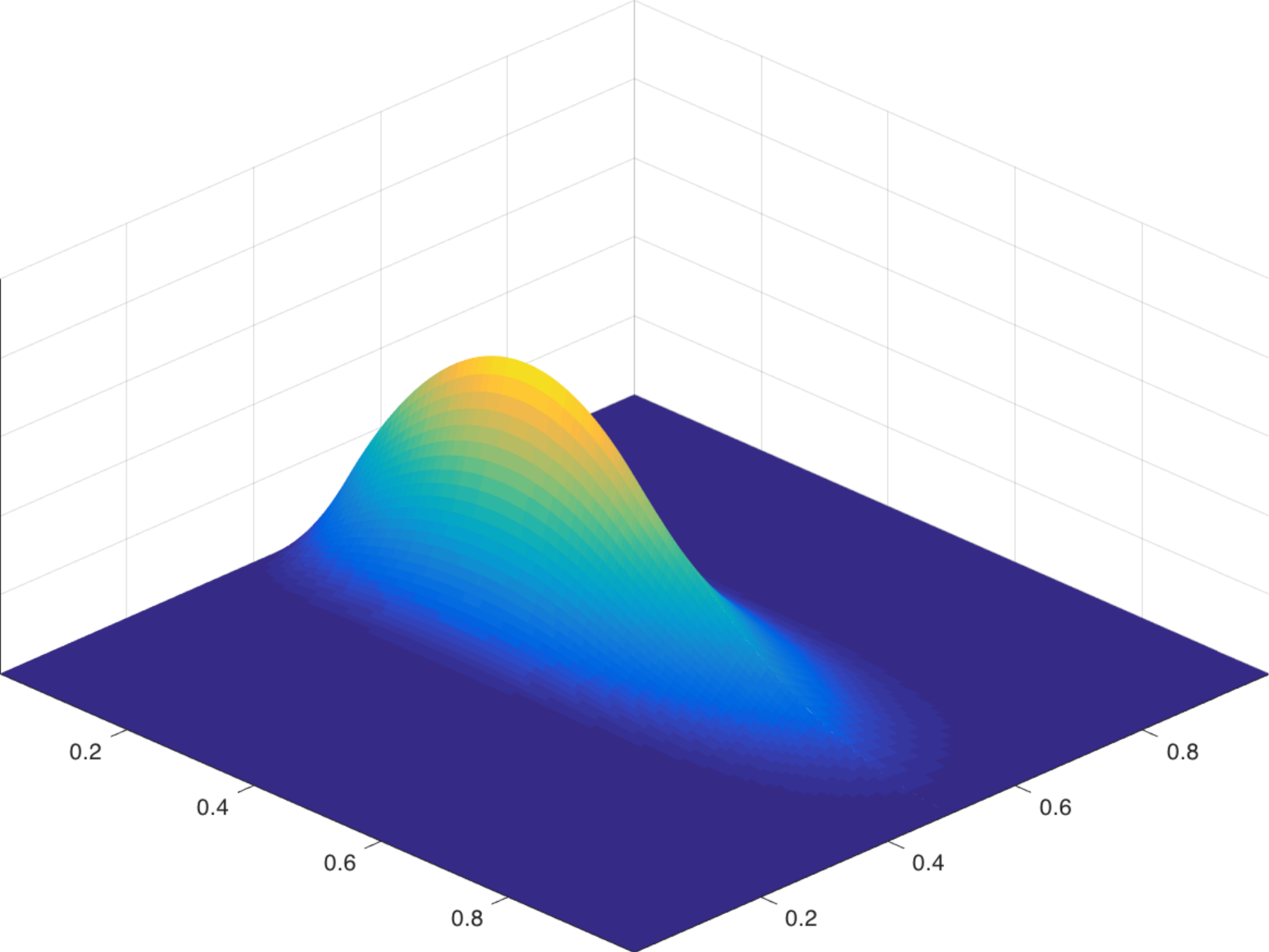} 
     \caption{Second interface basis.}	
    \end{subfigure}
   \begin{subfigure}{0.48\textwidth}
    \includegraphics[width=\textwidth]{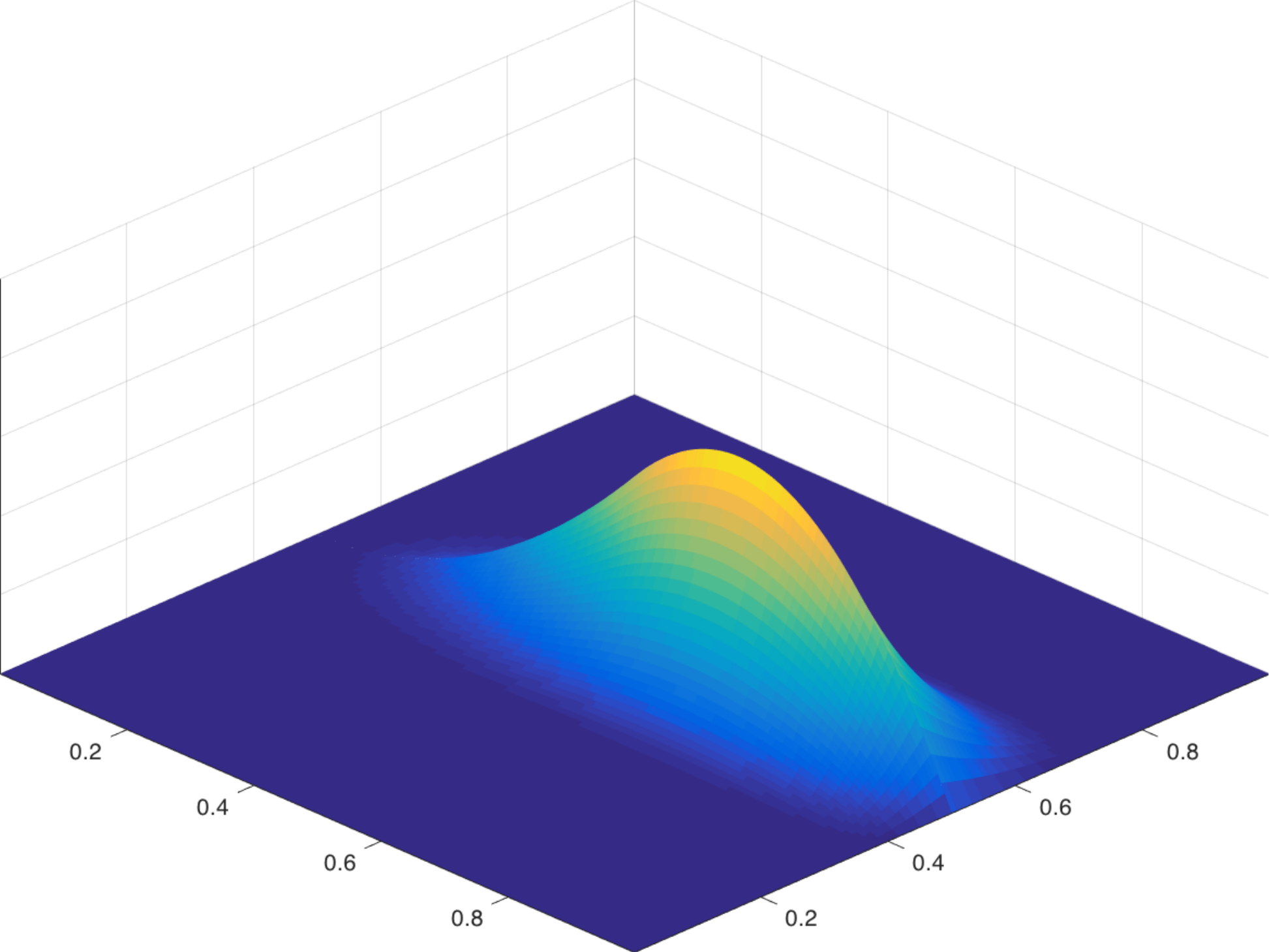} 
     \caption{Third interface basis.}	
    \end{subfigure}
   \begin{subfigure}{0.48\textwidth}
    \includegraphics[width=\textwidth]{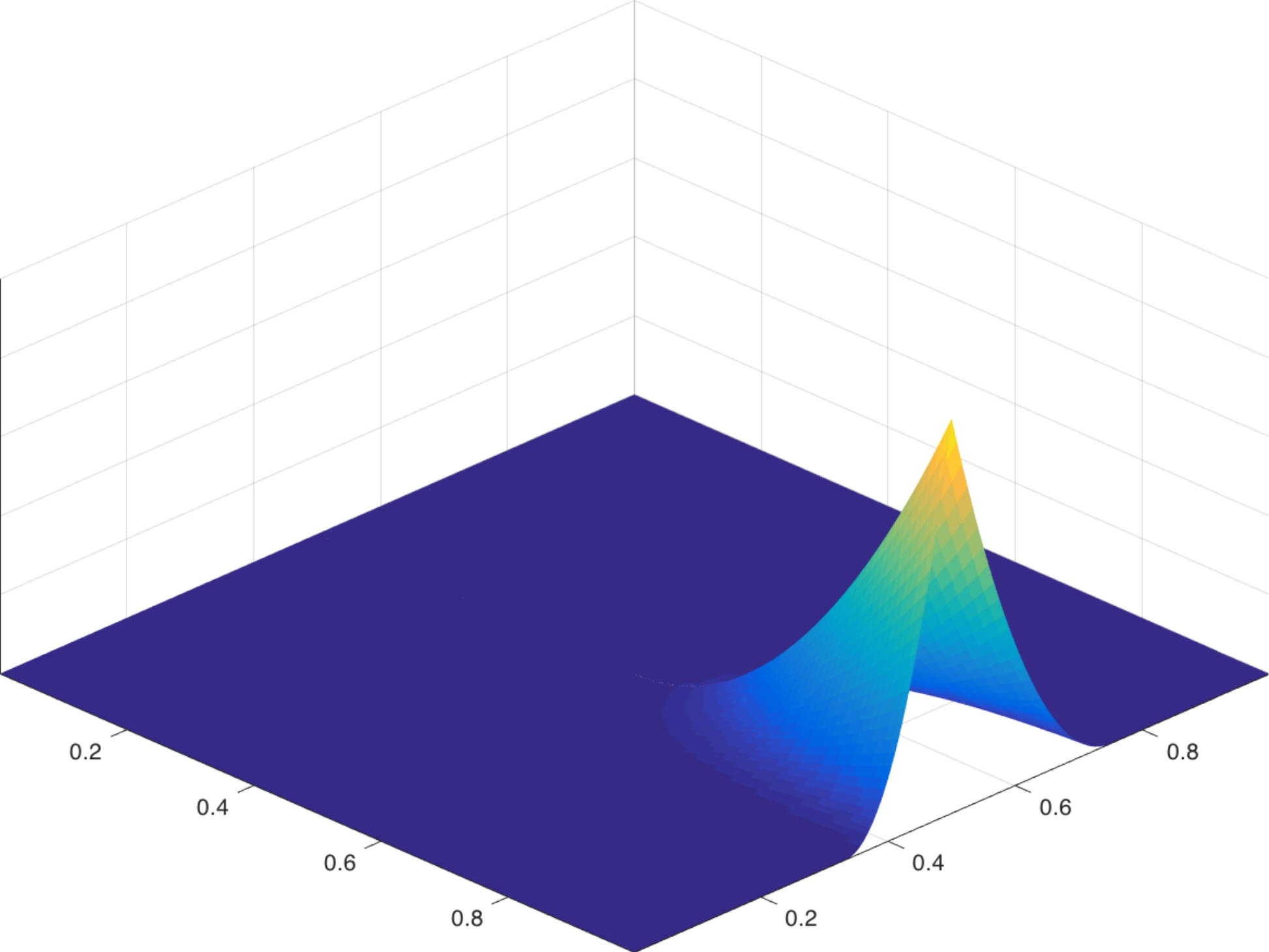} 
     \caption{Fourth interface basis.}	
    \end{subfigure}
     \caption{Weakly continuous basis functions along an interface corresponding to the mesh shown in Figure~\ref{fig:schematics_of_spliting_element}.}	
     \label{fig:patchinterfacebasis}
 \end{figure}

If the slave interface is refined, then (\ref{eq:global_basis_relation}) can be written as
\begin{align}\label{eq:refined_global_basis_relation}
  \mathbf{N}^m &= (\mathbf{G}_{\bar{N}^r, N^m})^\text{T} \mathbf{N}^r
\end{align}
where $\mathbf{N}^r$ is the refined slave interface basis vector. Similarly, (\ref{eq:elem_basis_relation}) can be
written as
\begin{align}
  \mathbf{N}^{m,e} &= (\mathbf{G}^e_{\bar{N}^r, N^m})^\text{T} \mathbf{R}^{e,r} \mathbf{B}^{e,r} = (\mathbf{G}^e_{\bar{N}^r, N^m})^\text{T} \mathbf{R}^{e,r} \mathbf{M}^{-\text{T}} \mathbf{B}^{e,s} 
\end{align}
where $\mathbf{R}^{e,r}$ is the standard element extraction operator defined on
the refined slave interface, and $\mathbf{M}$ is the Bernstein basis transformation
matrix defined in (\ref{eq:Mmatrix}). Therefore, the weakly continuous element extraction operator can be written as
\begin{align}
  \label{eq:refine_wg}
\tilde{\mathbf{R}}^e = (\mathbf{G}^e_{\bar{N}^r, N^m})^\text{T} \mathbf{R}^{e,r} \mathbf{M}^{-\text{T}}.
\end{align}
Figure~\ref{fig:basis_transf} shows the action of (\ref{eq:refine_wg}) for the interface element $e_{21}$ in
Figure~\ref{fig:schematics_of_spliting_element} and the resulting weakly continuous two-dimensional basis functions along the interface are shown in Figure~\ref{fig:patchinterfacebasis}. The full expressions for the weakly continuous element
extraction operators $\tilde{\mathbf{R}}^e$ are given in \ref{app:weak_extra_op}.

\section{Numerical results}
\label{sec:results}
We evaluate the performance of the \Bezier dual mortar method on several
benchmark problems. In all cases, dual basis refinement will be employed, where,
during the first step of refinement the master knots are projected into the
slave interface. Note that this initial refinement step is always possible for
any master/slave pairing. Subsequent refinements then utilize element splitting
of the slave interface. Note that refinement of the dual space does not increase
the problem number of degrees-of-freedom. To avoid the well-known mortar integral error in calculating
$\mathbf{G}_{\bar{N}, N^m}$ in~(\ref{discretedispcons}) when no refinement of
the dual basis is employed, the master knots are still
projected into the slave interface, and then the integration is performed on the
combined knot intervals as described in~\cite{Dornisch2017449}. We compare our method to a global dual
mortar method~\cite{NME:NME4918,Dornisch2017449}, where the global dual basis is
computed using $L^2$ projection.  
 \subsection{A manufactured solution on a square domain} \label{square_twopatches}
 \label{test1}
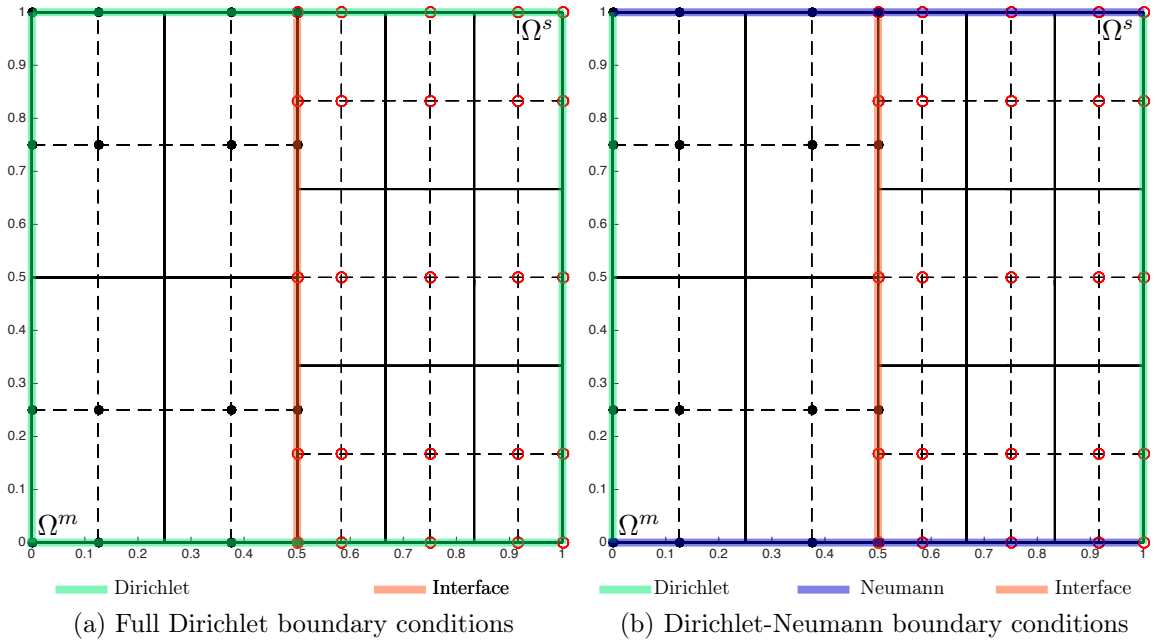
\begin{figure}[h]
   \centering
   \begin{subfigure}{0.46\textwidth}
     \begin{overpic}[scale=0.5]{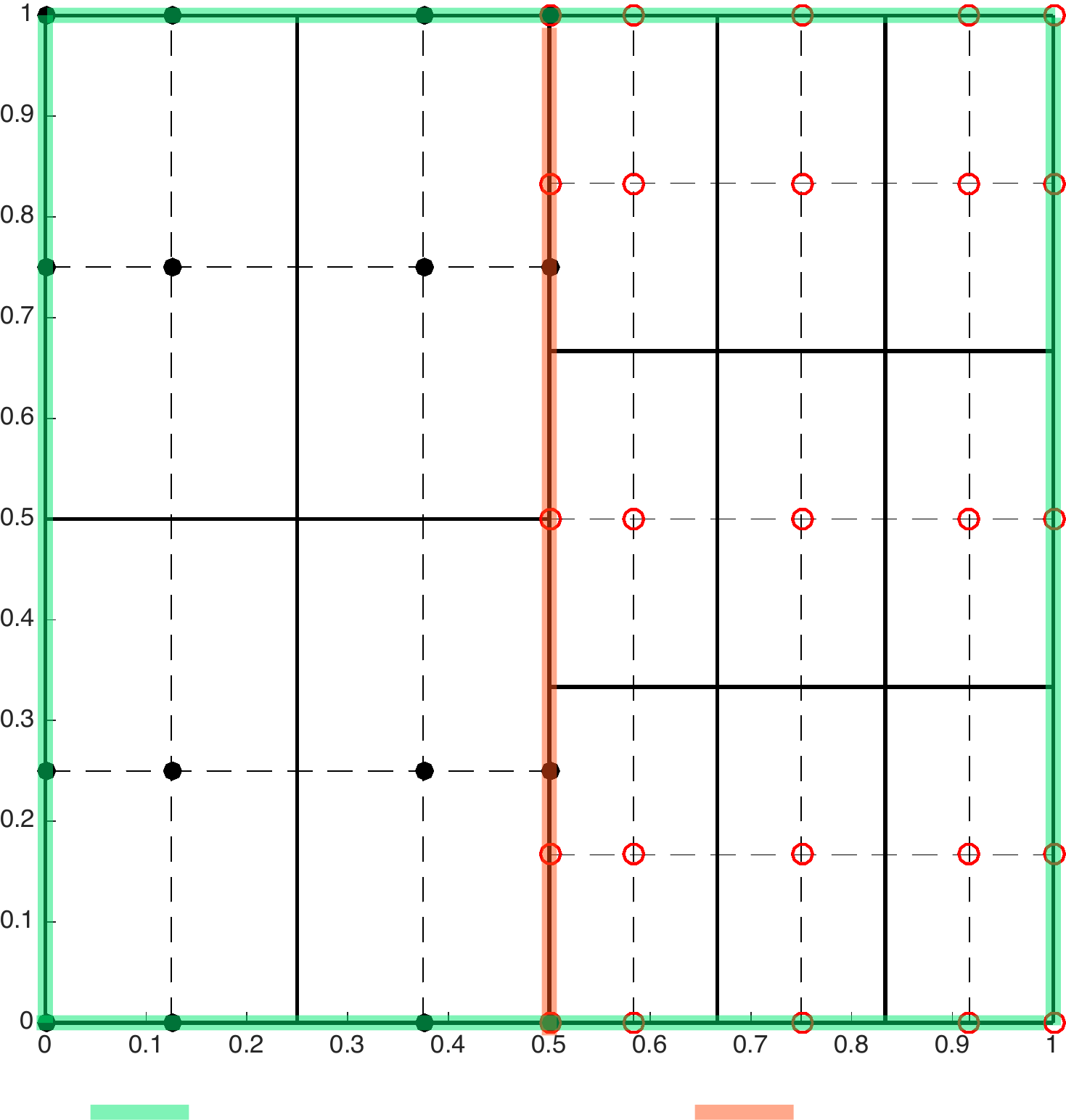}
    \put(18, 0){\scalebox{0.7}{Dirichlet}}
    \put(72, 0){\scalebox{0.7}{Interface}}
    \put(72, 0){\scalebox{0.7}{Interface}}
    \put(5, 10){\scalebox{1}{$\Omega^m$}}
    \put(87, 94){\scalebox{1}{$\Omega^s$}}
    \end{overpic}
    \caption{Full Dirichlet boundary conditions}
    \end{subfigure}
   \begin{subfigure}{0.46\textwidth}
     \begin{overpic}[scale=0.5]{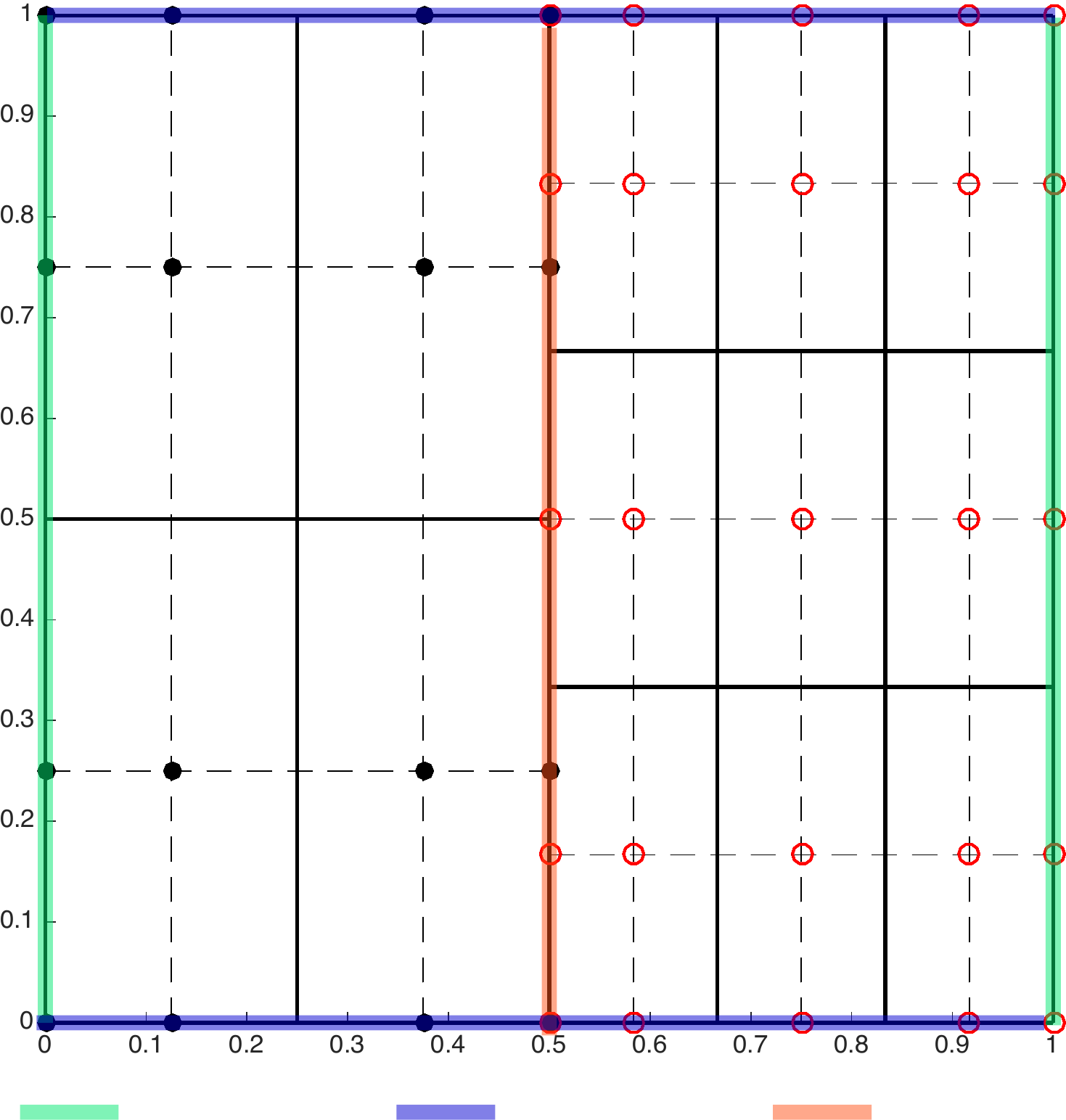}
    \put(11, 0){\scalebox{0.7}{Dirichlet}}
    \put(46, 0){\scalebox{0.7}{Neumann}}
    \put(79, 0){\scalebox{0.7}{Interface}}
    \put(5, 10){\scalebox{1}{$\Omega^m$}}
    \put(87, 94){\scalebox{1}{$\Omega^s$}}
    \end{overpic}
    \caption{Dirichlet-Neumann boundary conditions}
    \end{subfigure}
    \caption{Two quadratic maximally smooth nonconforming B-spline patches.}	
 \label{fig:thermalsquareschematic}
 \end{figure}

 We first solve the Laplace equation, $-\Delta u =  0$, on the square domain,
 $\Omega = \left(0,1\right )\times \left(0,1\right)$.  The domain is modeled
 with two maximally smooth quadratic B-spline patches where the left patch is
 the master and the right patch is the slave as shown in Figure
 \ref{fig:thermalsquareschematic}. Two
 different boundary conditions, shown in Figure
 \ref{fig:thermalsquareschematic}a and b, are considered to demonstrate that, in contrast
 to the global dual mortar method, the \Bezier dual mortar method does not
 suffer from the so-called crosspoint problem~\cite{Dornisch2017449,Brivadis2015292}. This superior performance is due to the locality of the dual basis
 functions. Both boundary conditions satisfy the manufactured solution,
 $u\left(x,y\right)= \sin(\pi y)\sinh( \pi x)$. The ratio of master to slave element size is initially chosen to be $2:3$. The master and slave interface boundaries are matched but the underlying meshes are nonconforming.

\begin{figure}[h]
\centering
    \begin{subfigure}{0.327\textwidth}
        \includegraphics[width=\textwidth]{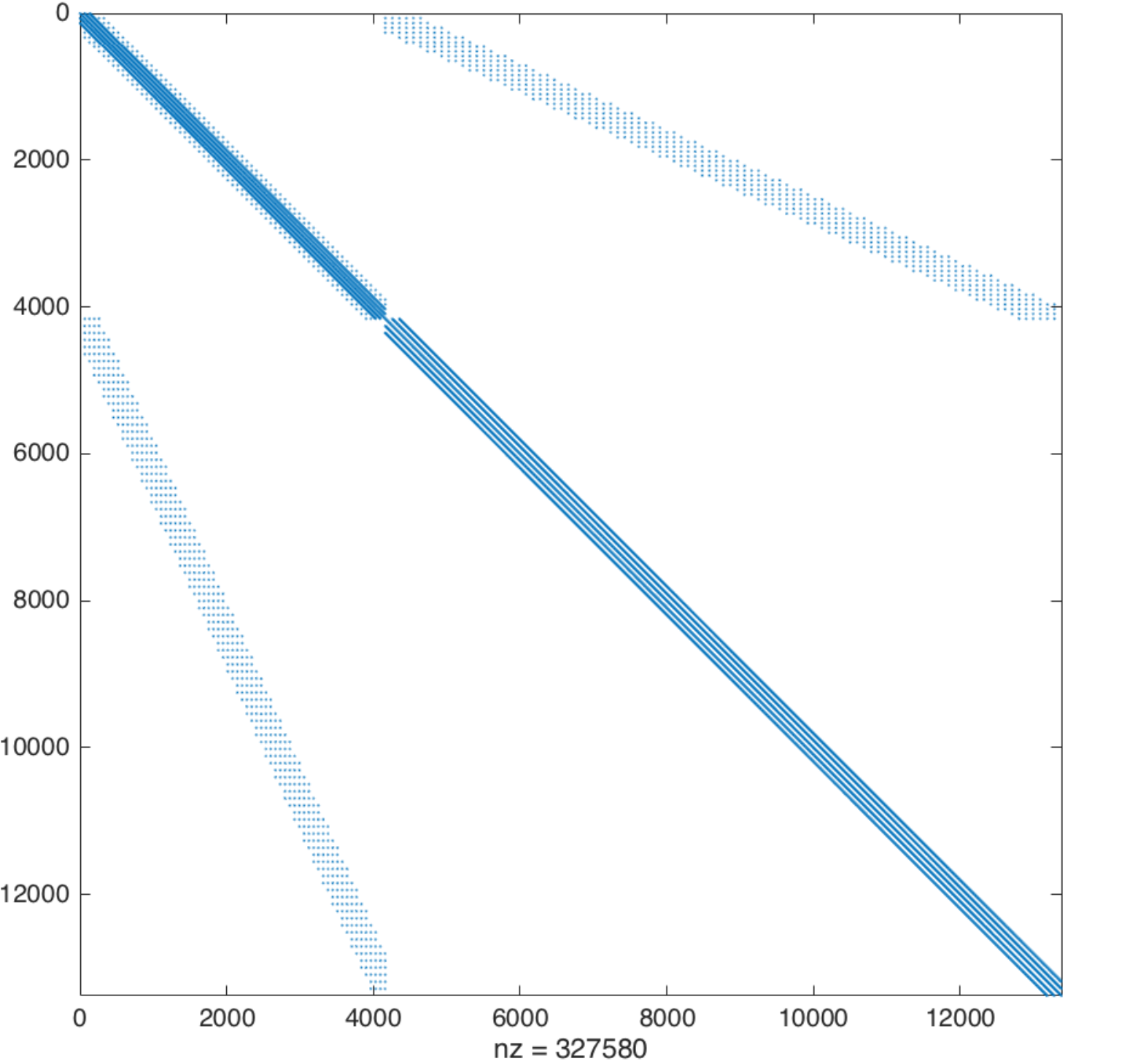}
        \caption{The proposed method}
    \end{subfigure}
    \begin{subfigure}{0.327\textwidth}
        \includegraphics[width=\textwidth]{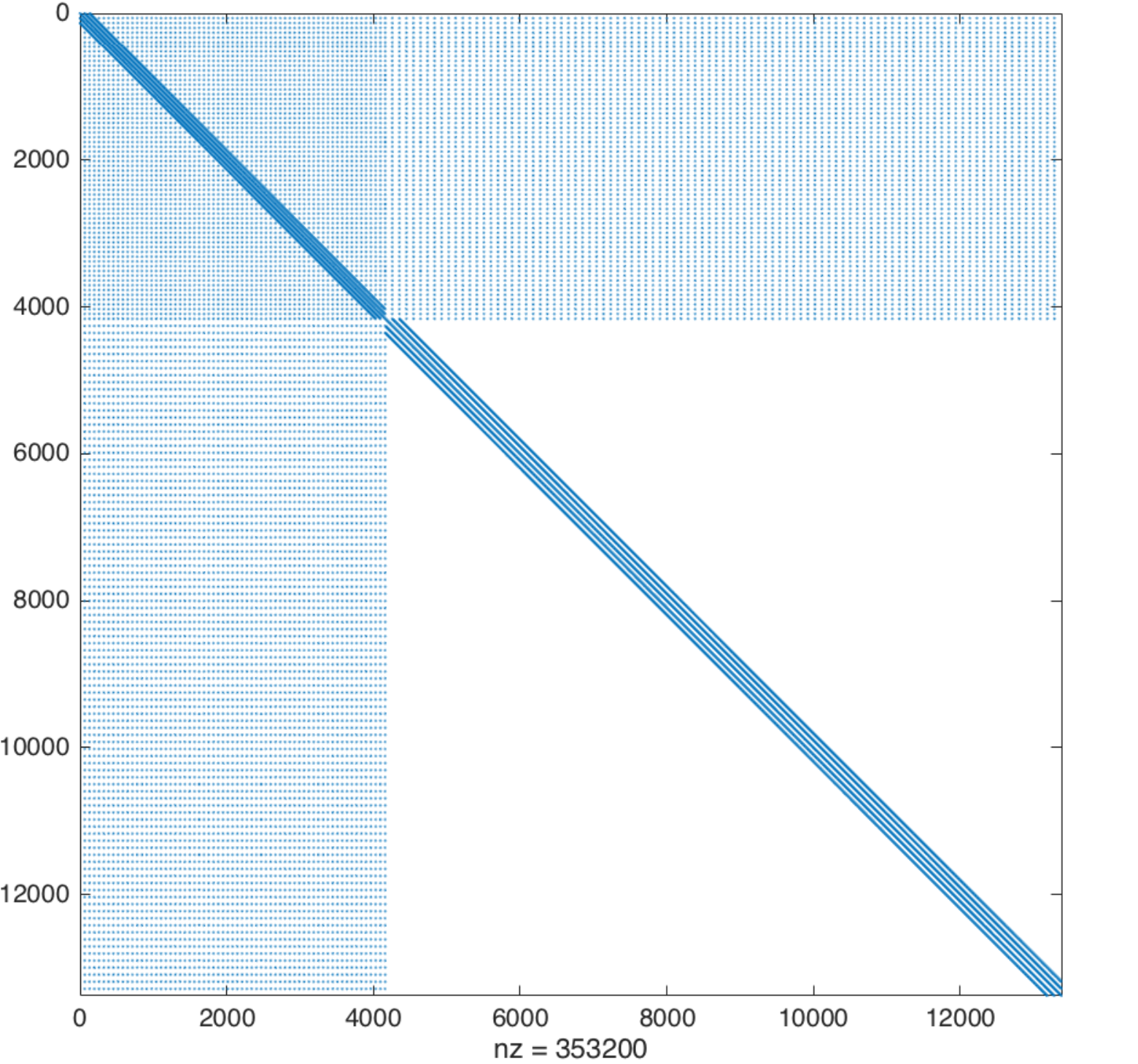}
        \caption{A global dual mortar method}
    \end{subfigure} %
    \begin{subfigure}{0.327\textwidth}
        \includegraphics[width=\textwidth]{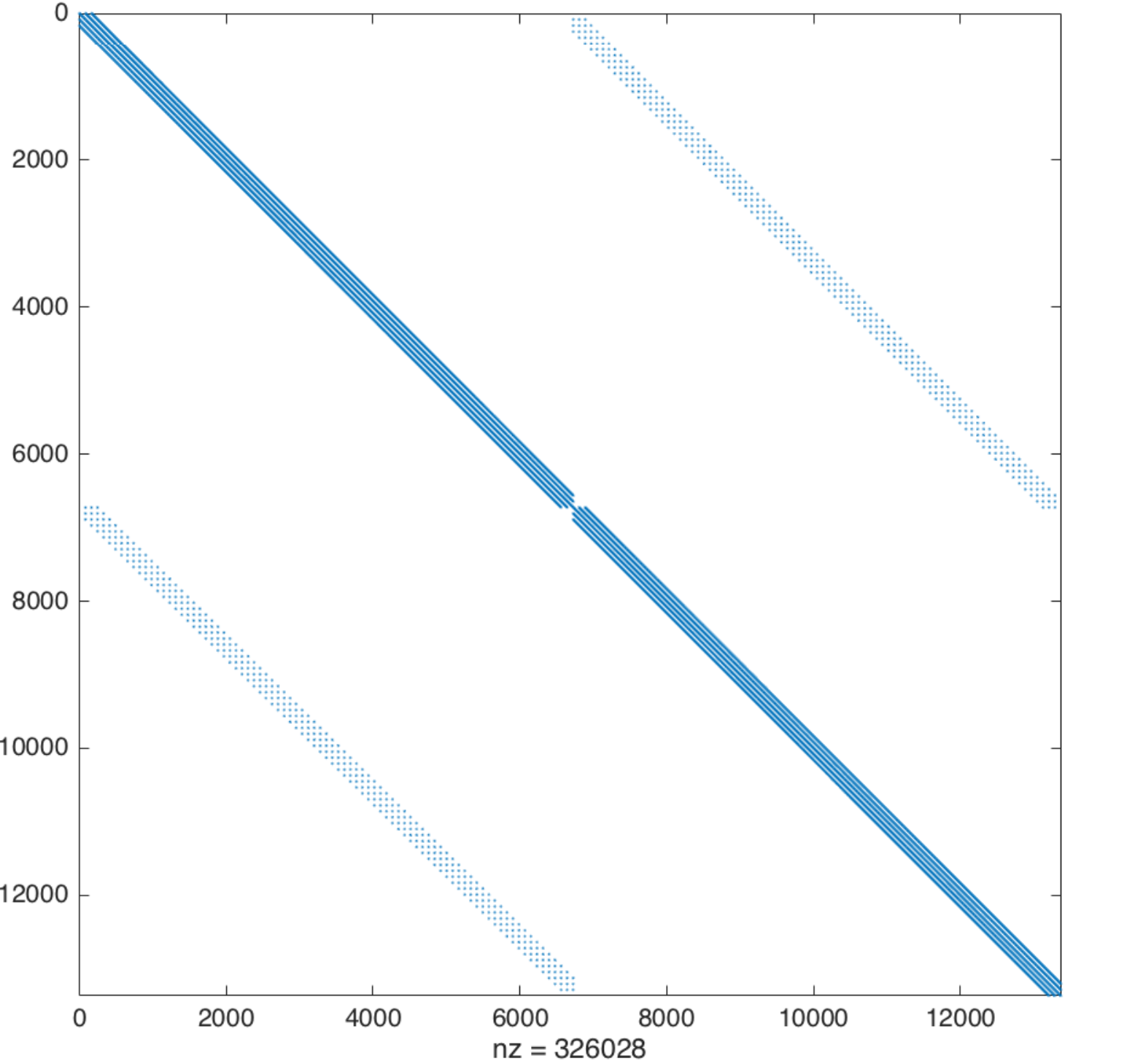}
        \caption{A conforming method}
    \end{subfigure}
    \caption{Stiffness matrix sparsity patterns for (a) the proposed method, (b) a global dual mortar method, and (c) a standard conforming method. The nonconforming examples (i.e., (a) and (b)) are generated from nonconforming meshes after uniformly refining the initial mesh shown in Figure~\ref{fig:thermalsquareschematic} four times. The resulting mesh has 13862 nodes. The conforming example in (c) is generated from a mesh with 13860 nodes.}	
    \label{fig:sparseK}
\end{figure}

 The sparsity patterns for the stiffness matrices for the proposed method, a global dual mortar method, and a standard conforming method are shown in Figure~\ref{fig:sparseK}a, \ref{fig:sparseK}b, and \ref{fig:sparseK}c, respectively, after four applications of uniform global refinement. It is clear that the proposed method generates a sparse stiffness matrix with only a slight increase in bandwidth when compared to a conforming method.

\begin{figure}[h]
 \centering
     \begin{subfigure}{0.75\textwidth}
        \includegraphics[width=\textwidth]{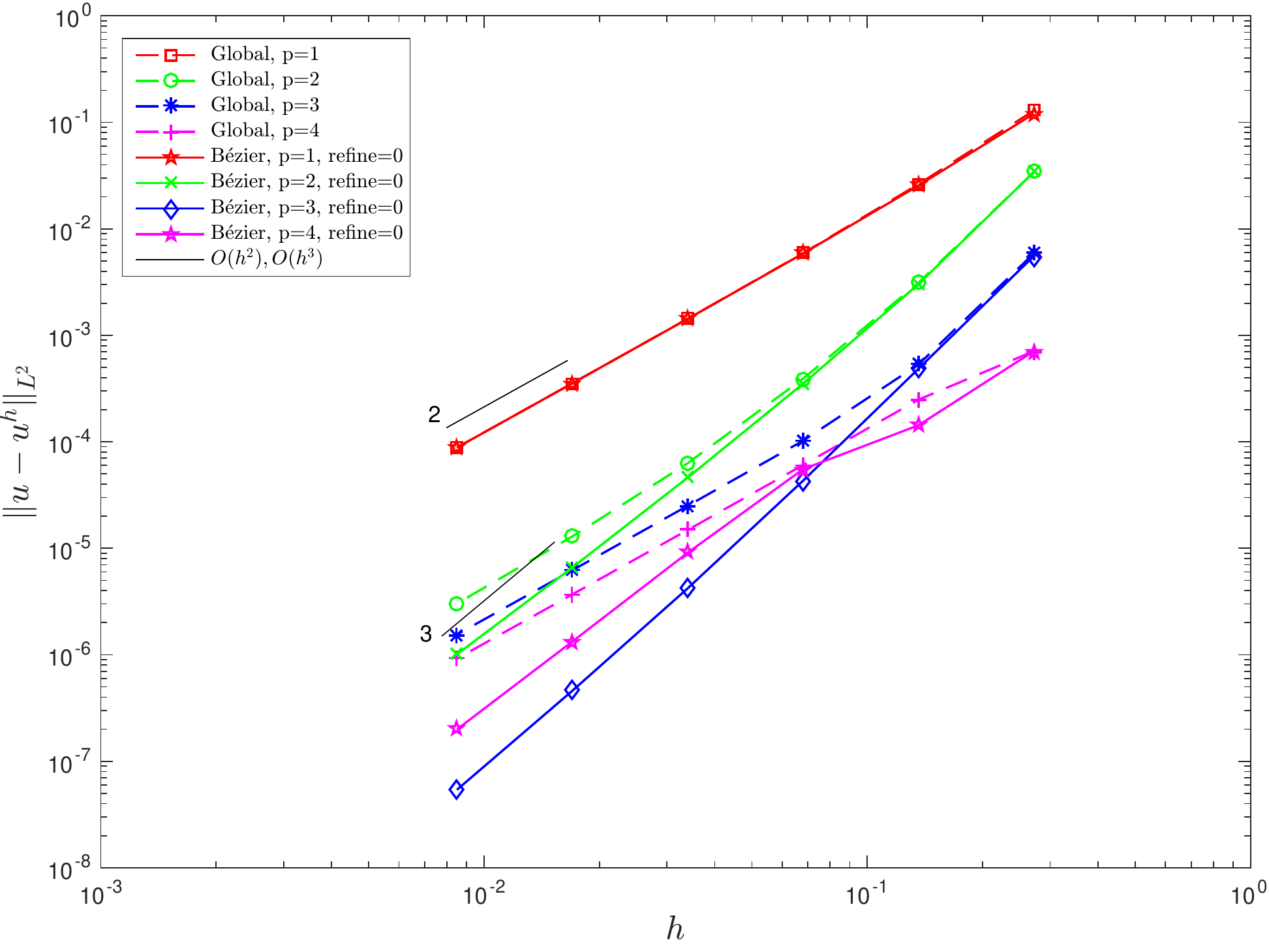}
        \caption{Master/slave mesh ratio, $m:s=2:3$}
    \end{subfigure} \\
    \begin{subfigure}{0.75\textwidth}
        \includegraphics[width=\textwidth]{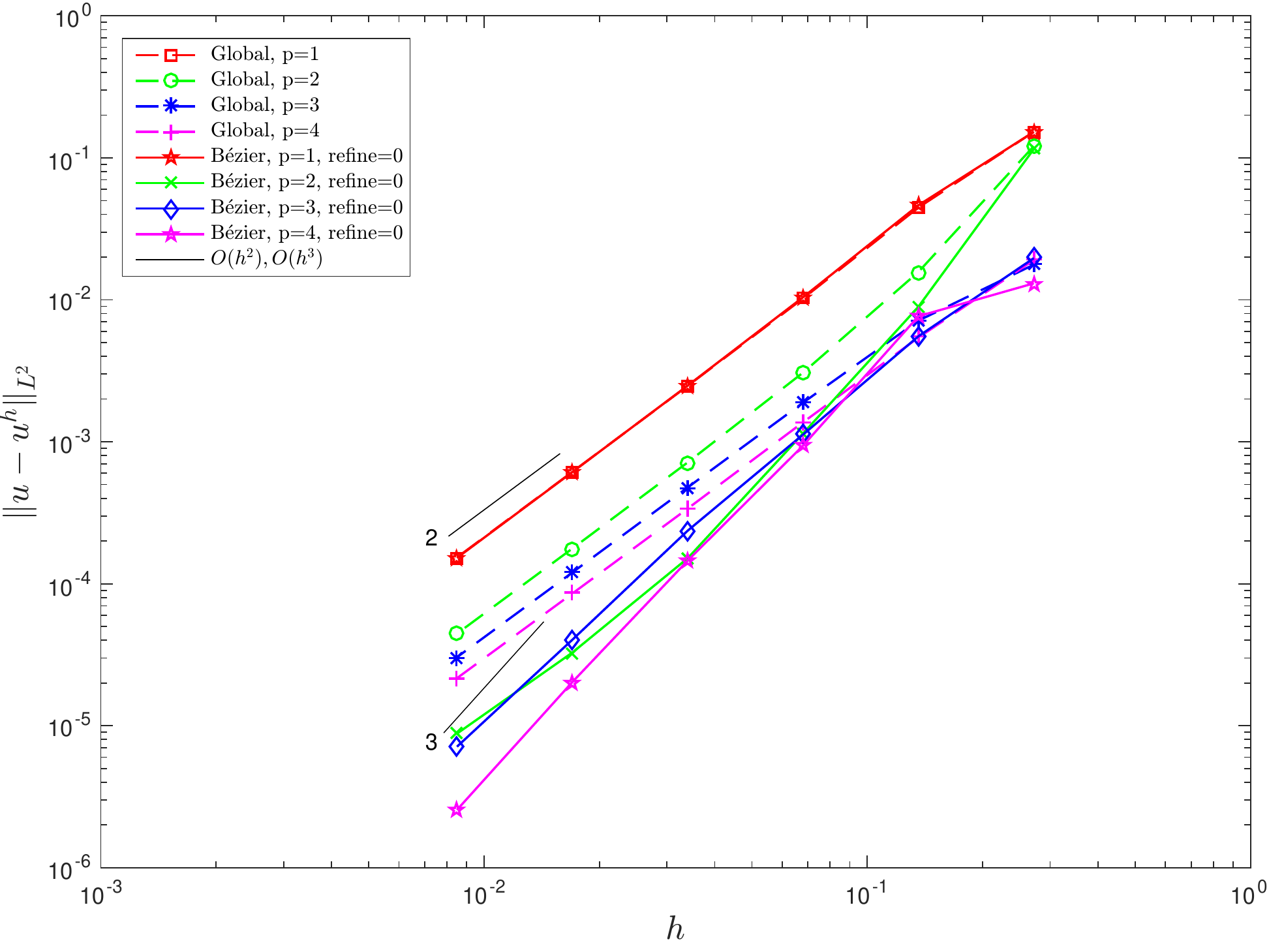}
        \caption{Master/slave mesh ration, $m:s=3:2$}
    \end{subfigure}
      \caption{Convergence rates for a square
        domain with two non-conforming patches, full Dirichlet boundary
        conditions (see Figure~\ref{fig:thermalsquareschematic}a) and matched parameterizations.}\label{fig:thermal2DsquareIGSE}
\end{figure}

We first present results for the boundary condition shown in
Figure \ref{fig:thermalsquareschematic}a. The
convergence rates of the displacement error in the $L^2$-norm for several
different degrees are shown in Figure~\ref{fig:thermal2DsquareIGSE}a. The
proposed approach is compared to a global dual mortar method. As can be seen,
the global dual mortar method only gets optimal rates for $p=1$, and for $p = 2,3,4$
the convergence rates are 2, which is suboptimal. This reduction in the rates
can be attributed to the crosspoint problem~\cite{Dornisch2017449,Brivadis2015292}. In other words, the Lagrange multiplier space is bigger
than the primal space due to the crosspoint, which, in this case, corresponds to the points where the interface and Dirichlet boundary conditions intersect. As a result, inf-sup stability is lost.
\textit{Without refinement of the dual basis}, the proposed method achieves optimal rates for $p=1$, and slightly deteriorated rates for $p=2$. For $p=3,4$, the convergence rates are
reduced but still converge faster than the global dual
mortar method. To demonstrate the insensitivity of the method to master and slave
selection, we change the mesh ratio to $m:s = 3:2$. The convergence rates are
shown in Figure \ref{fig:thermal2DsquareIGSE}b. As can be seen, the convergence
rates are close to the previous case, $m:s = 2:3$.

\begin{figure}[h]
 \centering
\begin{subfigure}{0.75\textwidth}
\includegraphics[width=\textwidth]{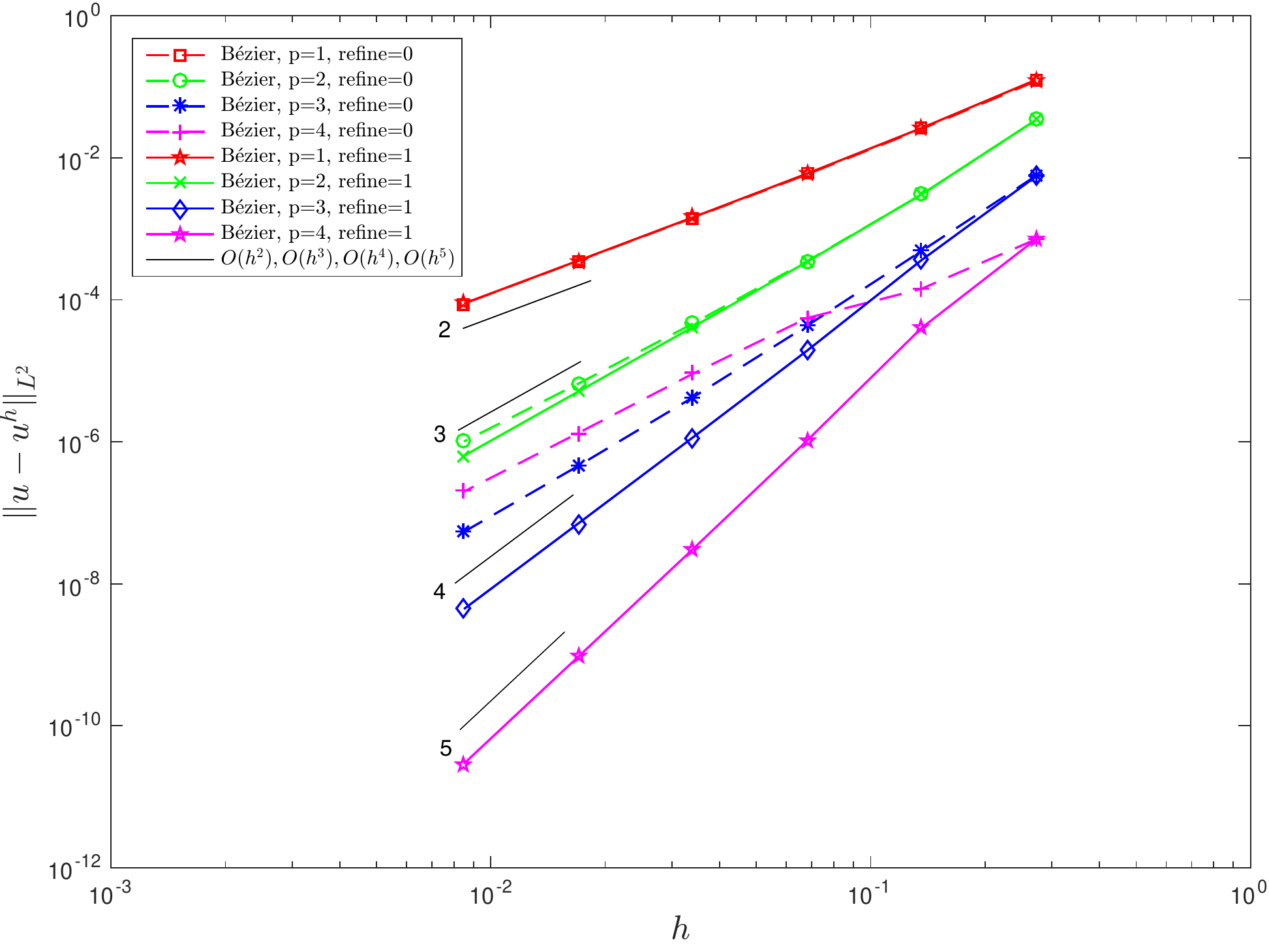} 
\caption{Master/slave mesh ratio, $m:s=2:3$}
\end{subfigure} \\
\begin{subfigure}{0.75\textwidth}
\includegraphics[width=\textwidth]{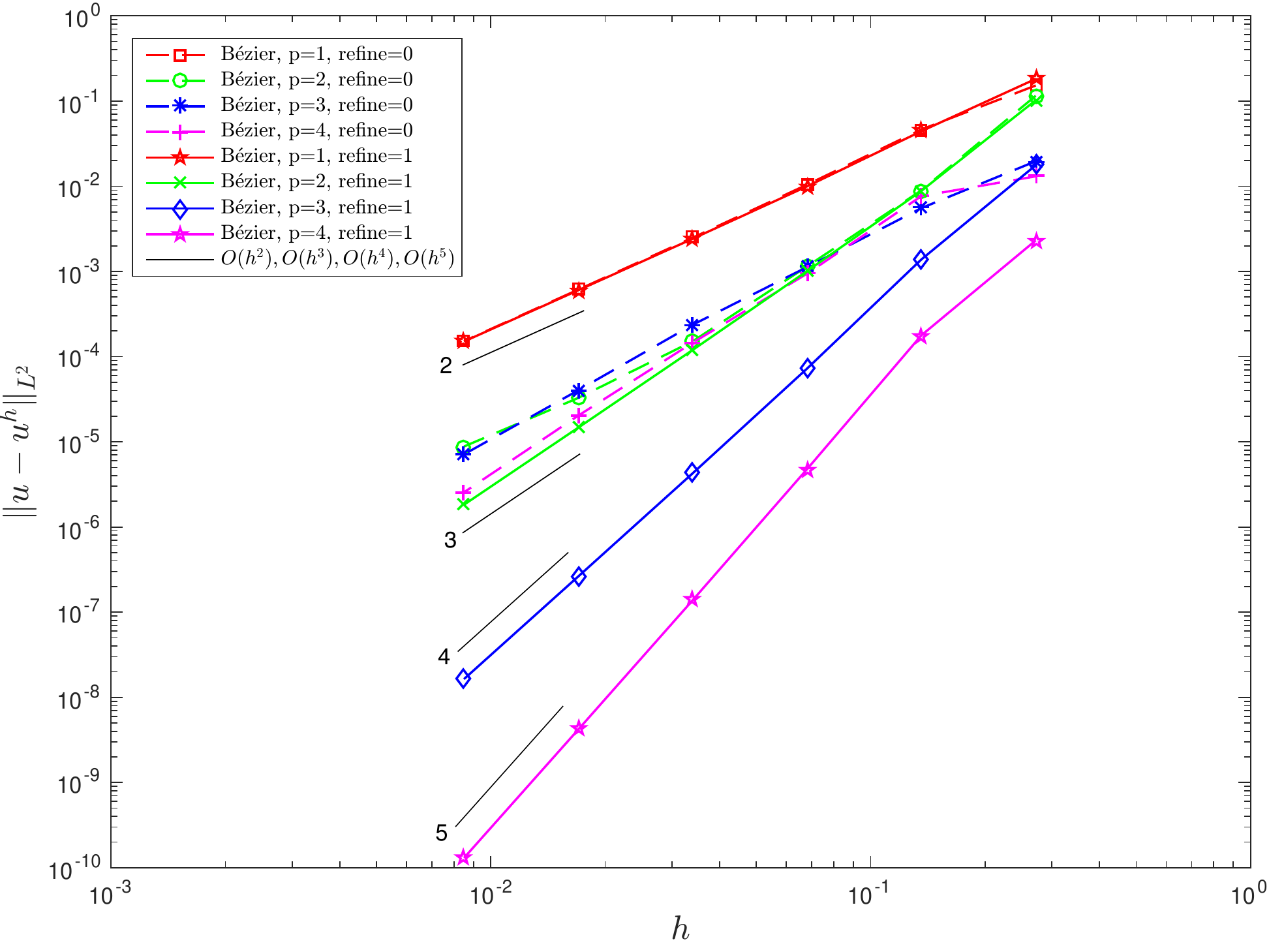}
\caption{Master/slave mesh ratio, $m:s=3:2$}
\end{subfigure}
\caption{Convergence rates for a square domain with two non-conforming patches,
  full Dirichlet boundary conditions (see Figure~\ref{fig:thermalsquareschematic}a) and matched parameterizations where the dual space is refined $n$ times, $n=0,1$.}	
\label{fig:thermal2DsquareIGSEmatchedrefine}
\end{figure}

We now refine the proposed dual space to improve the accuracy. The convergence rates are shown in Figure~\ref{fig:thermal2DsquareIGSEmatchedrefine}. As expected, with one
refinement of the interface dual basis the proposed method obtains optimal convergence rates for all degrees
$p=1,2, 3$, and $4$ for both mesh ratios $m:s=2:3$ and $m:s=3:2$. 

\begin{figure}[h]
 \centering
     \begin{subfigure}{0.75\textwidth}
        \includegraphics[width=\textwidth]{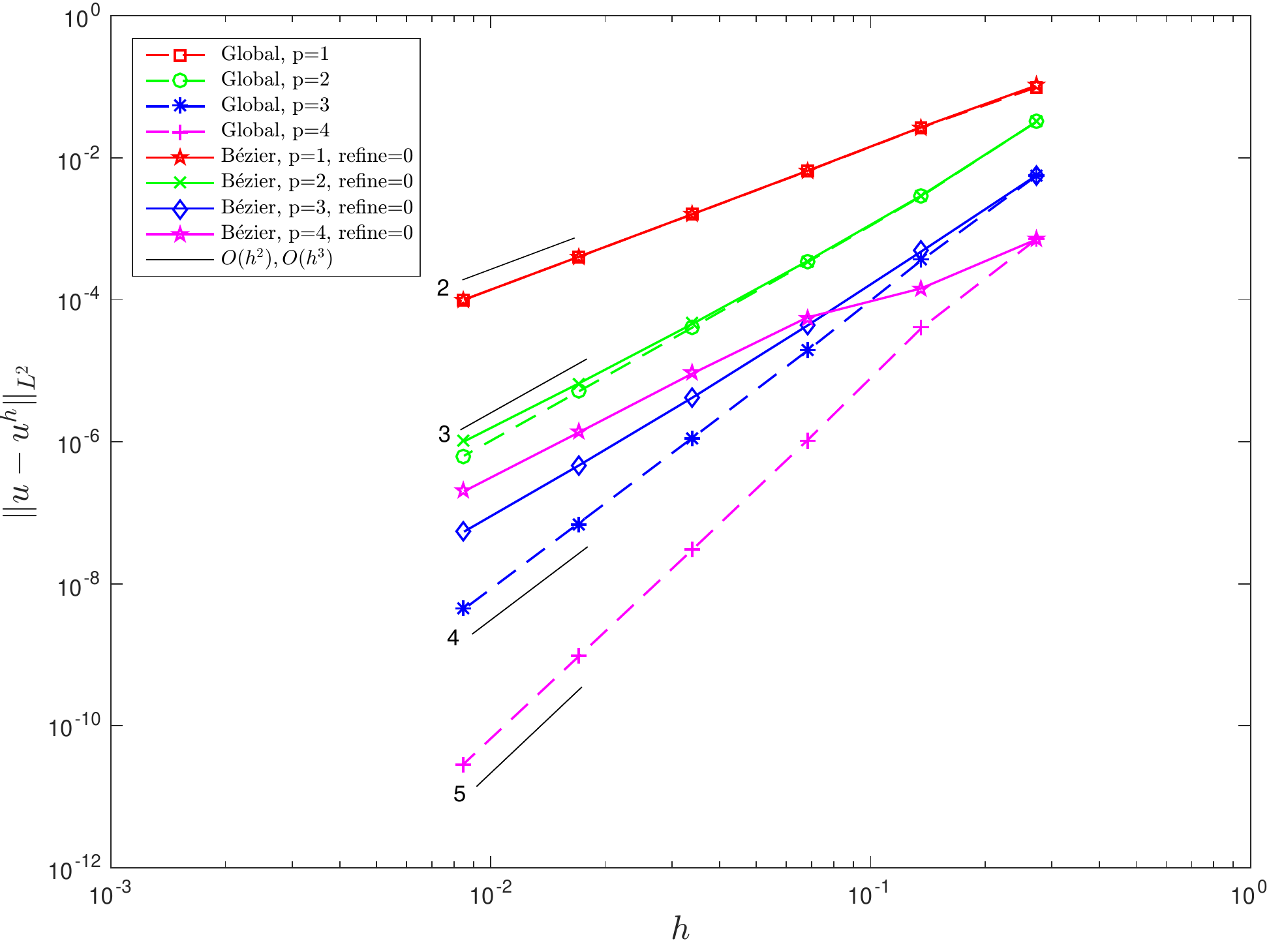}
        \caption{Master/slave mesh ratio, $m:s=2:3$}
    \end{subfigure} \\
    \begin{subfigure}{0.75\textwidth}
        \includegraphics[width=\textwidth]{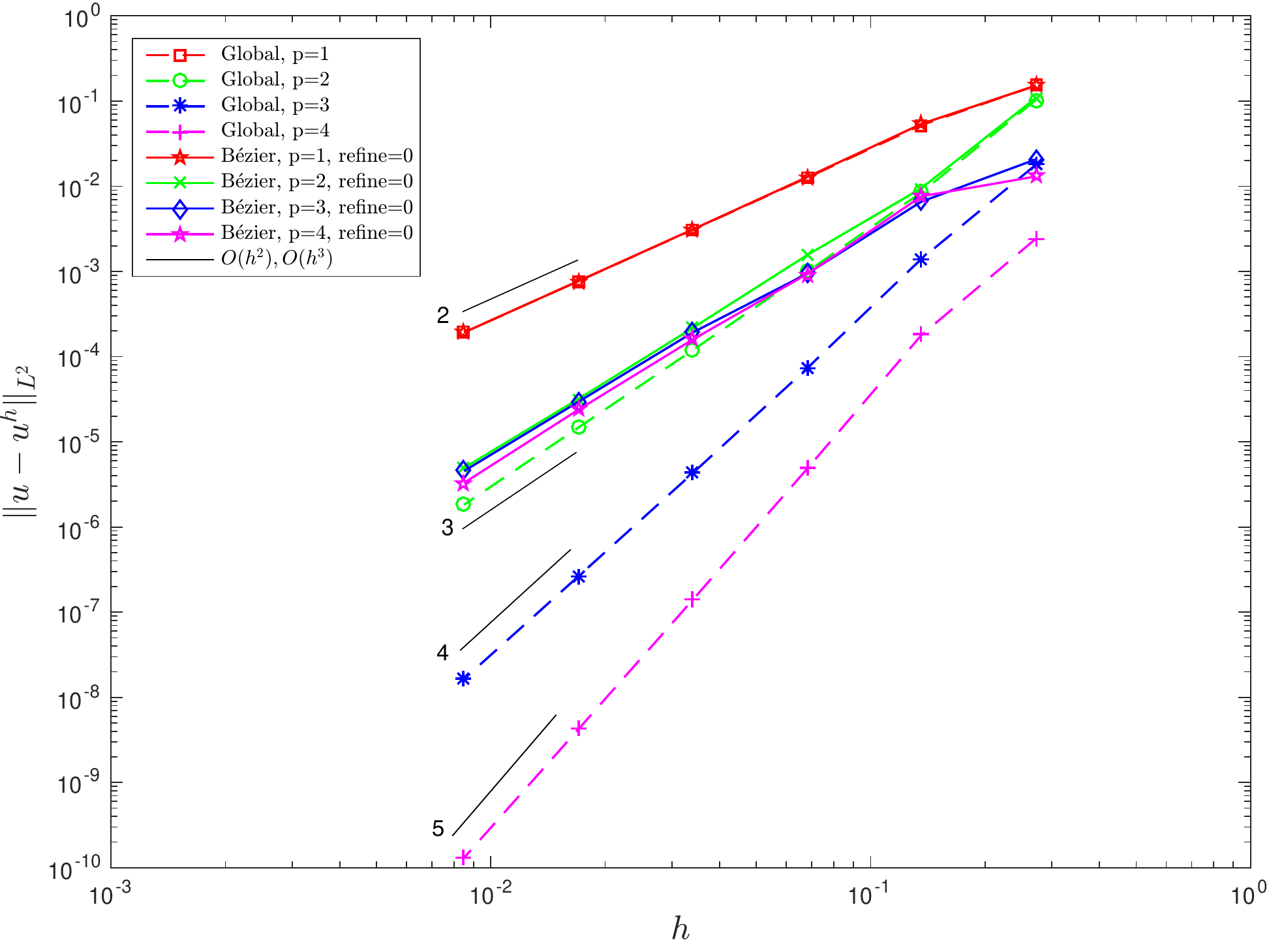}
        \caption{Master/slave mesh ration, $m:s=3:2$}
    \end{subfigure}
      \caption{Convergence rates for a square domain with two non-conforming
        patches, Dirichlet-Neumann
        boundary conditions (see Figure~\ref{fig:thermalsquareschematic}b),
        and matched parameterizations.}\label{fig:thermal2DsquareIGSENeumann}
\end{figure}

\begin{figure}[h]
 \centering
\begin{subfigure}{0.75\textwidth}
\includegraphics[width=\textwidth]{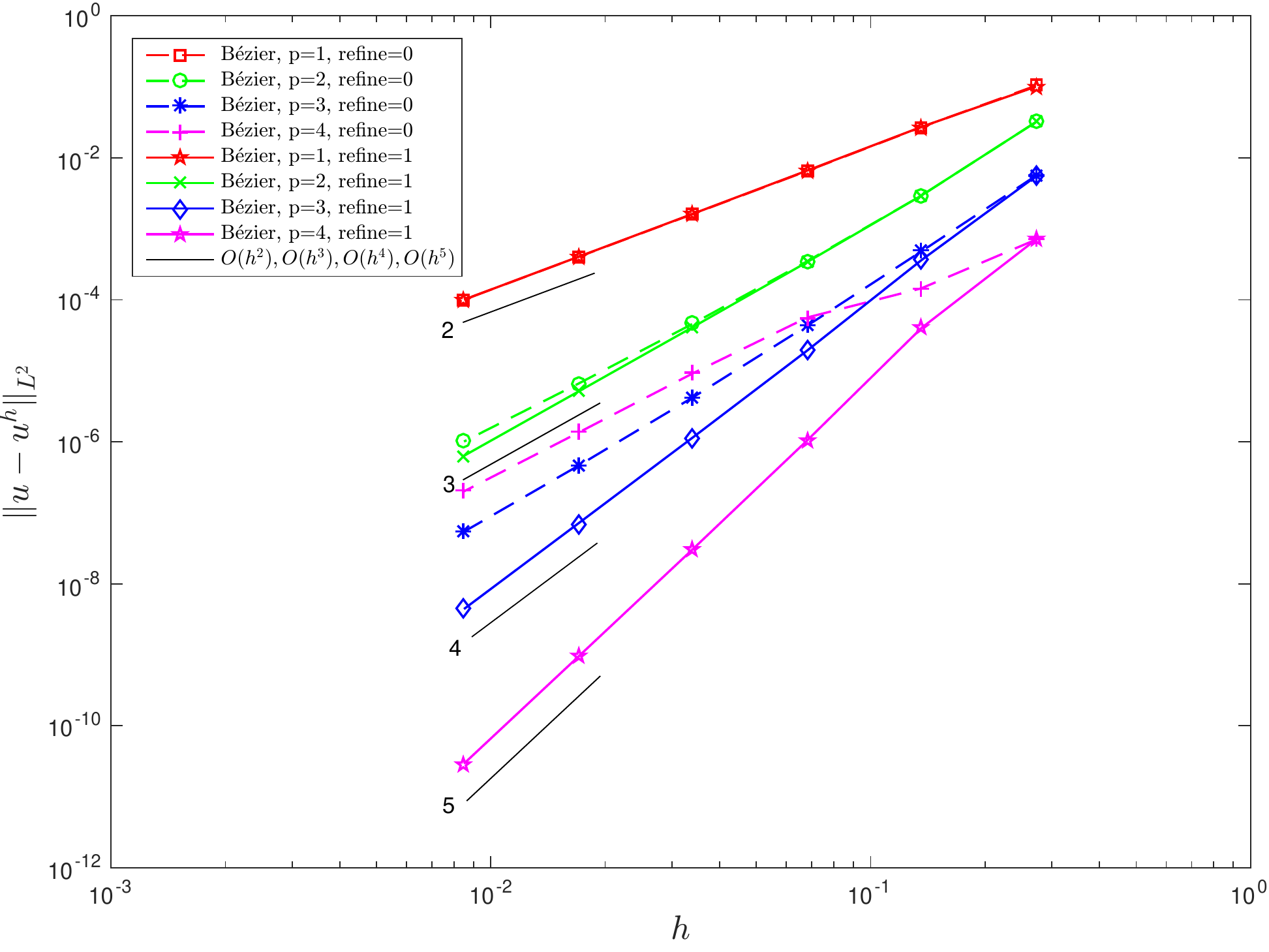} 
\caption{Master/slave mesh ratio, $m:s=2:3$}
\end{subfigure} \\
\begin{subfigure}{0.75\textwidth}
\includegraphics[width=\textwidth]{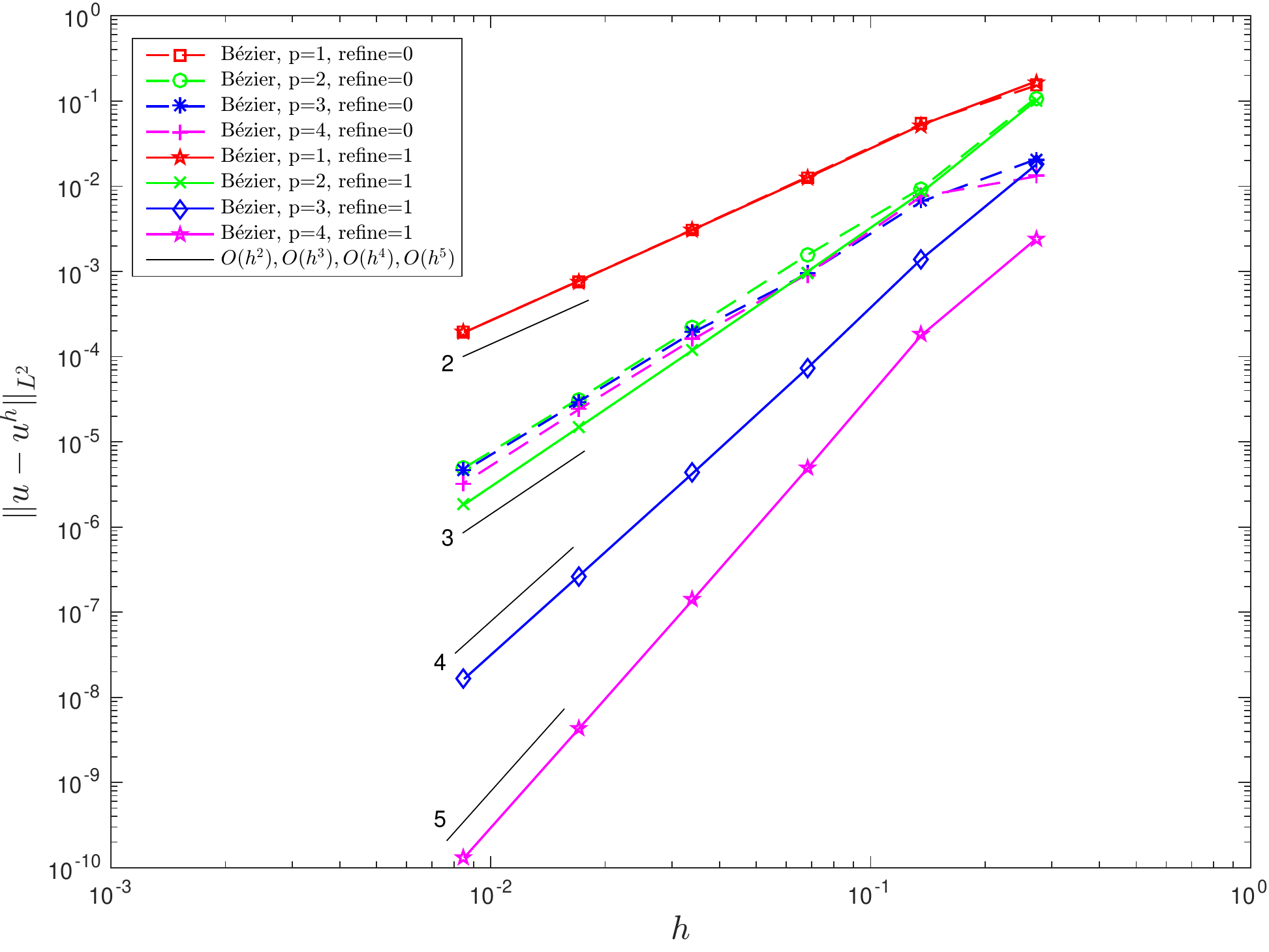}
\caption{Master/slave mesh ratio, $m:s=3:2$}
\end{subfigure}
\caption{Convergence rates for a square domain with two non-conforming patches,
  Dirichlet-Neumann boundary conditions (see
  Figure~\ref{fig:thermalsquareschematic}b),
  matched parameterizations where the dual space is refined $n$ times, $n=0,1$.}	
\label{fig:thermal2DsquareIGSEmatchedNeumannrefine}
\end{figure}

The second boundary condition case, shown in Figure~\ref{fig:thermalsquareschematic}b, allows for a comparison of the proposed
method with the global dual method without crosspoint pollution. As shown in Figure~\ref{fig:thermal2DsquareIGSENeumann}, the optimality of the
global dual mortar method can be observed for this case, while the proposed method
behaves in a similar fashion as in the first full Dirichlet boundary condition case. Again, with one
refinement, the \Bezier dual mortar method obtains optimal rates as shown in Figure~\ref{fig:thermal2DsquareIGSEmatchedNeumannrefine}. This demonstrates that the proposed method is relatively insensitive to crosspoint pollution. This superior behavior can be attributed to the locality of the dual basis. The reduced rates in the proposed method without refinement is due to the lack of higher-order polynomial reproduction in the dual basis.

\begin{figure}[h]
 \centering
\begin{subfigure}{0.75\textwidth}
\includegraphics[width=\textwidth]{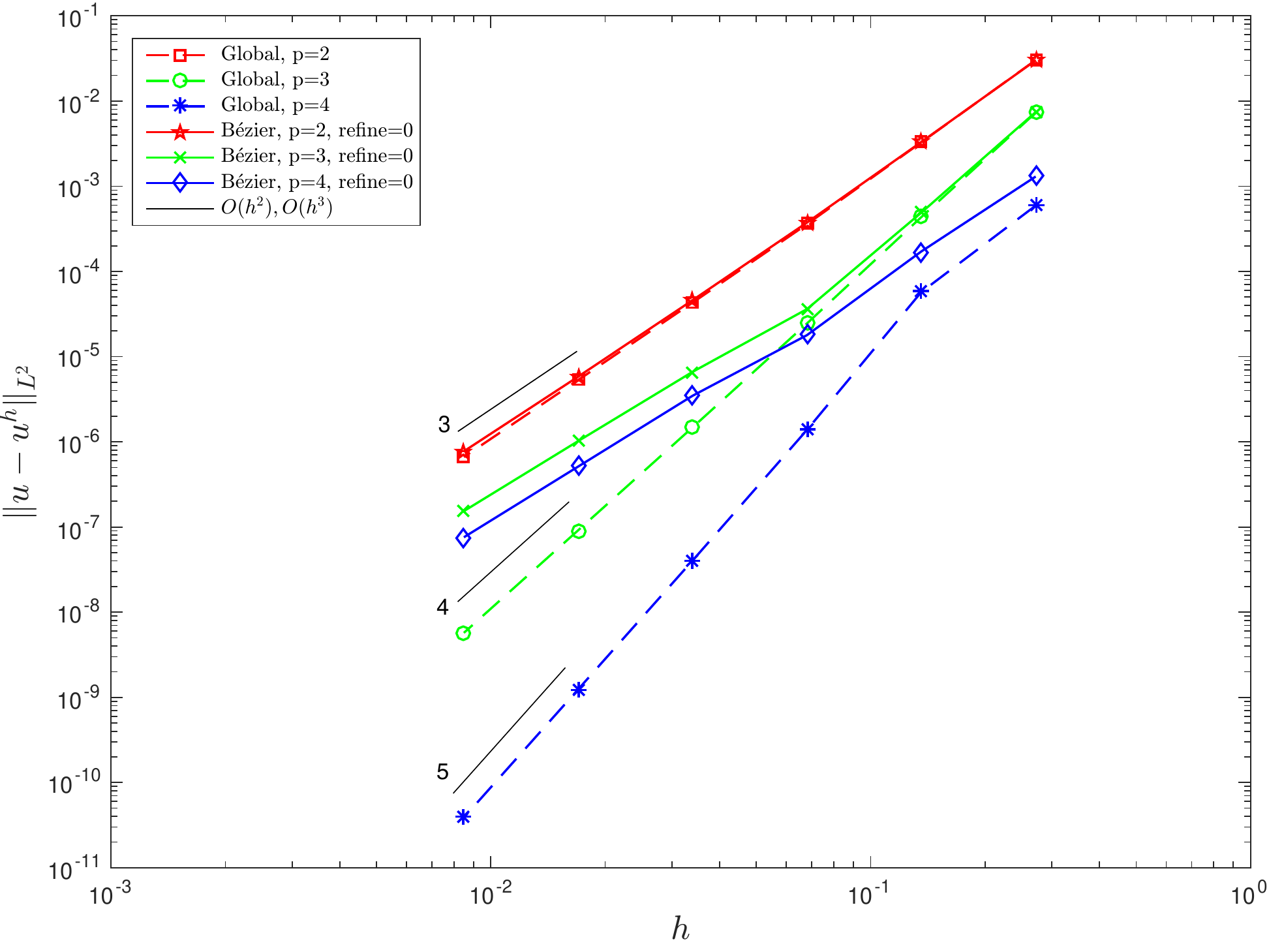} 
\end{subfigure}
\caption{Convergence rates for a square domain with two non-conforming patches, Dirichlet-Neumann boundary conditions (see Figure~\ref{fig:thermalsquareschematic}b) and mismatched parameterizations, master/slave mesh ratio, $m:s=2:3$.}	
\label{fig:thermal2DsquareIGSEmismatchedNeumann}
\end{figure}

\begin{figure}[h]
 \centering
\begin{subfigure}{0.75\textwidth}
\includegraphics[width=\textwidth]{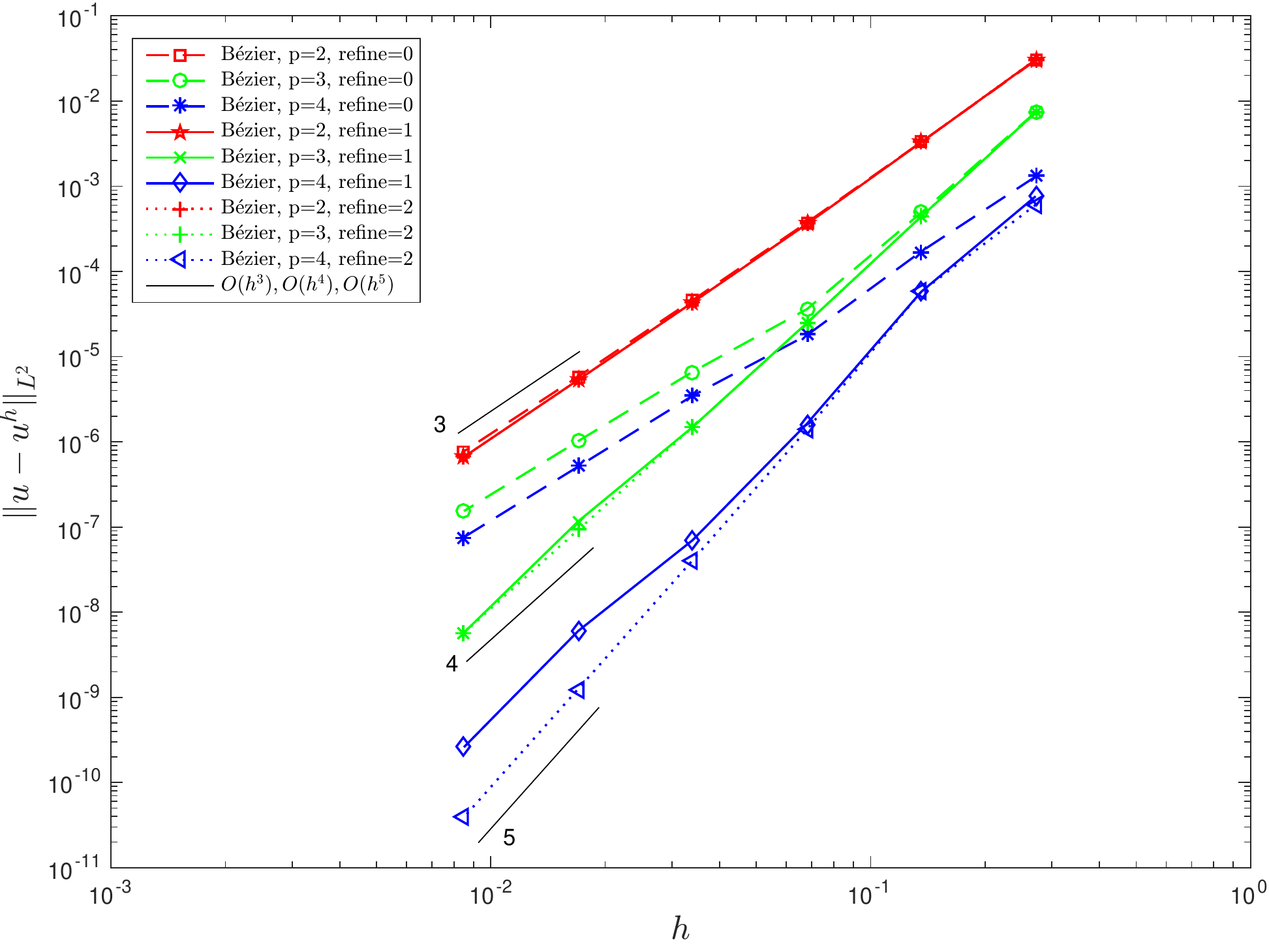} 
\end{subfigure} 
\caption{Convergence rates for a square domain with two non-conforming patches, Dirichlet-Neumann boundary conditions (see Figure~\ref{fig:thermalsquareschematic}b) and mismatched parameterizations where the dual space is refined $n$ times, $n=0,1,2$, master/slave mesh ratio, $m:s=2:3$.}	
\label{fig:thermal2DsquareIGSEmismatchedNeumannrefine}
\end{figure}

We next investigate the impact of mismatched parameterizations on the results.
Mismatched parameterizations can be created by perturbing the position of
the control points along the master and/or slave interfaces. To avoid crosspoint pollution we only consider the Dirichlet-Neumann boundary condition case. The convergence
results for mesh ratio $m:s=2:3$ without refinement are shown
in Figure~\ref{fig:thermal2DsquareIGSEmismatchedNeumann}. Again, the global dual
mortar method obtains the optimal rates for $p=2,3$ and $4$. Without refinement
of the dual basis, the \Bezier dual
mortar method behaves in a manner which is similar to the matched parameterization case.  
To improve solution behavior we refine the dual basis.
 The resulting convergence rates are shown in
 Figure~\ref{fig:thermal2DsquareIGSEmismatchedNeumannrefine}. Since the
geometric mapping is no longer linear, the continuity constraint
(\ref{eq:weakformconstraint}) cannot be imposed exactly by refining the dual
basis once. Therefore, optimal rates cannot be achieved. However, we can
improve the accuracy by simply refining the dual space additional times. As shown in Figure~\ref{fig:thermal2DsquareIGSEmismatchedNeumannrefine}, for $p=2,3$, uniformly refining the dual space once recovers optimal rates and refining twice recovers optimal rate for $p=4$. Recall that regardless of how many times the dual space is refined the number of global degrees-of-freedom remains fixed.

\subsection{A manufactured solution on an annular domain}
 \label{test2}

\begin{figure}[h]
  \centering
  \begin{overpic} [scale=0.48]{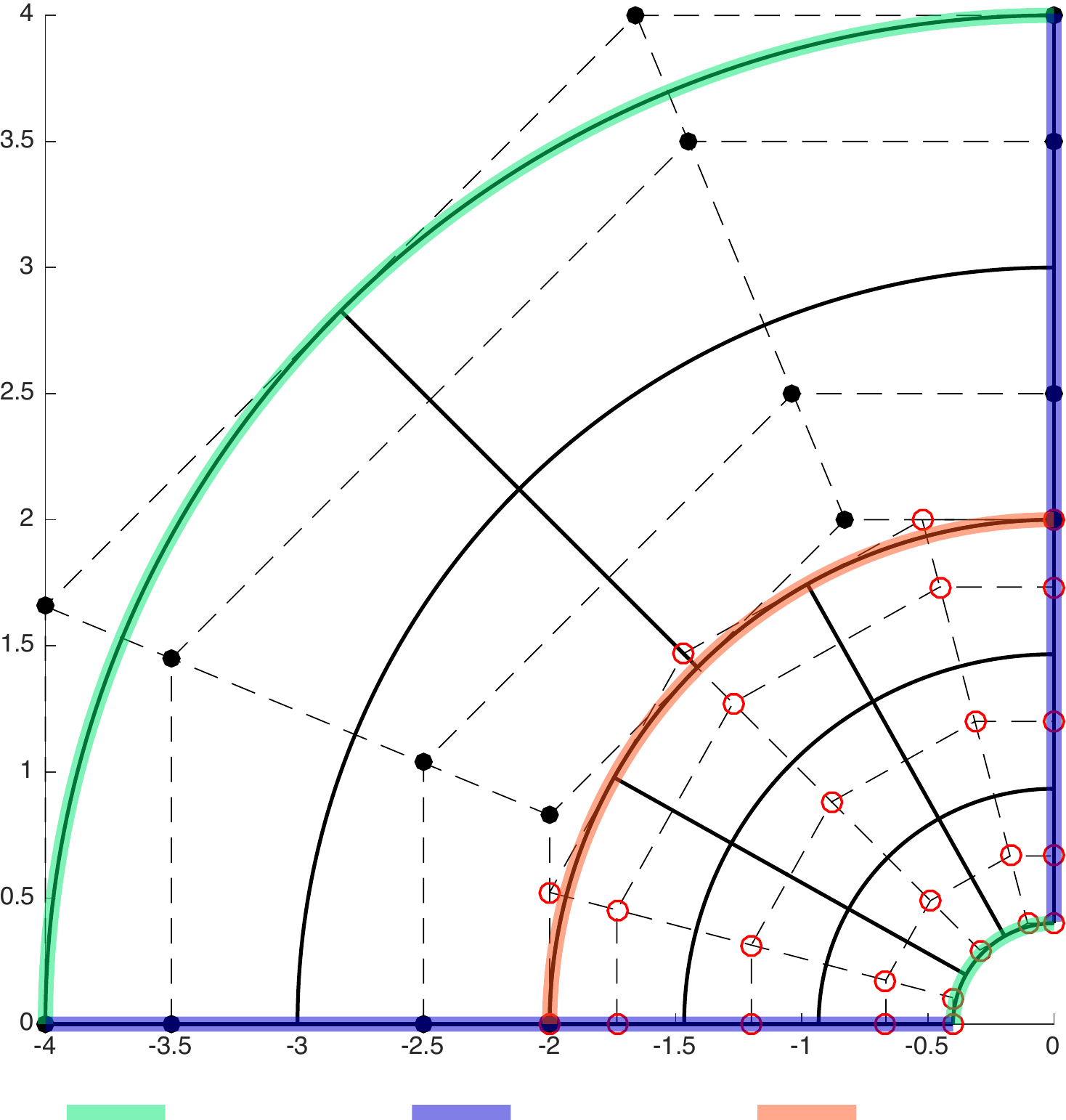}
    \put( 15, 0 ){\scalebox{0.7}{Dirichlet}}
    \put( 46, 0 ){\scalebox{0.7}{Neumann}}
    \put( 78, 0 ){\scalebox{0.7}{Interface}}
    \put( 86, 92 ){\scalebox{1}{$\Omega^m$}}
    \put( 96, 21 ){\scalebox{1}{$\Omega^s$}}
    \end{overpic}
\caption{An annular domain composed of two quadratic nonconforming NURBS patches.}
\label{fig:initial_mesh_circular}
\end{figure}

We now solve the Poisson equation, $-\Delta u =  f$, on the annular domain,
$\Omega= \big \{  (r,\phi) \, | \, 0.4 \leq r \leq 4, \allowbreak \pi/2\leq
\phi\leq \pi \big \}$ . This example tests the effectiveness of rational dual
basis functions. The domain is composed of two NURBS patches as shown in
Figure~\ref{fig:initial_mesh_circular}. The internal force and the boundary
conditions correspond to the manufactured solution, $u\left(x,y\right)= \sin(\pi
x)\sin(  \pi y)$. As shown in Figure~\ref{fig:initial_mesh_circular}, there are
no crosspoints in this problem. Note that we only consider matched parameterizations in this example.

\begin{figure}
  \centering
\begin{subfigure}{0.65\textwidth}
\includegraphics[width=\textwidth]{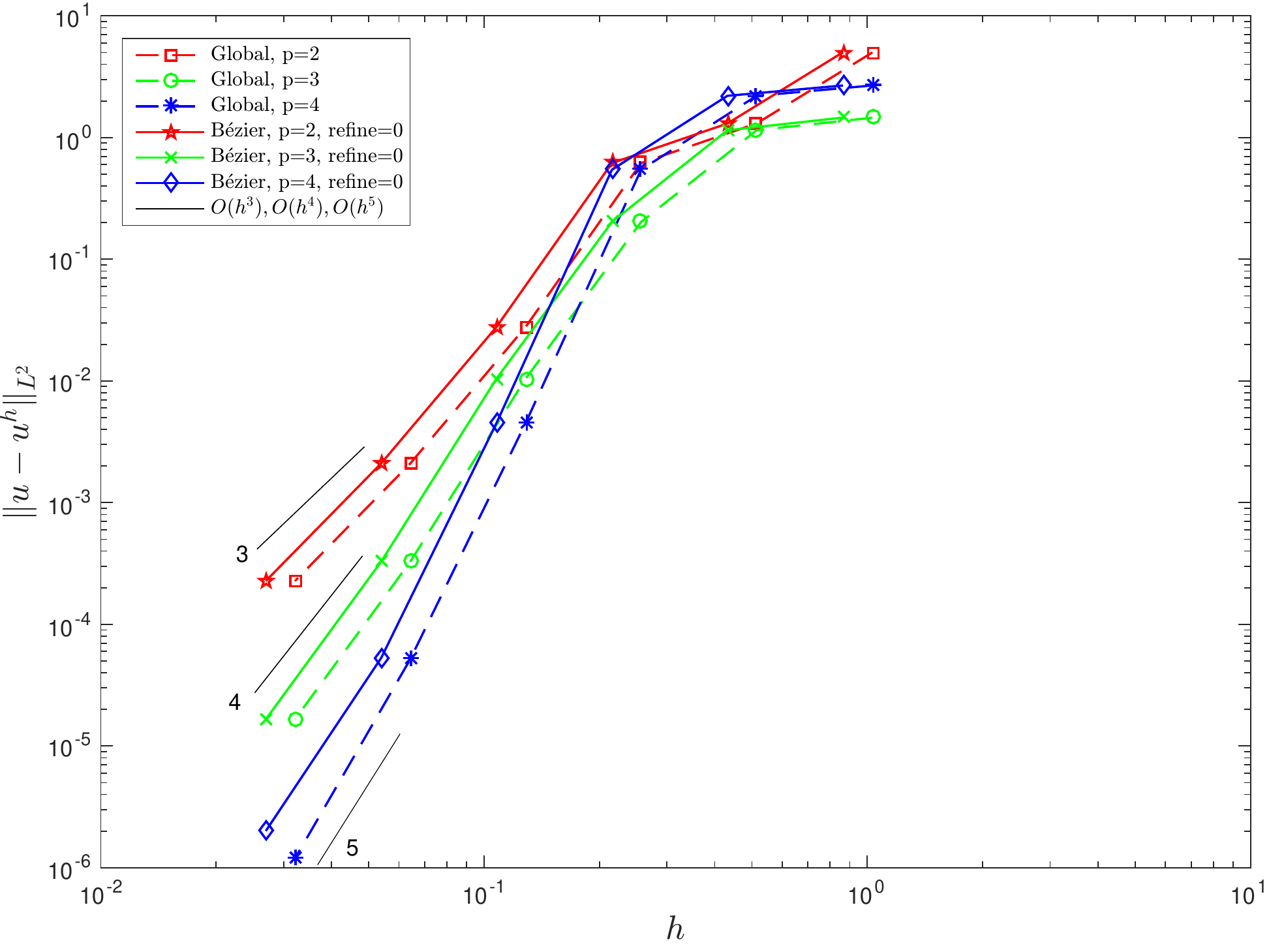} 
\caption{Master/slave mesh ratio, $m:s = 2:3$}
\end{subfigure} \\
\begin{subfigure}{0.65\textwidth}
\includegraphics[width=\textwidth]{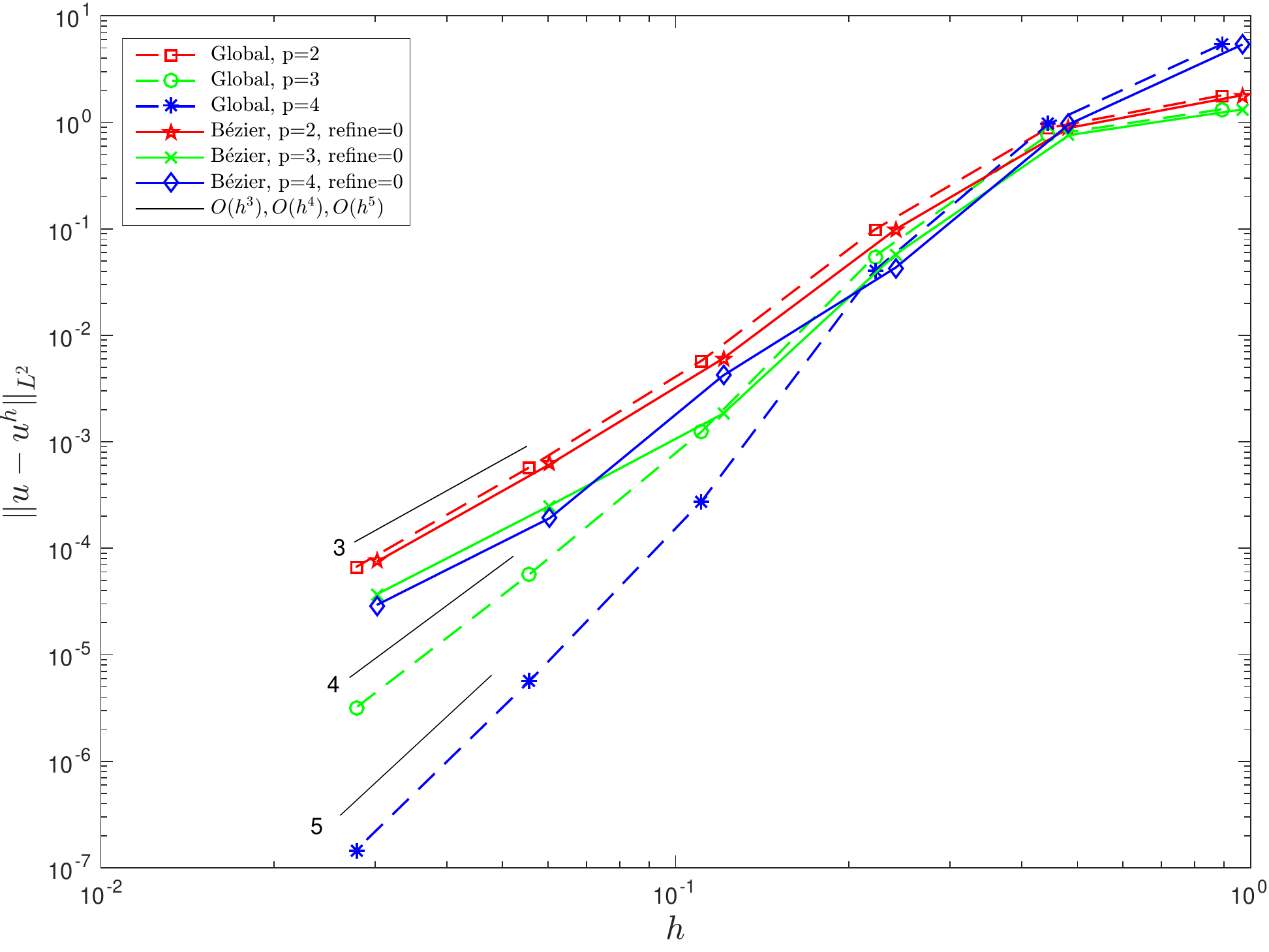}
\caption{Master/slave mesh ratio, $m:s = 3:2$}
\end{subfigure}
\caption{Convergence rates for an annular domain with two non-conforming NURBS
  patches, Dirichlet-Neumann boundary conditions (see Figure~\ref{fig:initial_mesh_circular}) and  matched parameterizations.}
\label{fig:thermalcircularms23ms32ratesNeumann}
\end{figure}

\begin{figure}
  \centering
\begin{subfigure}{0.65\textwidth}
\includegraphics[width=\textwidth]{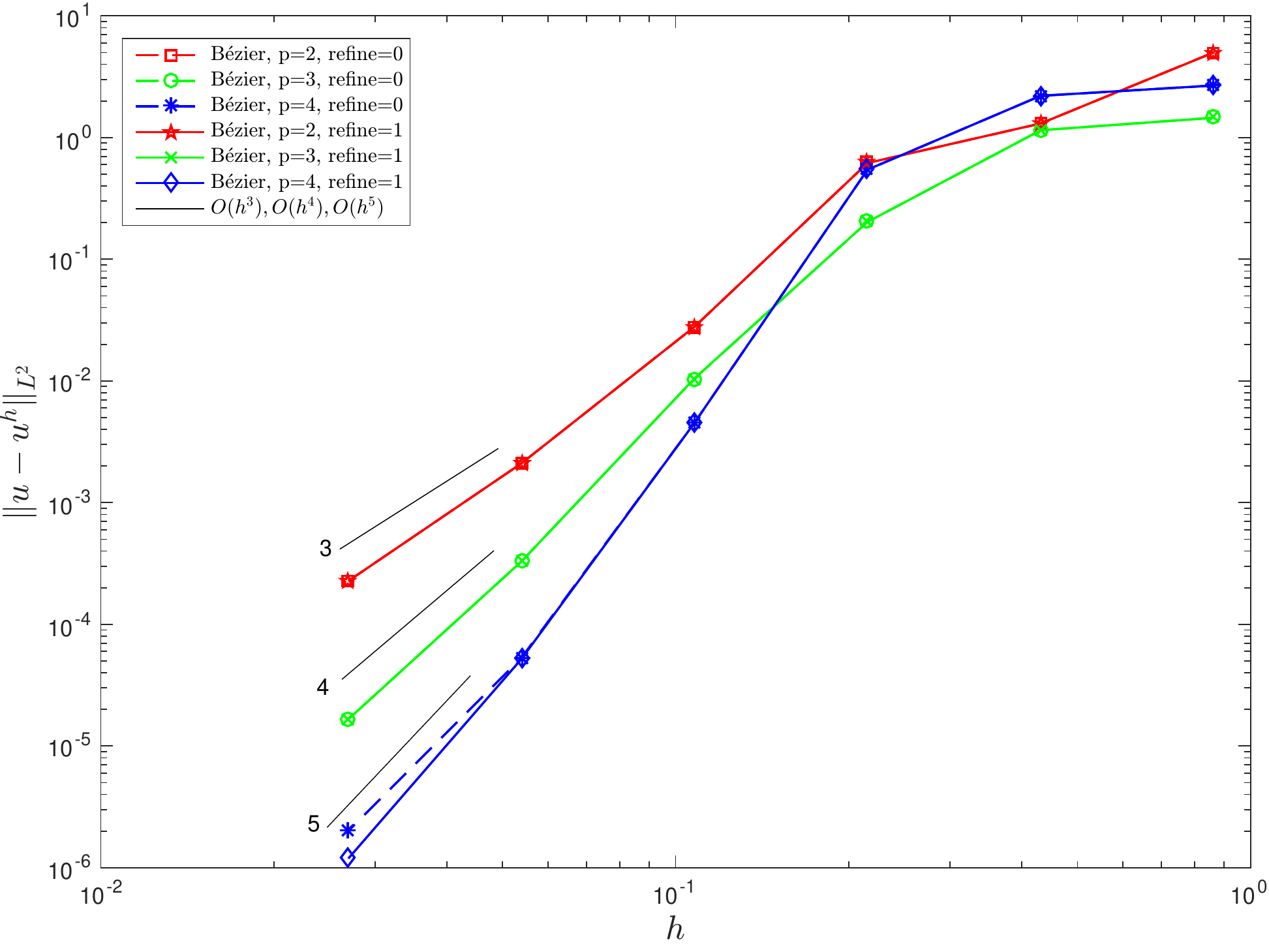} 
\caption{Master/slave mesh ratio, $m:s = 2:3$}
\end{subfigure} \\
\begin{subfigure}{0.65\textwidth}
\includegraphics[width=\textwidth]{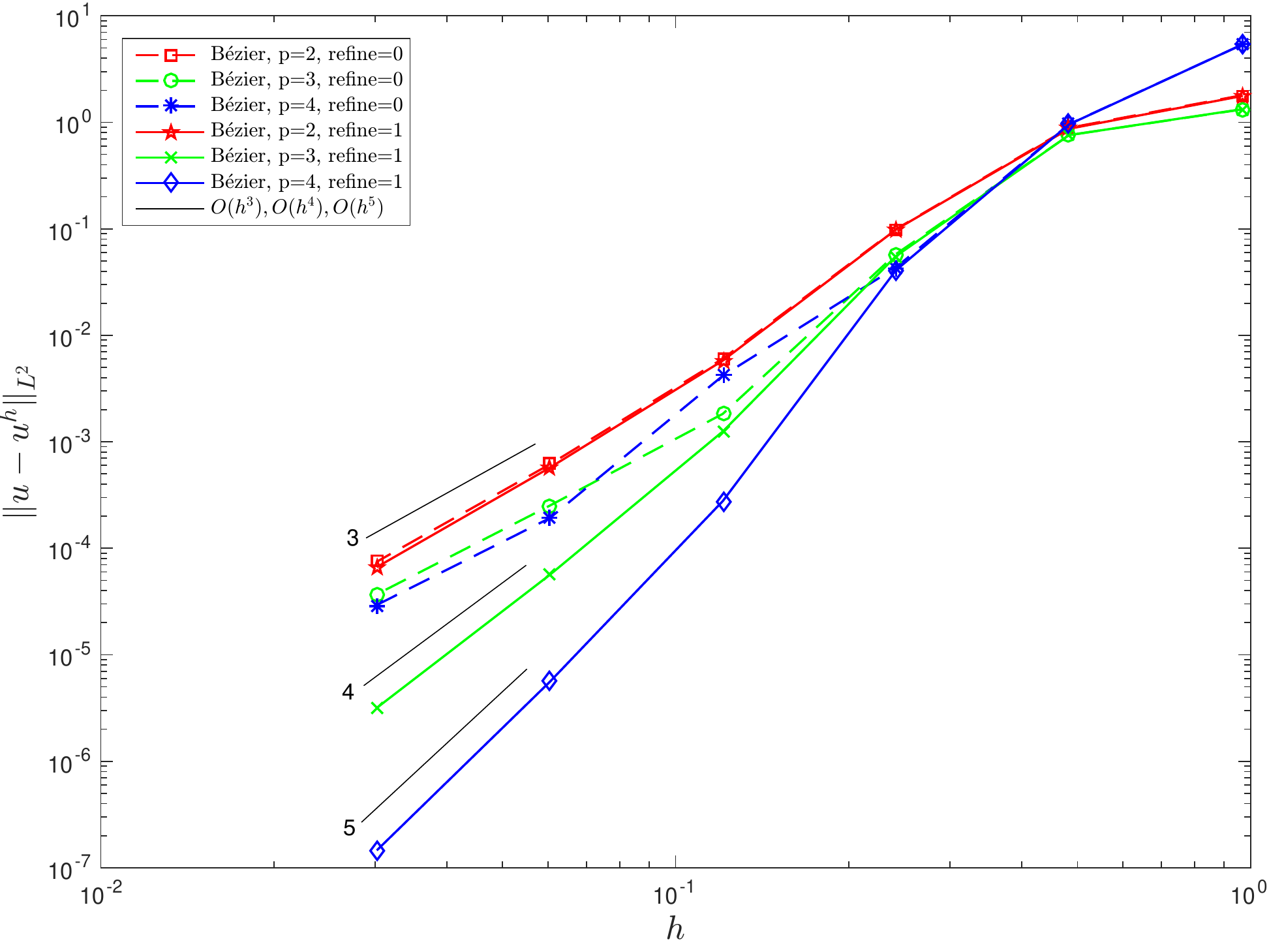}
\caption{Master/slave mesh ratio, $m:s = 3:2$}
\end{subfigure}
\caption{Convergence rates for an annular domain with two non-conforming NURBS
  patches, Dirichlet-Neumann boundary conditions (see
  Figure~\ref{fig:initial_mesh_circular}) and refined matched
  parameterizations.}
\label{fig:thermalcircularms23ms32refinedratesNeumann}
\end{figure}

The convergence rates in the $L^2$-norm of the displacement are shown in Figure
\ref{fig:thermalcircularms23ms32ratesNeumann} for $p=2,3,4$, without refining
the dual basis. As can be seen, the global dual mortar method achieves the optimal
rates for all degrees and mesh ratios $m:s=2:3$ and $m:s=3:2$. The proposed method
achieves the optimal rates for $p=2,3$ and $m:s=2:3$, and slightly deteriorated
convergence rate for $p=4$.  For mesh ratio $m:s=2:3$, the proposed method
experiences reduced convergence rates for $p=3,4$. However, one refinement
recovers optimal rates for $p=4$ with mesh ratio $m:s = 2:3$, and $p=3$ and
$4$ with mesh ratio $m:s = 3:2$, as shown in Figure~\ref{fig:thermalcircularms23ms32refinedratesNeumann}.

\FloatBarrier
\subsection{Infinite elastic plate with a circular hole} \label{subsec:infiniteplate}
 \label{test3}
 \begin{figure}[h]
  \centering
 \includegraphics[width=0.5\textwidth]{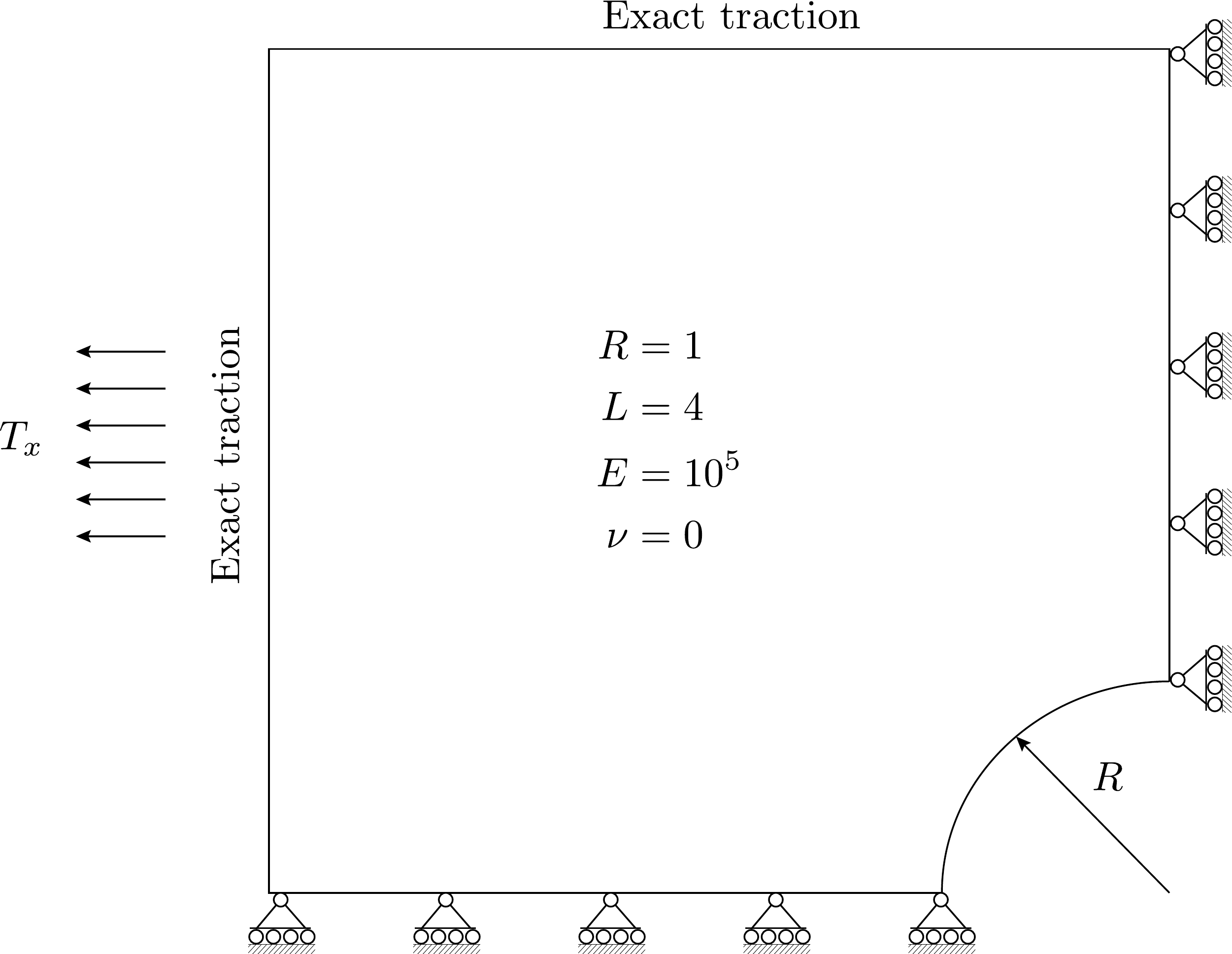} 
 \caption{A schematic for the infinite elastic plate with a circular hole benchmark.}	
 \label{fig:infiniteplateschematics}
 \end{figure}

We next simulate the classical infinite elastic plate with a circular hole benchmark problem. In this case, we apply a constant traction in the $x$-direction at infinity. Due to symmetry, only one quarter of the plate is modeled as shown in Figure~\ref{fig:infiniteplateschematics}, where $T_x$ is the traction, $R$ is the radius of the hole, $L$ is the length of each side of the plate, $E$ is Young's modulus, and $\nu$ is Poisson's ratio. An analytical solution to this problem can be found in~\cite{Cottrell:2009rp} and is reproduced here for completeness 
\begin{align}
  \sigma_{rr}(r,\theta) &= \frac{T_x}{2} \left ( 1 - \frac{ R^2 }{ r^2 } \right ) + \frac{T_x}{2} \left ( 1 - 4\frac{R^2}{r^2} + 3\frac{R^4}{r^4} \right ) \text{cos} 2\theta,\\
  \sigma_{\theta \theta}(r,\theta) &= \frac{T_x}{2} \left ( 1 + \frac{ R^2 }{ r^2 } \right ) - \frac{T_x}{2} \left ( 1 + \frac{3}{4}\frac{R^4}{r^4} \right ) \text{cos} 2\theta,\\
  \sigma_{r\theta}(r,\theta) &= - \frac{T_x}{2} \left ( 1 + 2\frac{R^2}{r^2} - 3\frac{R^4}{r^4} \right ) \text{sin} 2\theta.
\end{align}

 \begin{figure}[h]
  \centering
  \begin{subfigure}{0.48\textwidth}
  \begin{overpic}[scale=0.46]{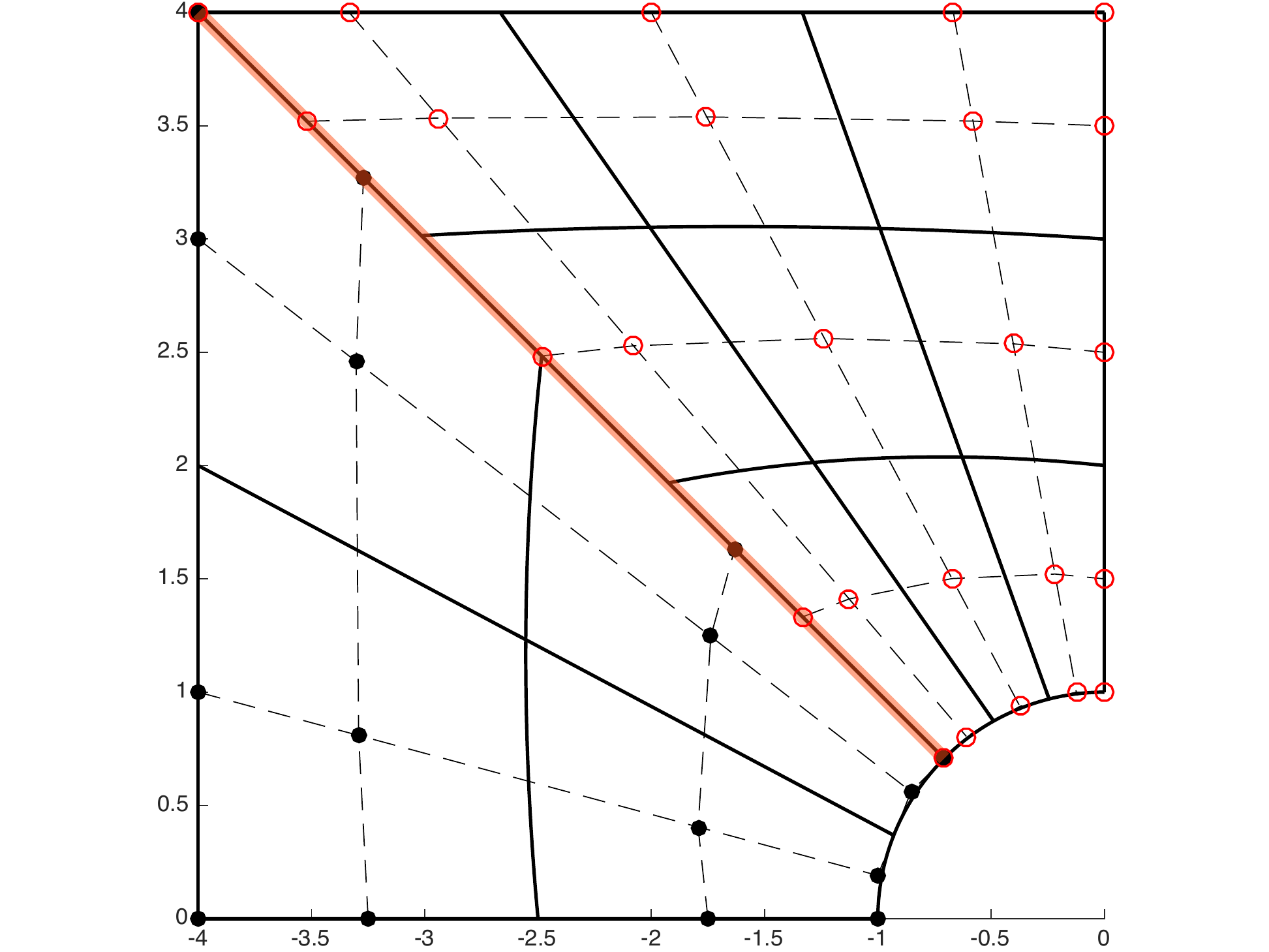} 
    \put(17, 4){$\Omega^m$}
    \put(80, 68){$\Omega^s$}
    \end{overpic}
 \caption{Matched parameterizations}
 \end{subfigure}
 \begin{subfigure}{0.48\textwidth}
   \begin{overpic}[scale=0.46]{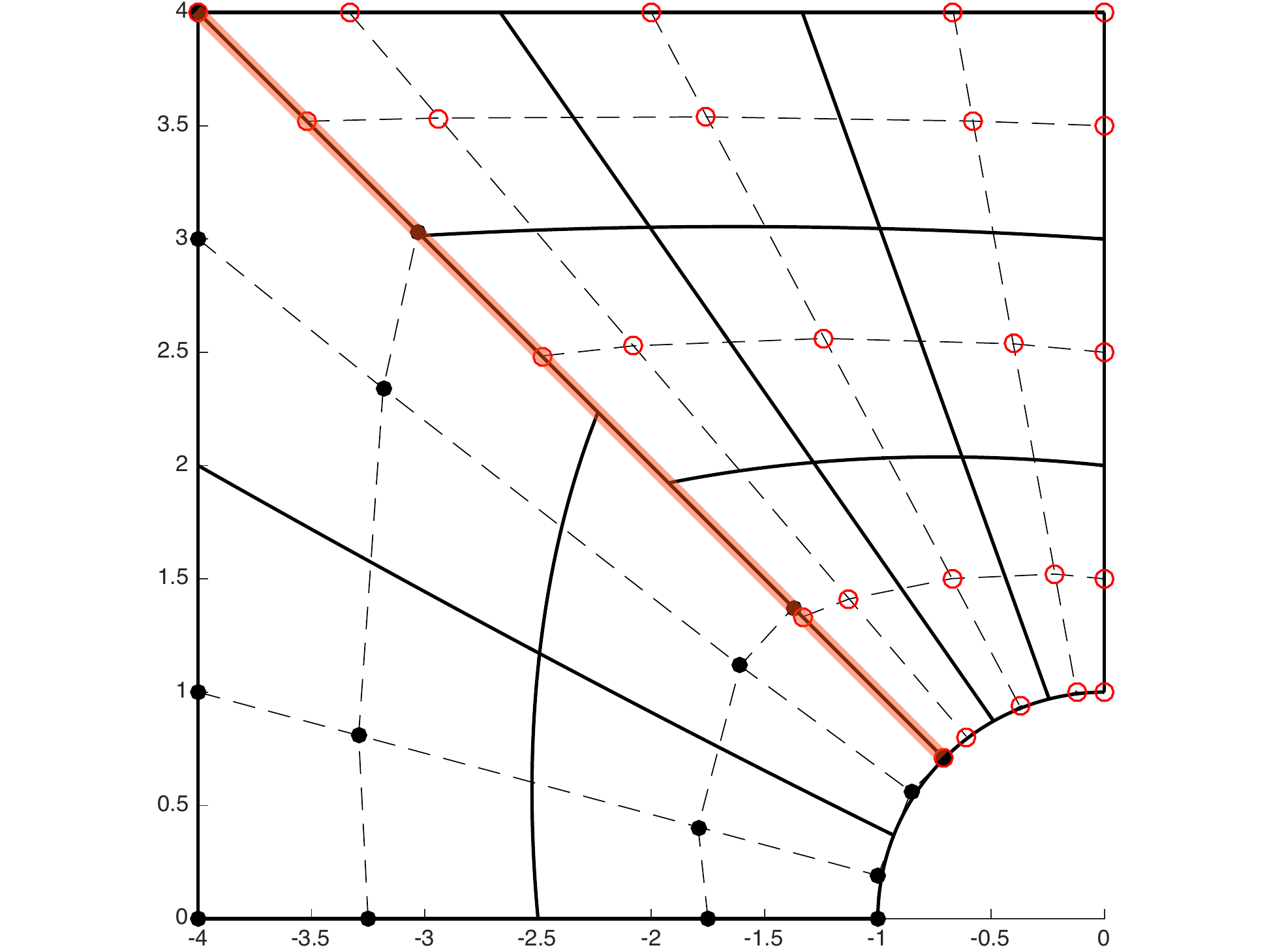} 
    \put(17, 4){$\Omega^m$}
    \put(80, 68){$\Omega^s$}
 \end{overpic}
 \caption{Mismatched parameterizations}
 \end{subfigure}
 \caption{NURBS meshes for a quarter plate with a hole.}	
 \label{fig:solidinfiniteplatepatch2schematic}
 \end{figure}
 
\begin{figure}[h]
 \centering
  \begin{subfigure}{0.75\textwidth}
 \includegraphics[width=\textwidth]{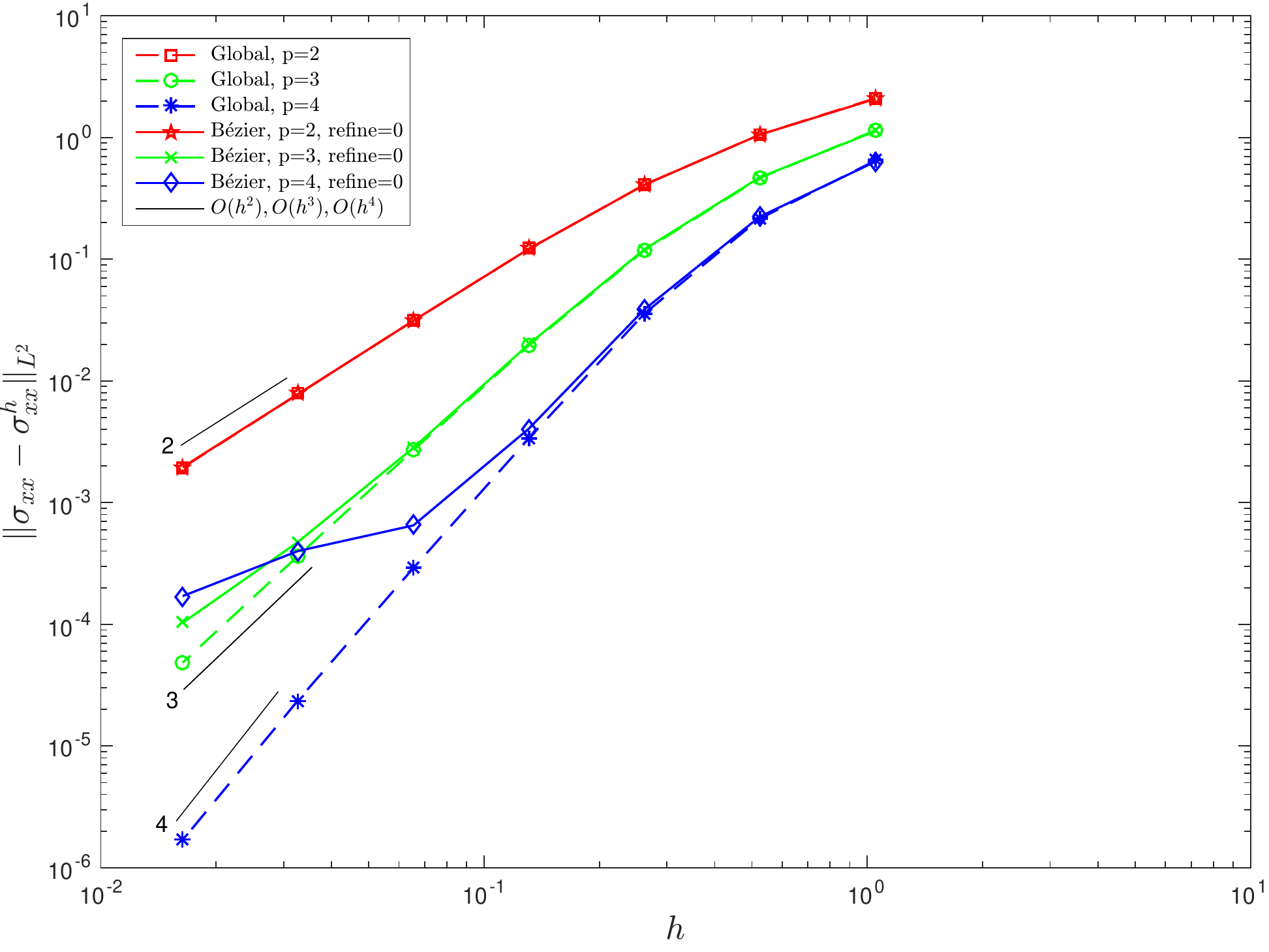} 
 \caption{Master/slave mesh ratio, $m:s=2:3$}
 \end{subfigure} \\
 \begin{subfigure}{0.75\textwidth}
 \includegraphics[width=\textwidth]{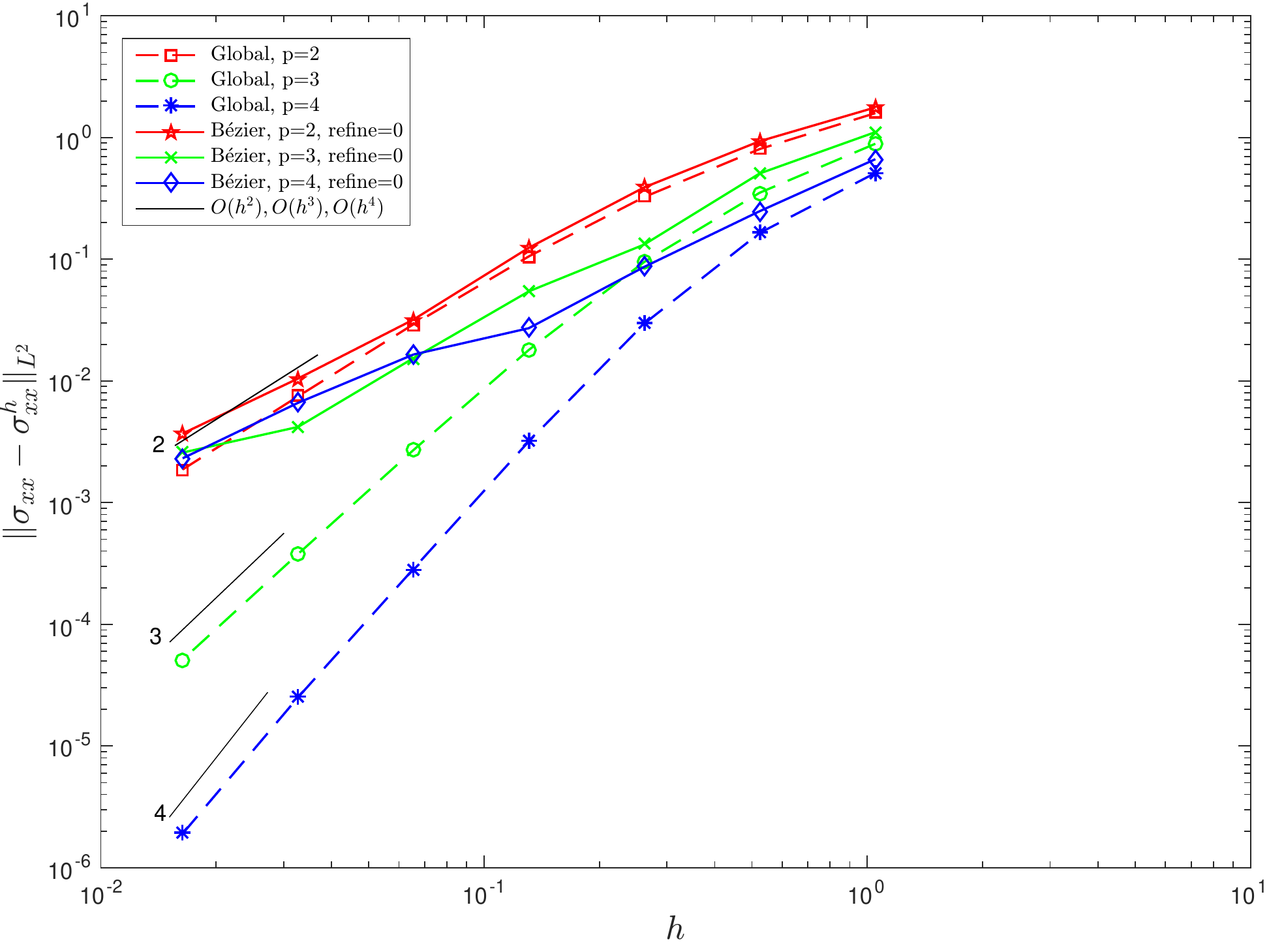} 
 \caption{Master/slave mesh ratio, $m:s=3:2$}
 \end{subfigure}
\caption{Stress convergence rates for a quarter plate with a circular hole
  decomposed into two nonconforming NURBS patches with matched parameterizations
  (see Figure~\ref{fig:solidinfiniteplatepatch2schematic}a).}	
\label{fig:solid2DinfiniteplateIGSE}
\end{figure}

\begin{figure}[h]
 \centering
  \begin{subfigure}{0.75\textwidth}
 \includegraphics[width=\textwidth]{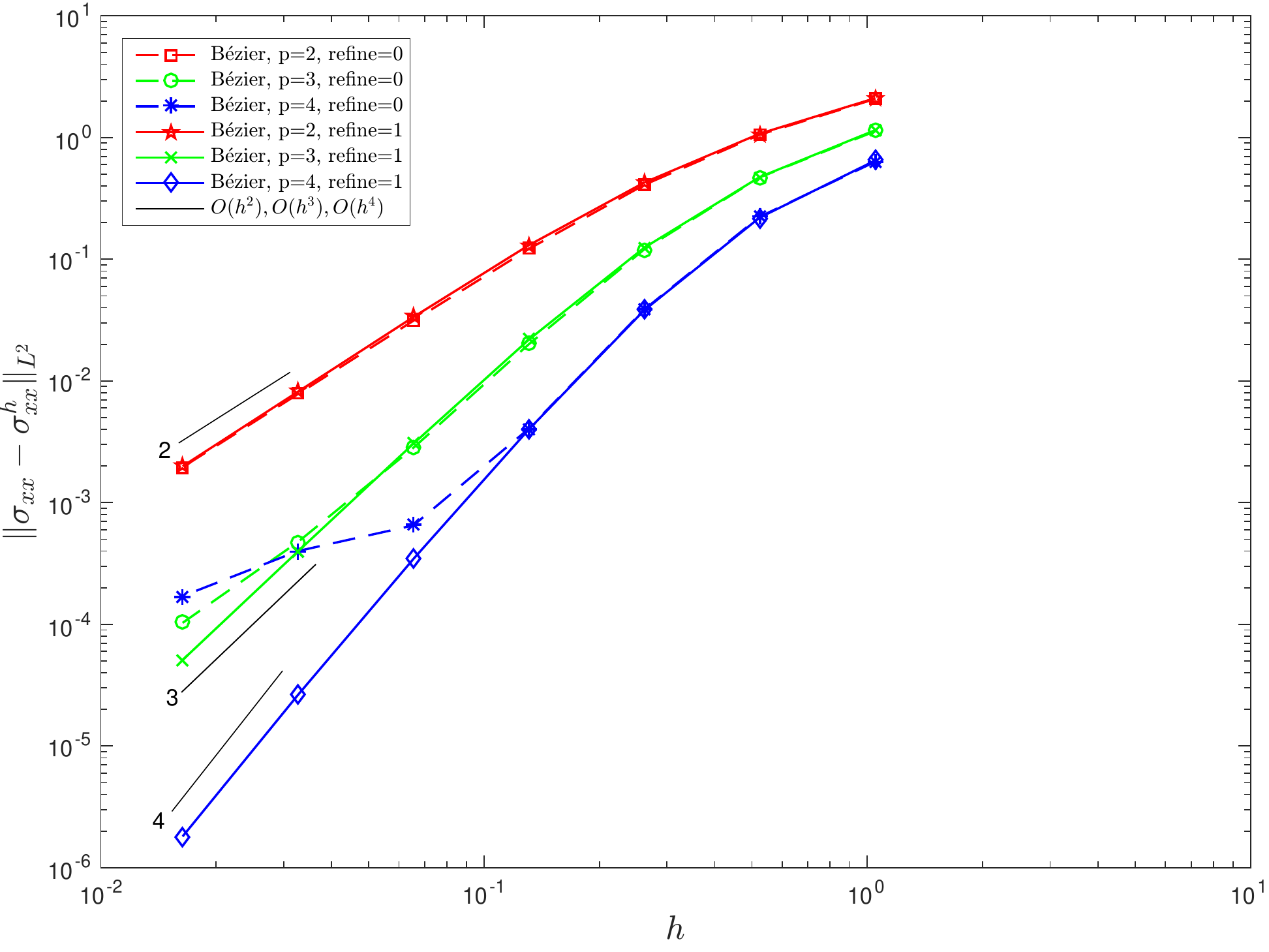} 
 \caption{Master/slave mesh ratio, $m:s=2:3$}
 \end{subfigure} \\
 \begin{subfigure}{0.75\textwidth}
 \includegraphics[width=\textwidth]{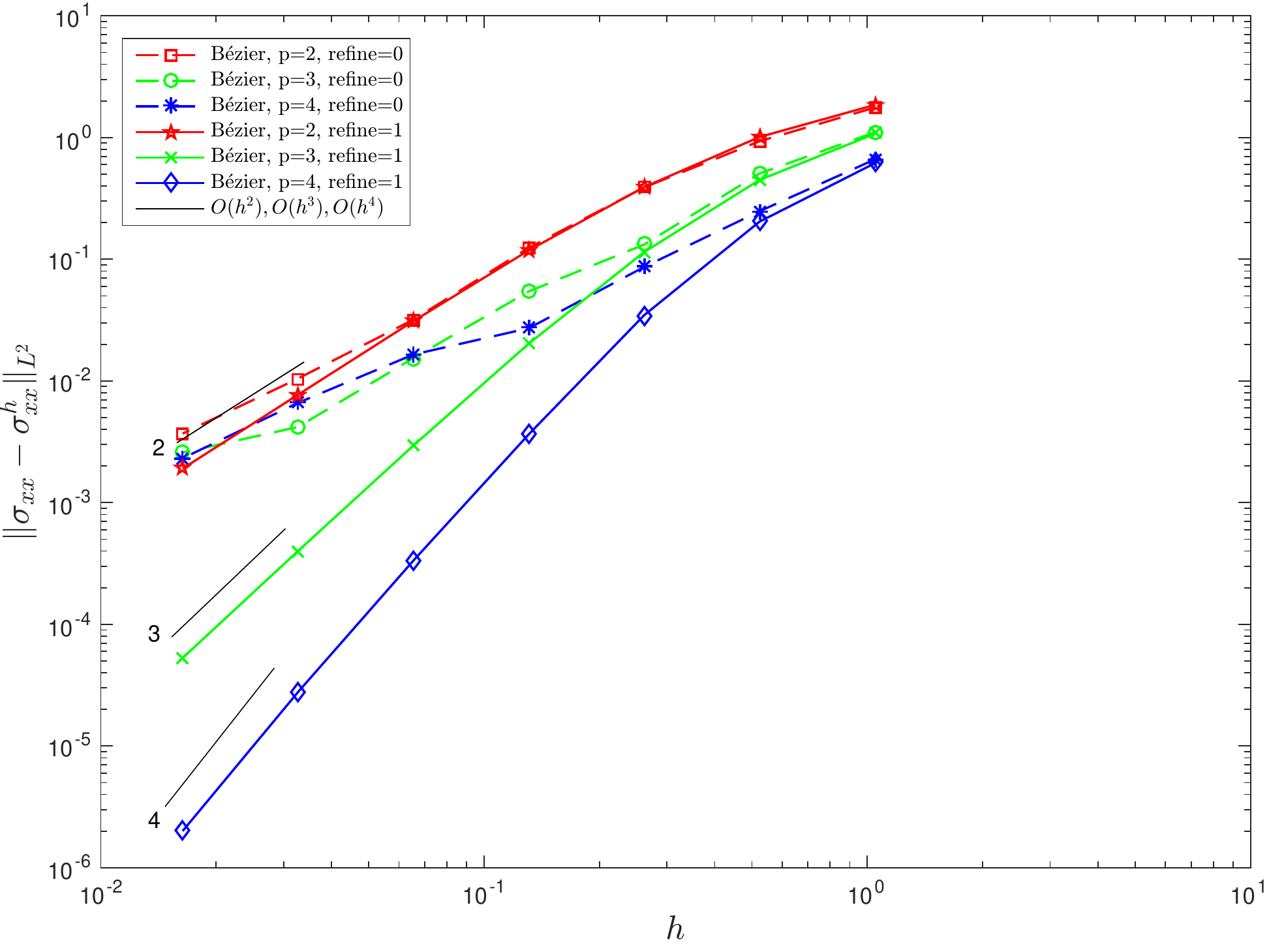} 
 \caption{Master/slave mesh ratio, $m:s=3:2$}
 \end{subfigure}
\caption{Stress convergence rates for a quarter plate with a circular hole
  decomposed into two nonconforming NURBS patches with refined matched
  parameterizations (see Figure~\ref{fig:solidinfiniteplatepatch2schematic}a).}
\label{fig:solid2DinfiniteplateIGSEmatchedrefine}
\end{figure}

\begin{figure}[h]
 \centering
  \begin{subfigure}{0.75\textwidth}
 \includegraphics[width=\textwidth]{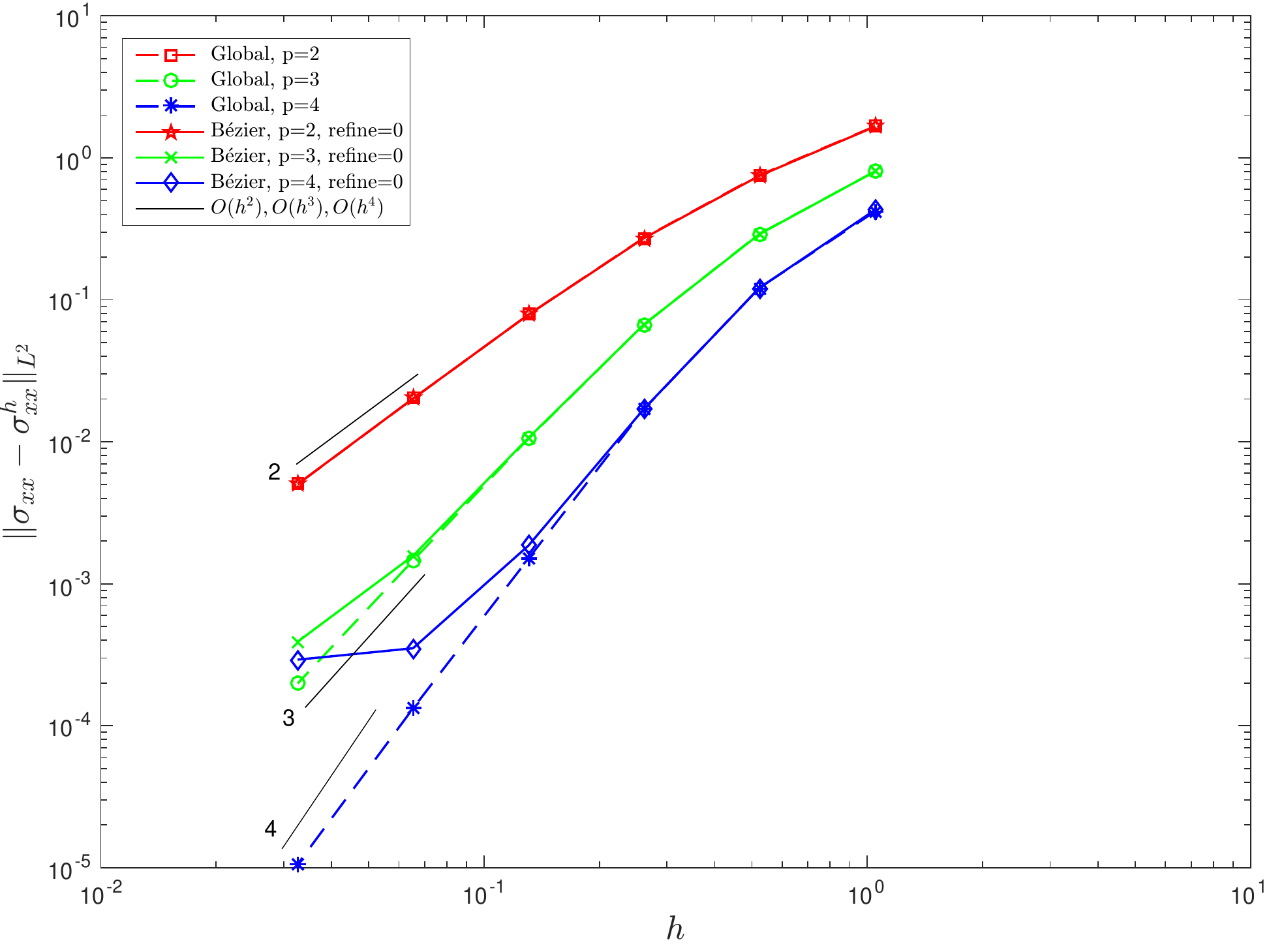} 
 \caption{Master/slave mesh ratio, $m:s=2:3$}
 \end{subfigure} \\
 \begin{subfigure}{0.75\textwidth}
 \includegraphics[width=\textwidth]{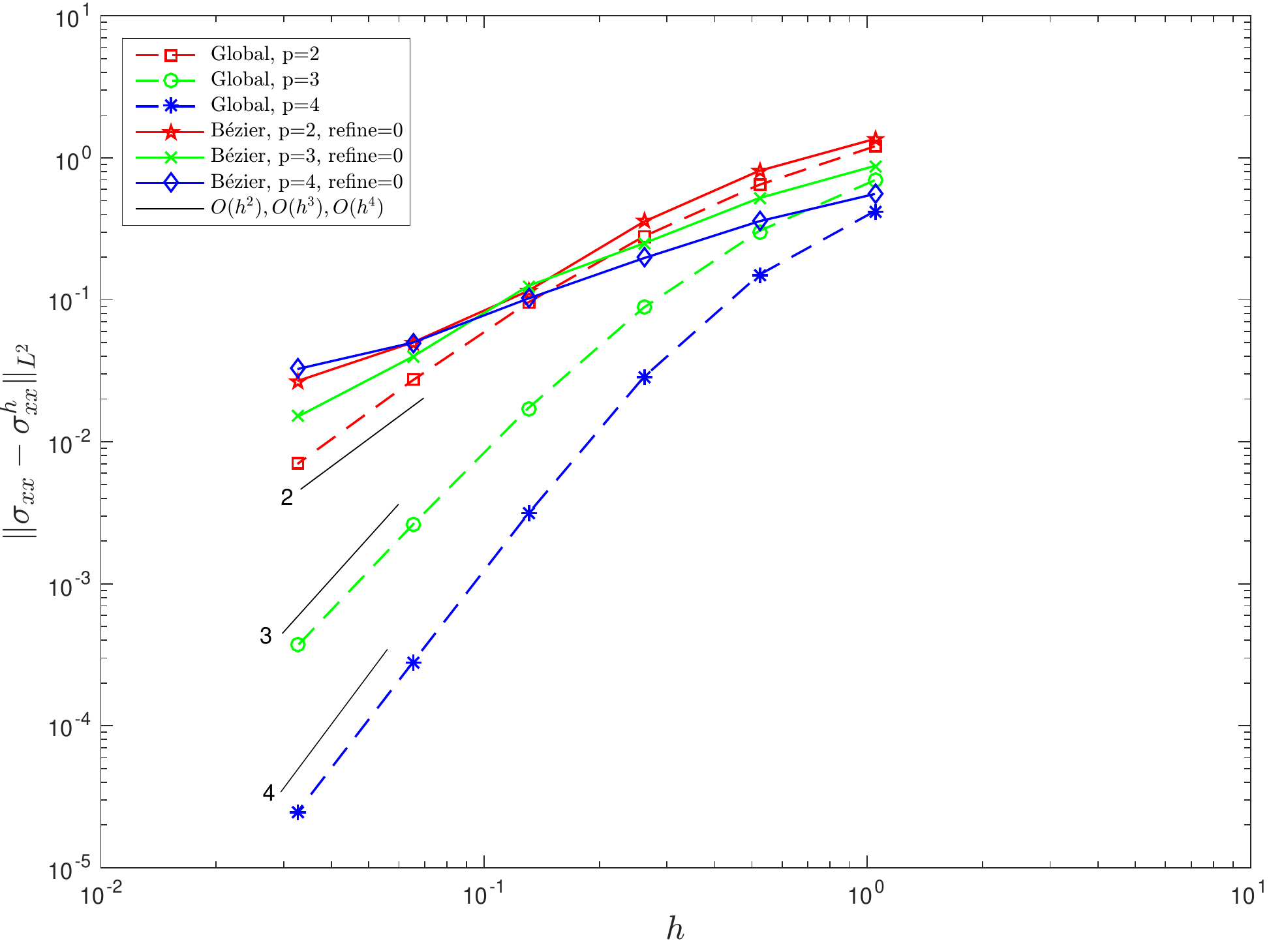} 
 \caption{Master/slave mesh ratio, $m:s=3:2$}
 \end{subfigure}
\caption{Stress convergence rates for a quarter
  plate with a circular hole decomposed into two nonconforming NURBS patches
  with mismatched parameterizations (see Figure~\ref{fig:solidinfiniteplatepatch2schematic}b).}
\label{fig:solid2DinfiniteplateIGSEmismatched}
\end{figure}

\begin{figure}[h]
 \centering
  \begin{subfigure}{0.75\textwidth}
 \includegraphics[width=\textwidth]{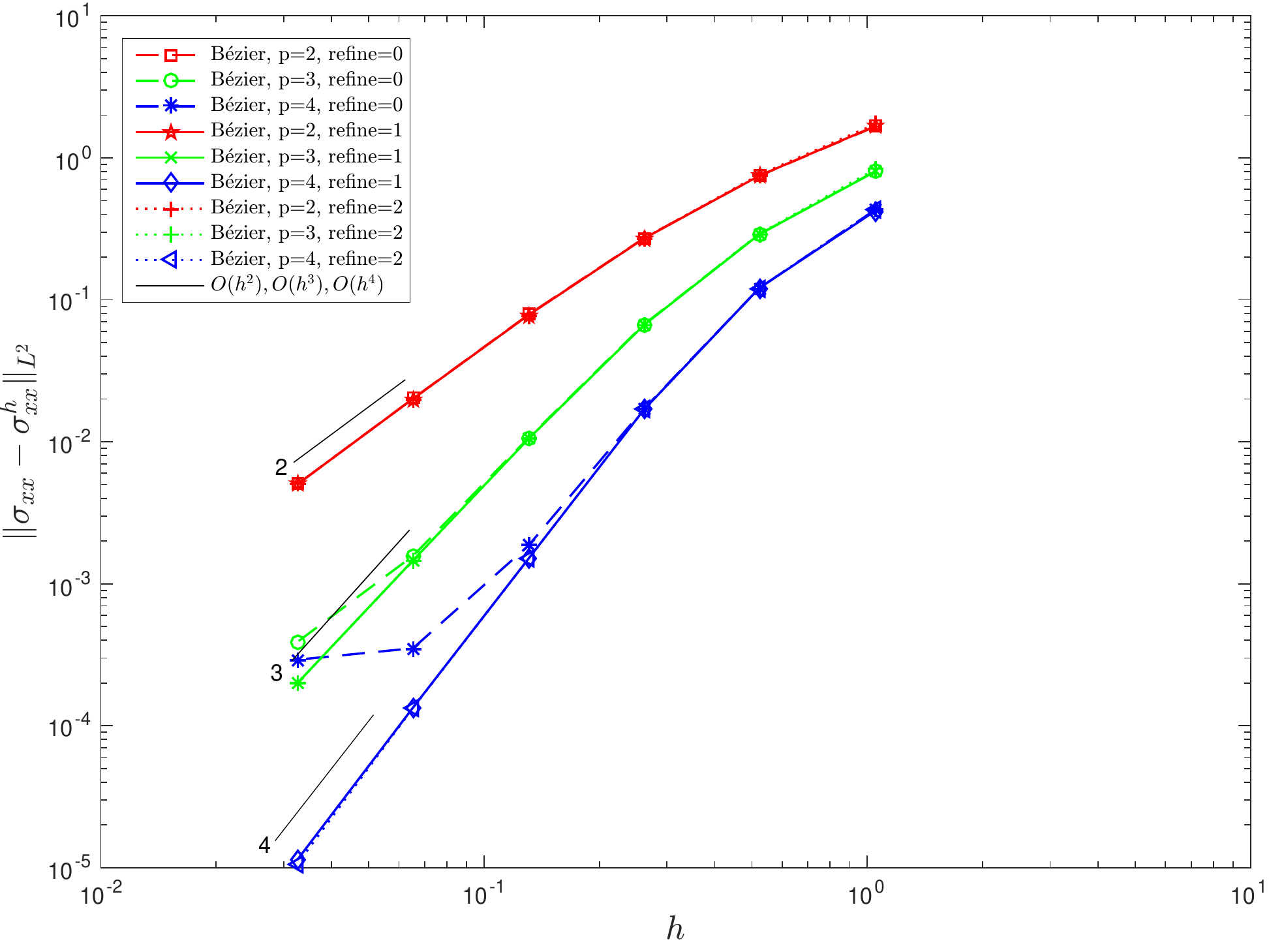} 
 \caption{Master/slave mesh ratio, $m:s=2:3$}
 \end{subfigure} \\
 \begin{subfigure}{0.75\textwidth}
 \includegraphics[width=\textwidth]{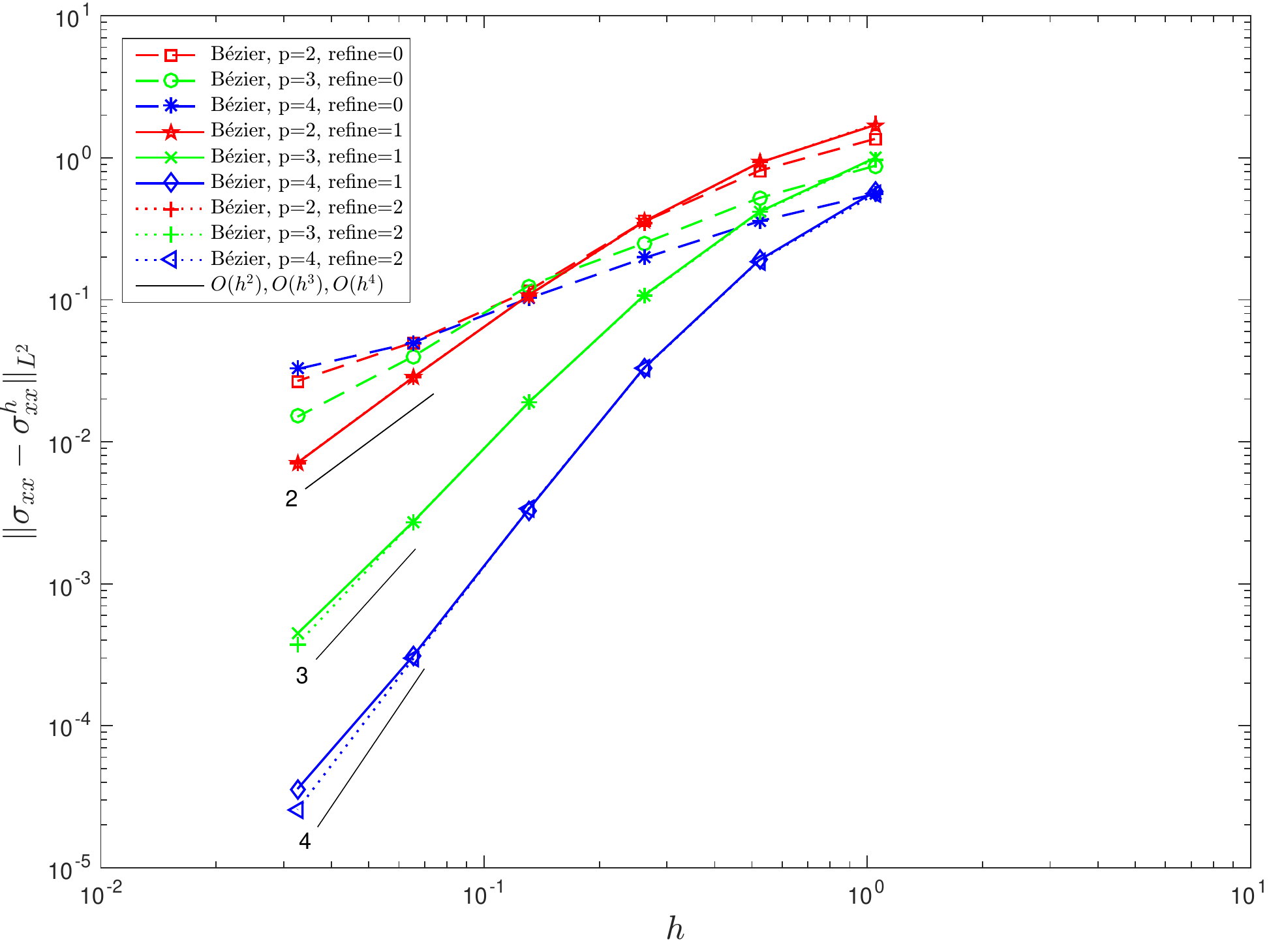} 
 \caption{Master/slave mesh ratio, $m:s=3:2$}
 \end{subfigure}
\caption{Stress convergence rates for a quarter plate with a circular hole
  decomposed into two nonconforming NURBS patches with refined mismatched
  parameterizations (see Figure~\ref{fig:solidinfiniteplatepatch2schematic}b).}
\label{fig:solid2DinfiniteplateIGSEmismatchedrefine}
\end{figure}

 As shown in Figures~\ref{fig:solidinfiniteplatepatch2schematic}a and b, we first
 decompose the geometry into two patches with matched and mismatched
 parameterizations, respectively. The convergence
 rates of the stress component $\sigma_{xx}$ in the $L^2$-norm are
 optimal for the global dual mortar method for all cases in
 Figure~\ref{fig:solid2DinfiniteplateIGSE} and
 \ref{fig:solid2DinfiniteplateIGSEmismatched} due to the absence of crosspoints.
 For the matched parameterization case without refinement, the \Bezier dual mortar method only achieves optimal rate for $p=2$, $m:s = 2:3$. With one refinement of the dual basis, the proposed method recovers optimal rates for all cases as
 shown in Figure~\ref{fig:solid2DinfiniteplateIGSEmatchedrefine}. For mismatched parameterizations without refinement, the \Bezier dual
 mortar method exhibits similar reduced convergence rates as for the matched parameterization case as shown in Figure~\ref{fig:solid2DinfiniteplateIGSEmismatched}.
 However, as shown in Figure~\ref{fig:solid2DinfiniteplateIGSEmismatchedrefine}a,
 after one refinement the proposed method recovers the optimal rates for all
 degrees,
 $m:s = 2:3$, and refining the dual basis twice results in optimal rates
 for $p=2,3$, $m:s = 3:2$, as shown in Figure~\ref{fig:solid2DinfiniteplateIGSEmismatchedrefine}b.

\begin{figure}[h]
  \centering
  \begin{subfigure}{0.48\textwidth}
    \begin{overpic}[scale=0.46]{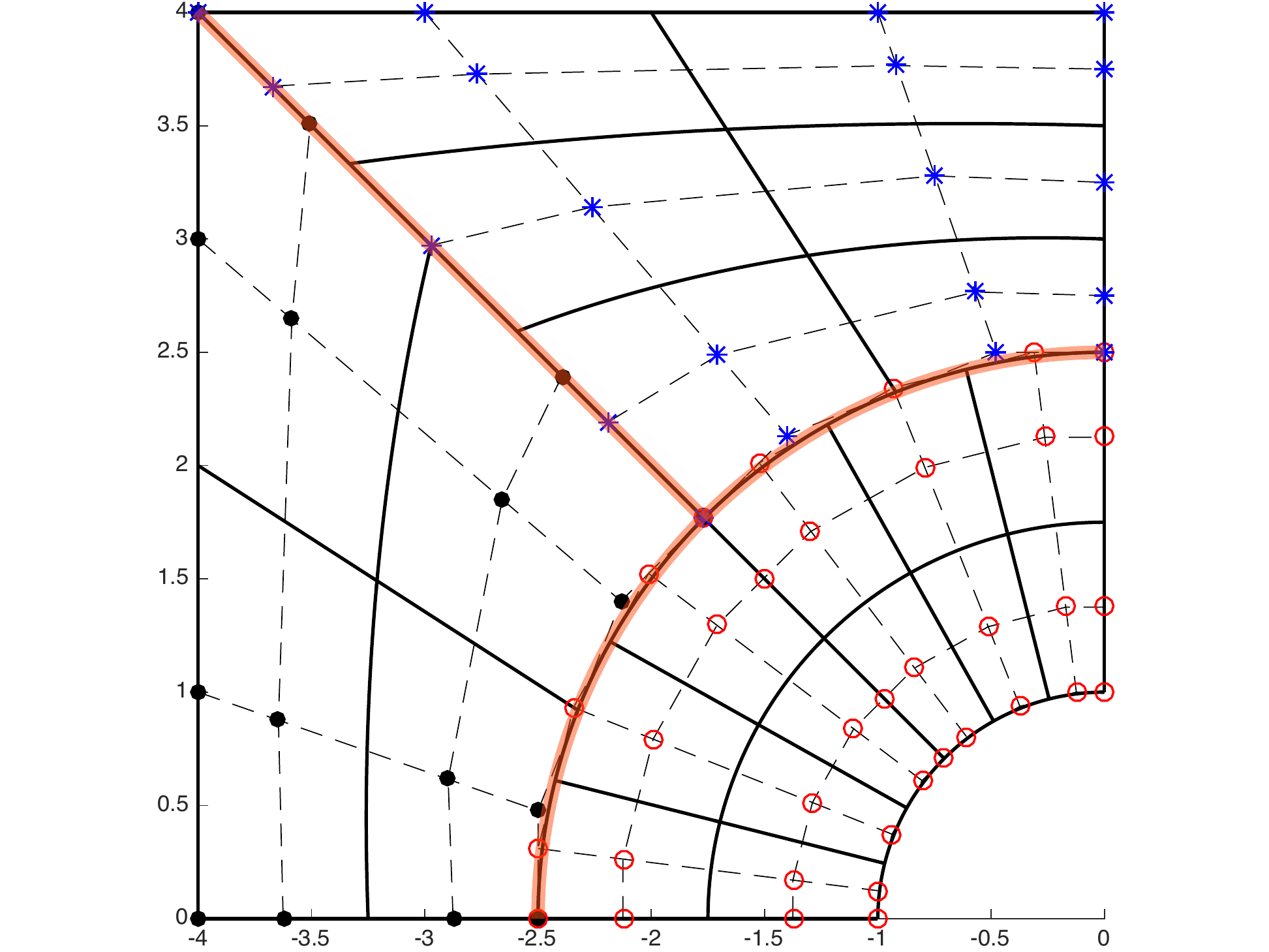}
    \put(17, 4){$\Omega^m$}
    \put(73, 68){$\Omega^s/$}
    \put(76, 10){$\Omega^m$}
    \put(80.5, 68){$\Omega^m$}
      \end{overpic}
 \caption{Matched parameterizations}
 \end{subfigure}
 \begin{subfigure}{0.48\textwidth}
   \begin{overpic}[scale=0.46]{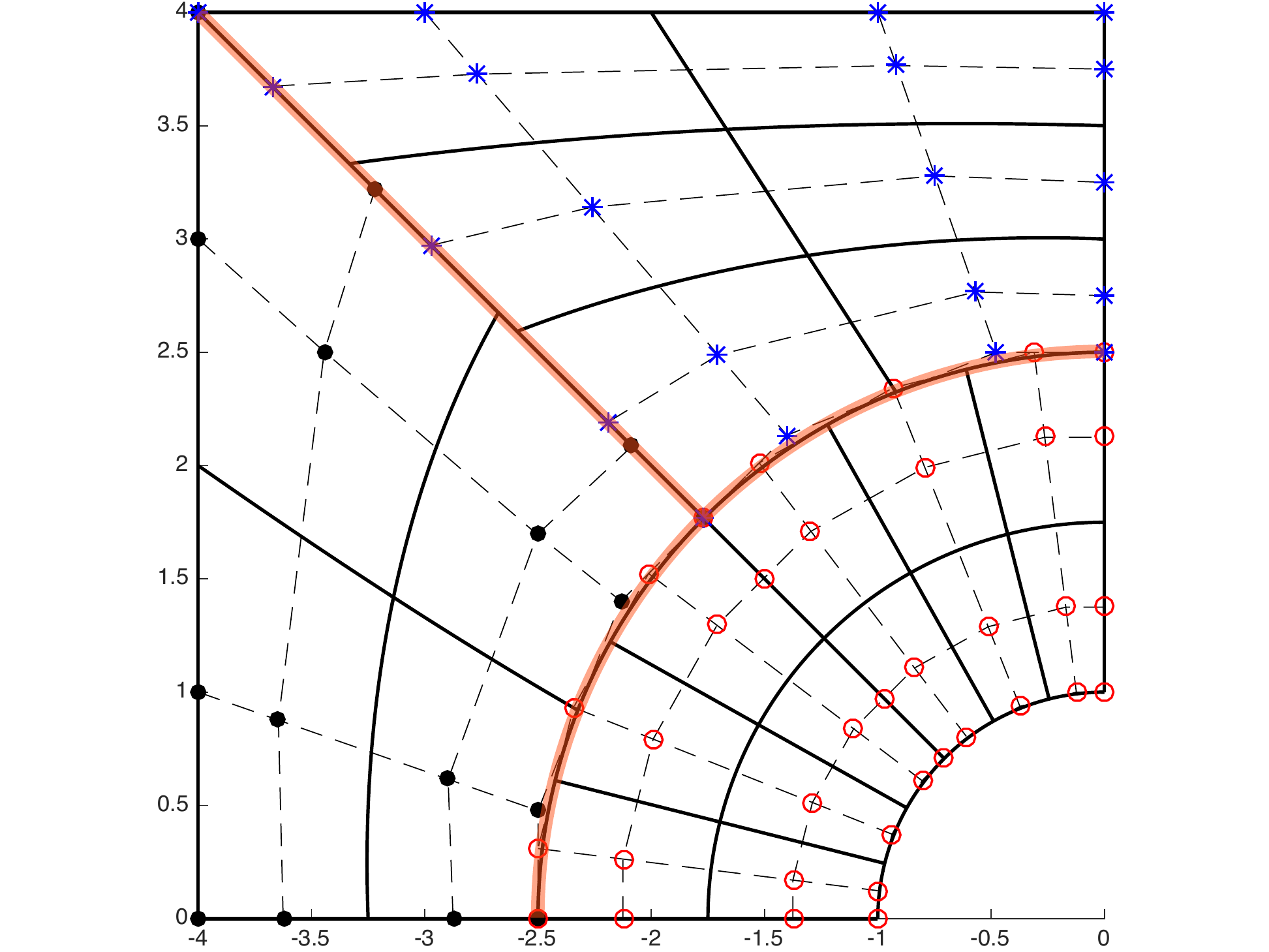} 
    \put(17, 4){$\Omega^m$}
    \put(73, 68){$\Omega^s/$}
    \put(76, 10){$\Omega^m$}
    \put(80.5, 68){$\Omega^m$}
     \end{overpic}
   \caption{Mismatched parameterizations}
   \end{subfigure}
 \caption{NURBS meshes for a quarter plate with a hole.}
 \label{fig:solidinfiniteplatepatch3schematic}
 \end{figure}

\begin{figure}[h]
\centering
\begin{subfigure}{0.75\textwidth}
\includegraphics[width=\textwidth]{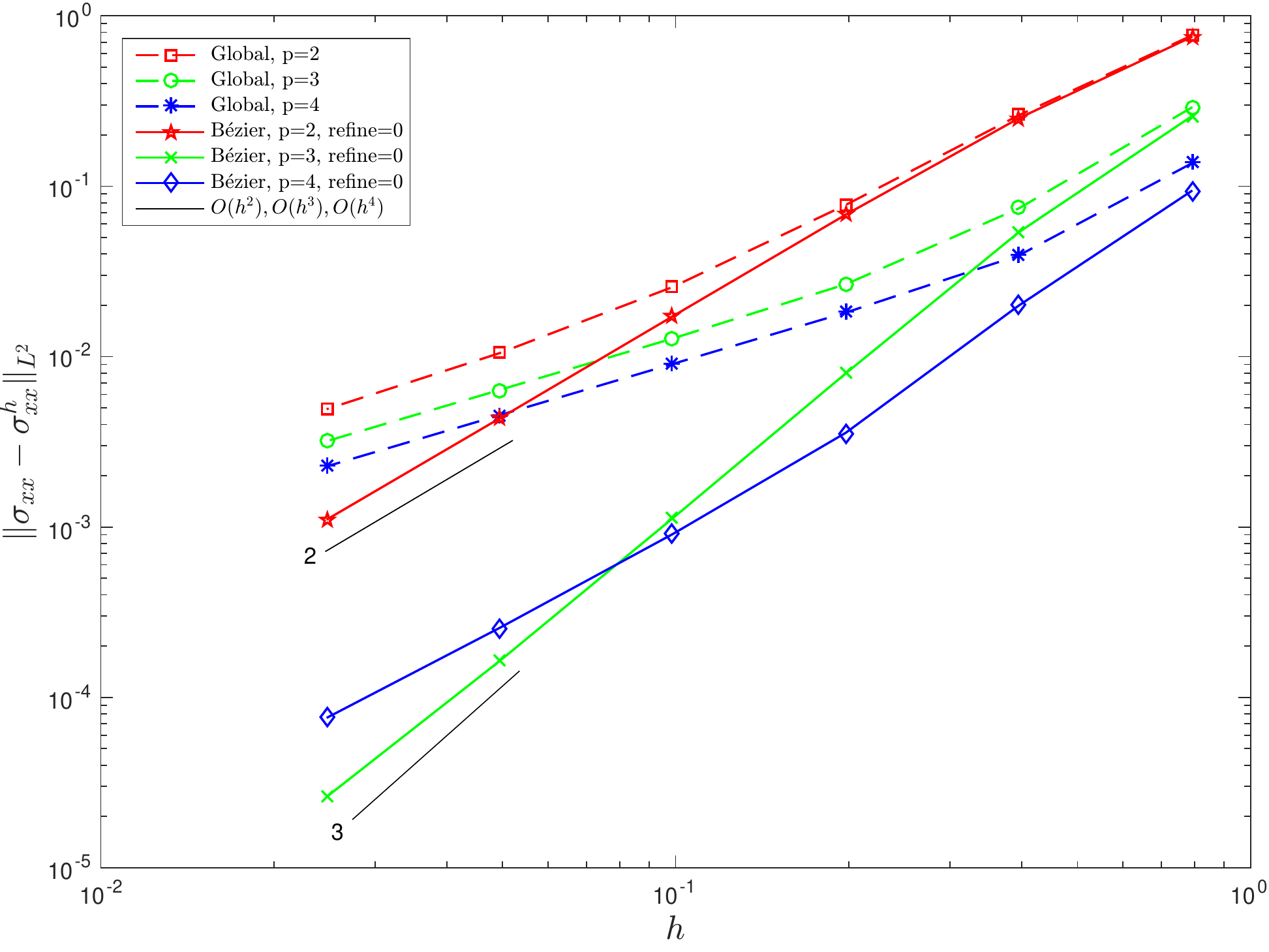} 
\caption{Without refinement}
\end{subfigure}\\
\begin{subfigure}{0.75\textwidth}
\includegraphics[width=\textwidth]{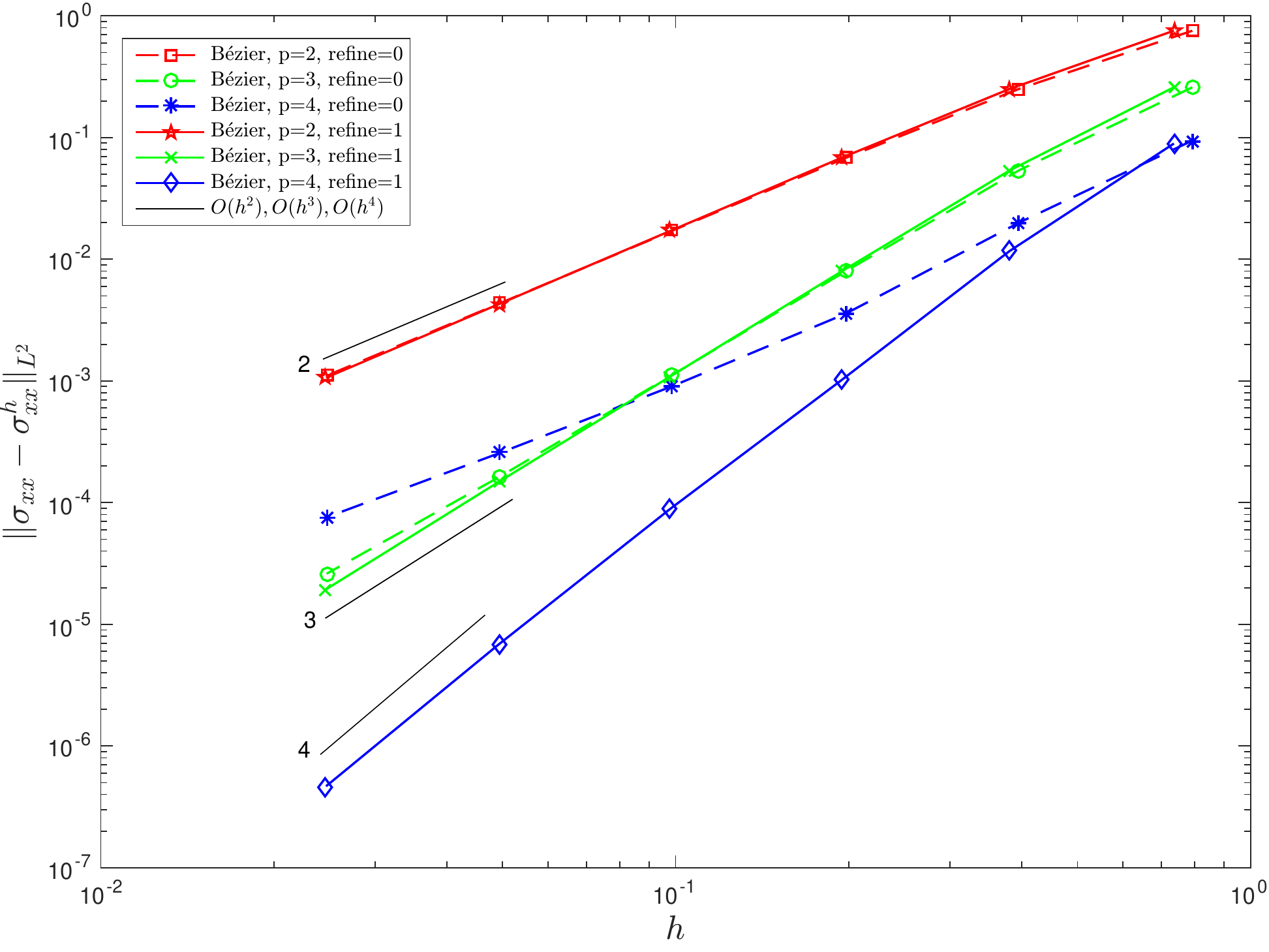} 
\caption{With refinement}
\end{subfigure} 
\caption{Stress convergence rates for a quarter plate with a circular hole
  decomposed into three nonconforming NURBS patches with matched
  parameterizations (see Figure~\ref{fig:solidinfiniteplatepatch3schematic}a).}
\label{fig:solid2Dinfiniteplatepatch3IGSE}
\end{figure}

\begin{figure}[h]
 \centering
  \begin{subfigure}{0.75\textwidth}
 \includegraphics[width=\textwidth]{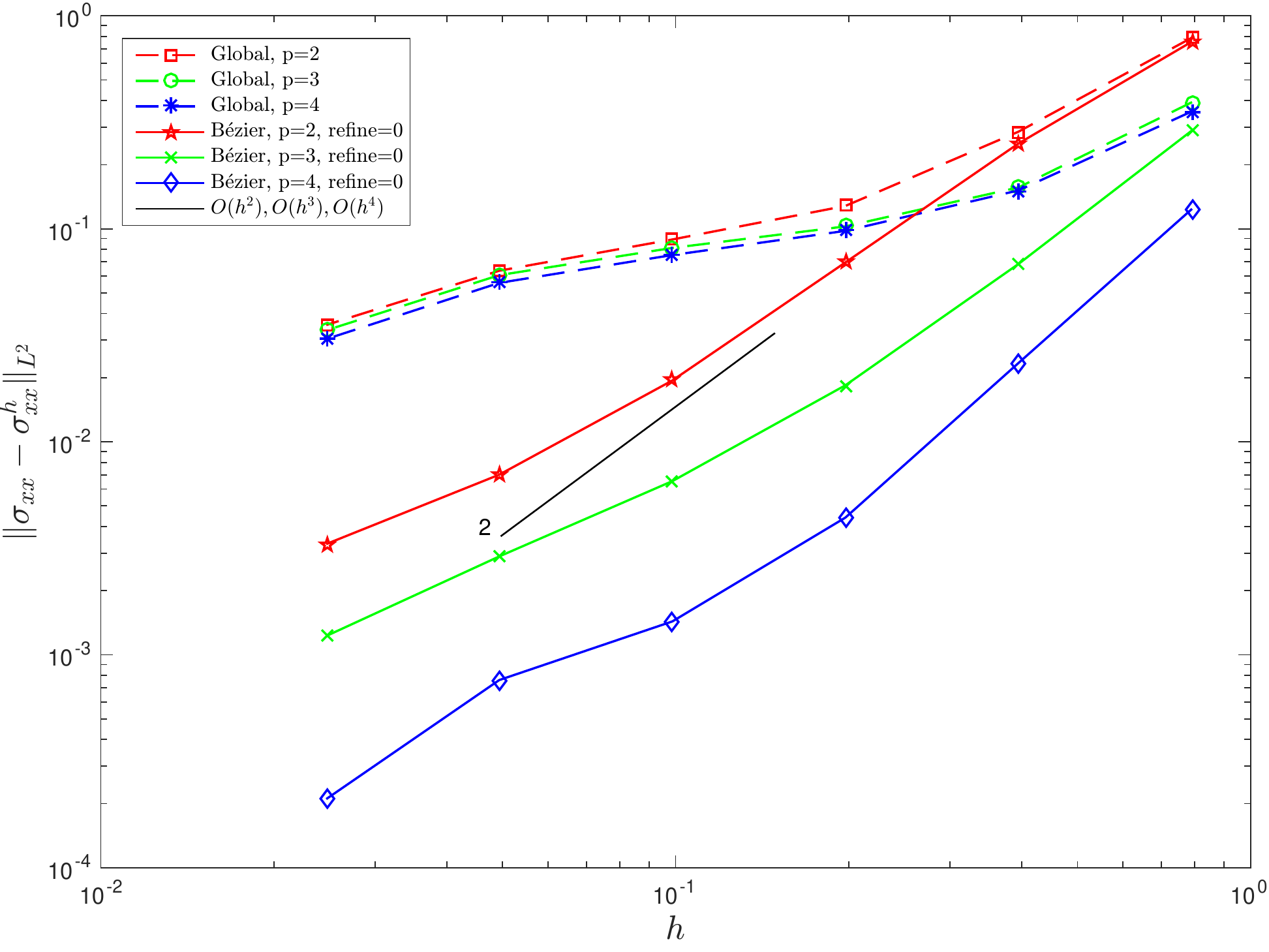} 
 \caption{Without refinement}
 \end{subfigure}\\
 \begin{subfigure}{0.75\textwidth}
 \includegraphics[width=\textwidth]{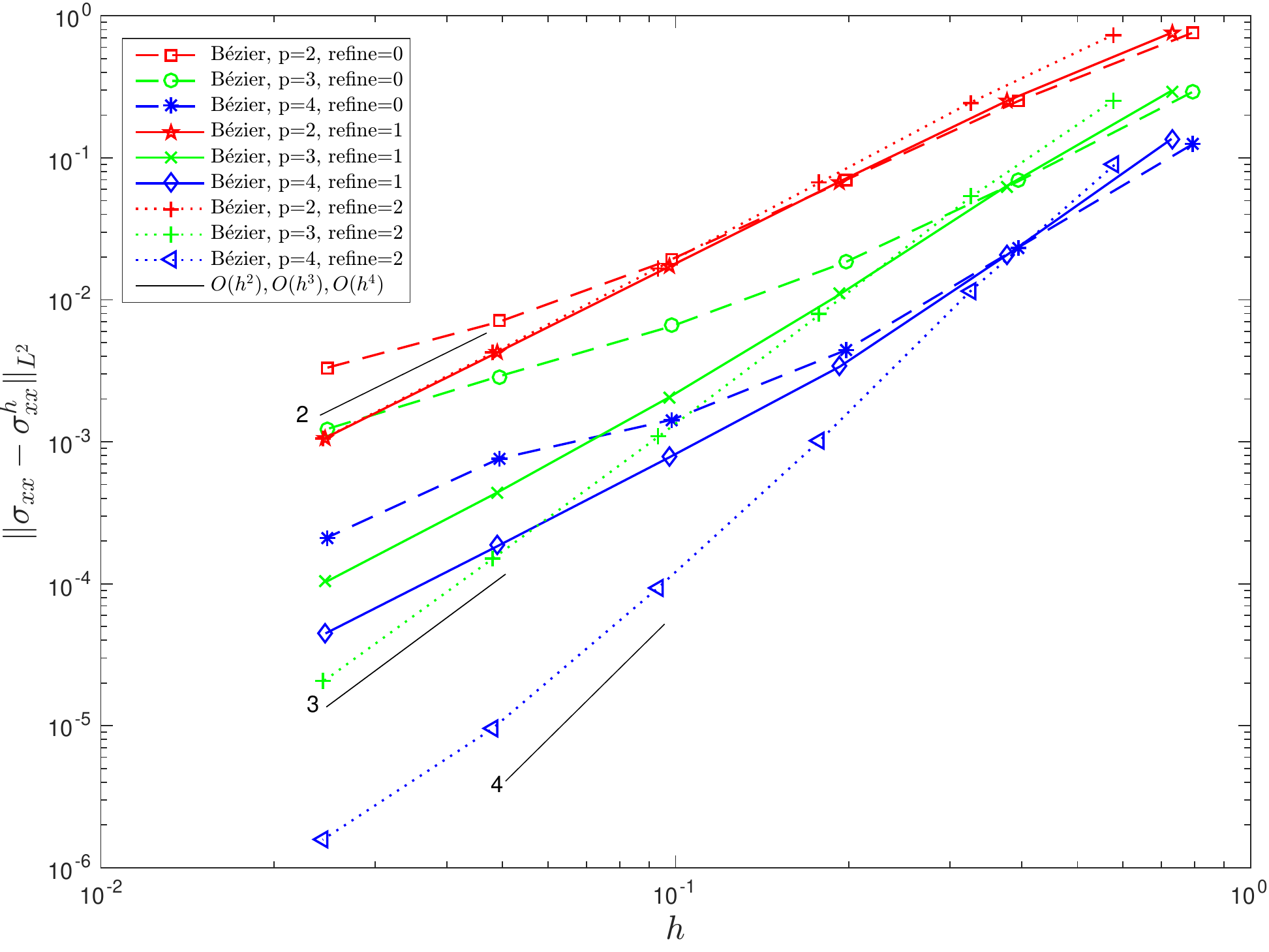} 
 \caption{With refinements}
 \end{subfigure}
\caption{Stress convergence rates for a quarter
  plate with a circular hole decomposed into three nonconforming NURBS patches
  with mismatched parameterizations (see Figure~\ref{fig:solidinfiniteplatepatch3schematic}b).}
\label{fig:solid2Dinfiniteplatepatch3IGSEmismatched}
\end{figure}

\begin{figure}[h]
  \begin{subfigure}{0.48\textwidth}
 \includegraphics[width=\textwidth, trim=1cm 1cm 1cm 0.5cm, clip=true]{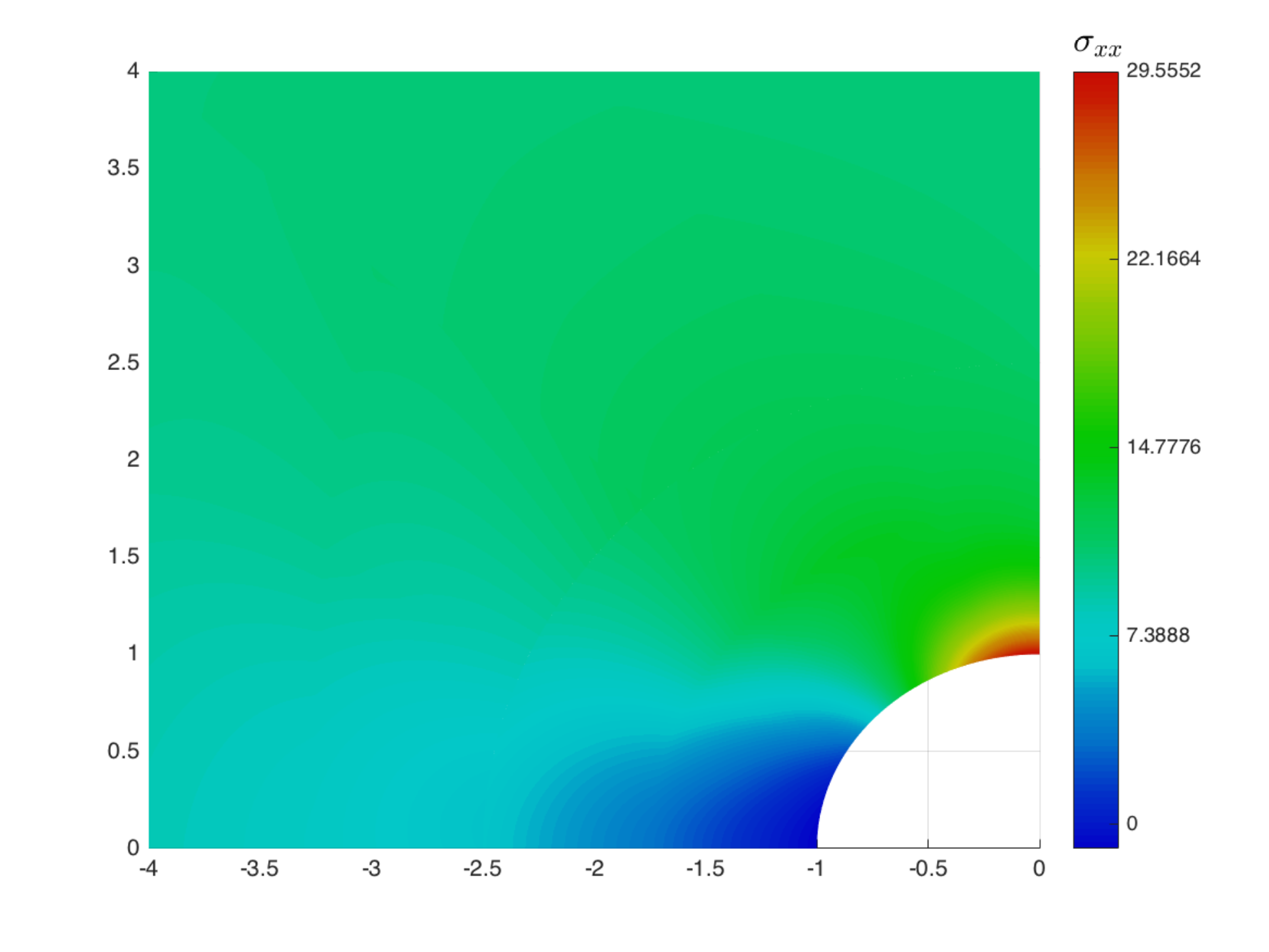} 
 \caption{ Initial mesh shown in Figure~\ref{fig:solidinfiniteplatepatch3schematic}a }
 \end{subfigure} 
  \begin{subfigure}{0.48\textwidth}
 \includegraphics[width=\textwidth, trim=1cm 1cm 1cm 0.5cm, clip=true]{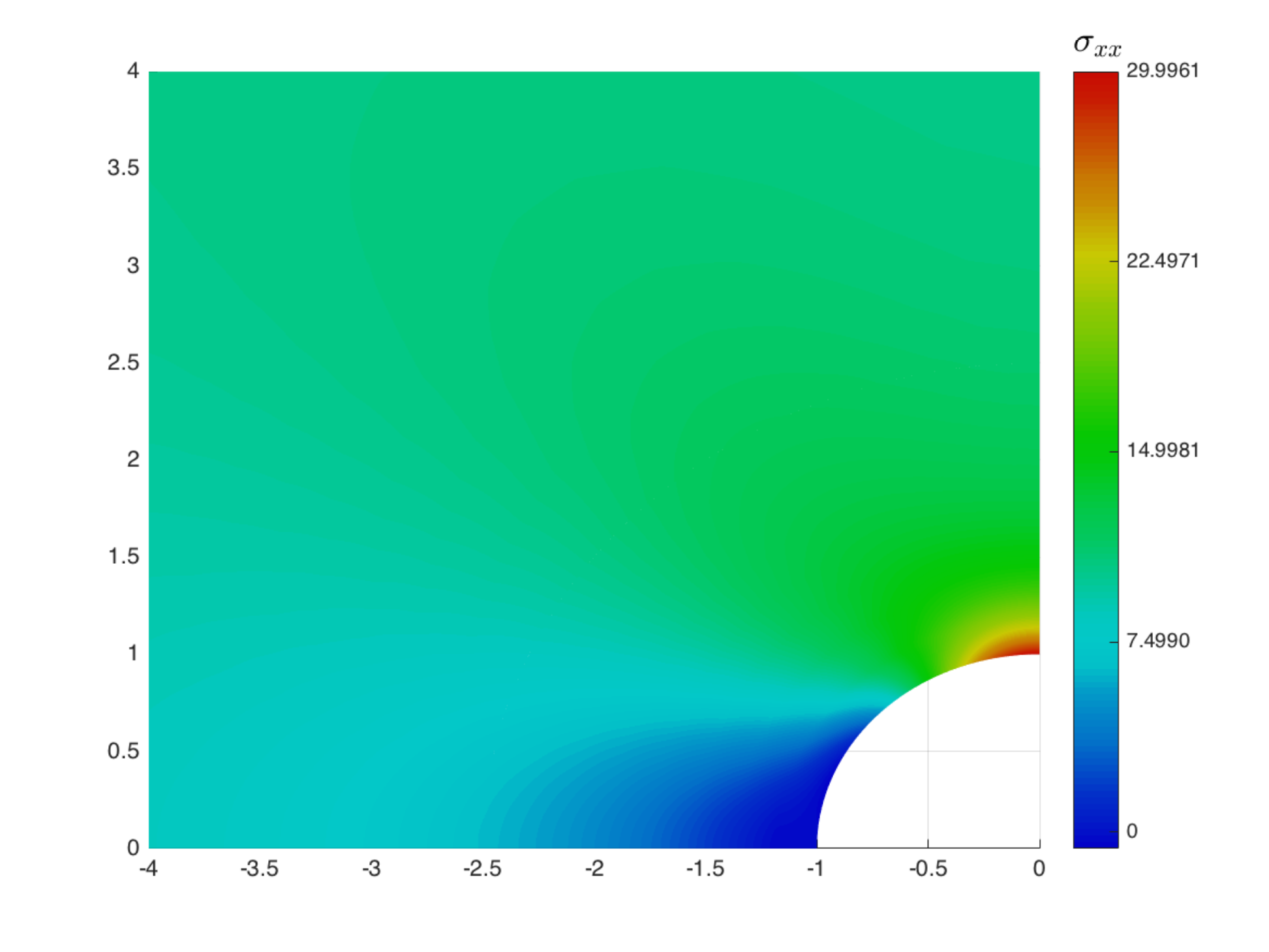} 
 \caption{ Mesh after 2 uniform refinements on (a)}
 \end{subfigure}\\
  \begin{subfigure}{0.48\textwidth}
 \includegraphics[width=\textwidth, trim=1cm 1cm 1cm 0.5cm, clip=true]{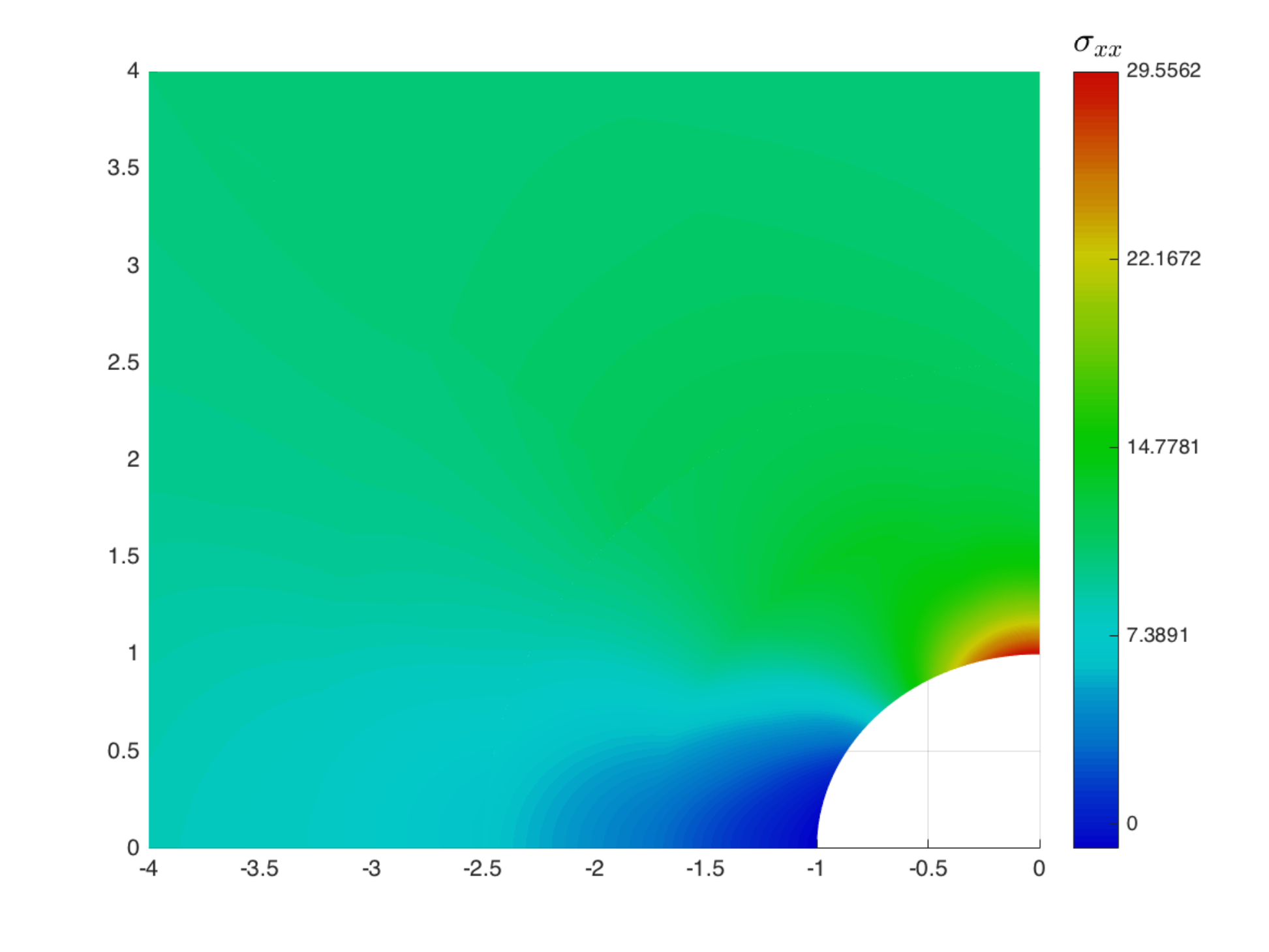} 
 \caption{ Initial mesh shown in Figure~\ref{fig:solidinfiniteplatepatch3schematic}b }
 \end{subfigure} 
  \begin{subfigure}{0.48\textwidth}
 \includegraphics[width=\textwidth, trim=1cm 1cm 1cm 0.5cm, clip=true]{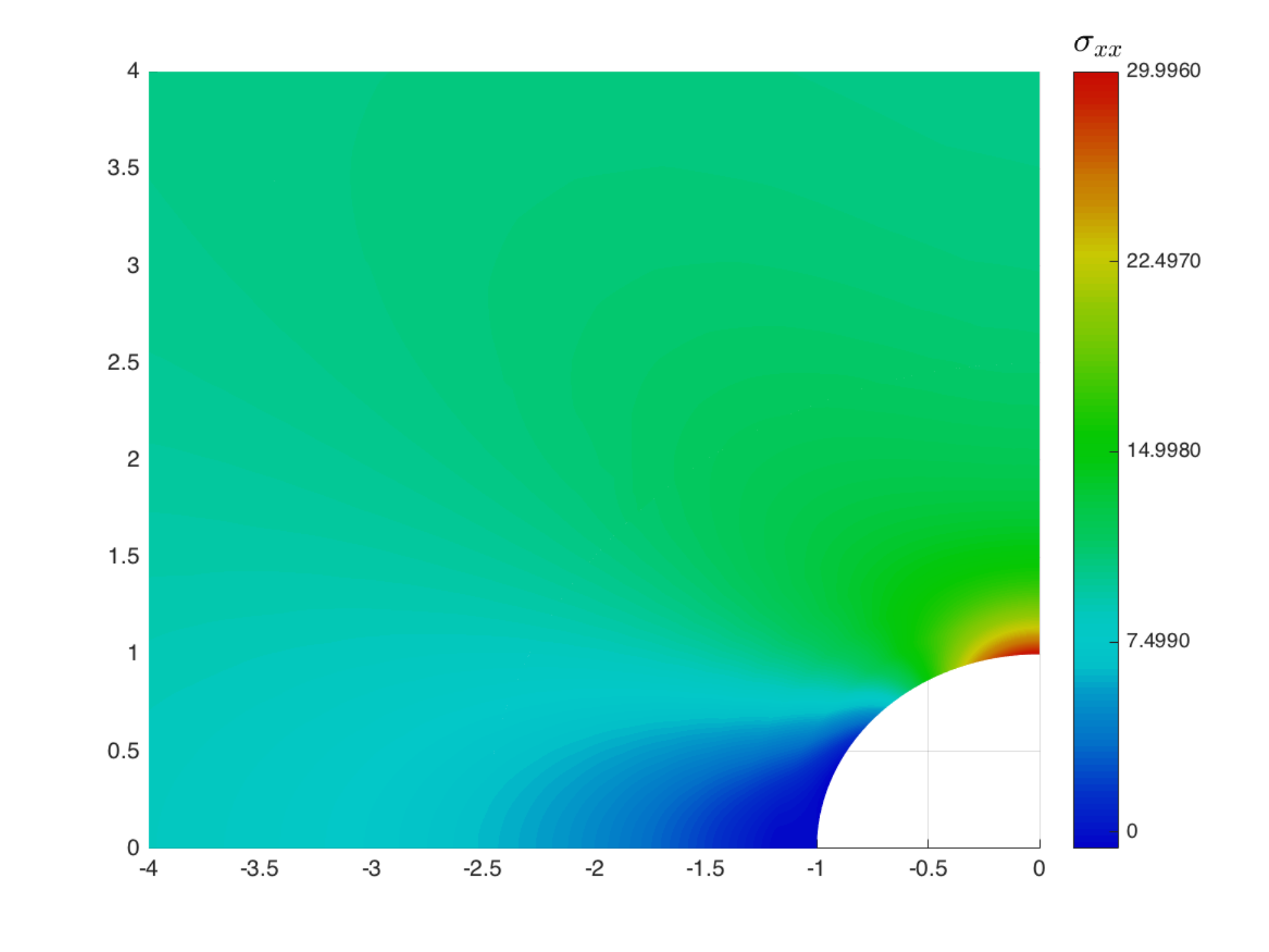} 
 \caption{ Mesh after 2 uniform refinements on (c)}
 \end{subfigure}\\
\caption{Stress $\sigma_{xx}$, $p=2$, no refinement of the dual basis space.}
\label{fig:solid2Dinfiniteplatepatch3mismatchedparameterizationstress}
\end{figure}

To assess the ability of the method to handle multiple patch coupling the
geometry is decomposed into three nonconforming NURBS patches with matched and
mismatched parameterizations as shown in
Figures~\ref{fig:solidinfiniteplatepatch3schematic}a and b, respectively.
As shown in Figures~\ref{fig:solid2Dinfiniteplatepatch3IGSE}a and \ref{fig:solid2Dinfiniteplatepatch3IGSEmismatched}a, the global dual mortar method
suffers from severely deteriorated convergence rates for $p=2,3$ and $4$ for both
matched and mismatched parameterizations due to two types of crosspoints, i.e., the
interface/Dirichlet intersections
and the interface/interface intersection in the middle. For matched
parameterizations without refinement, the \Bezier dual mortar method achieves optimal rate for
$p=2$, and slightly deteriorated rate for $p=3$ as
shown in Figure~\ref{fig:solid2Dinfiniteplatepatch3IGSE}a. Refining once recovers the optimal rates for $p=3,4$
(see Figure~\ref{fig:solid2Dinfiniteplatepatch3IGSE}b). For mismatched
parameterizations, optimal rates are achieved for $p=2$ by refining the dual
space once and for $p=3$ by refining the dual space twice (see
Figure~\ref{fig:solid2Dinfiniteplatepatch3IGSEmismatched}b).
Figure~\ref{fig:solid2Dinfiniteplatepatch3mismatchedparameterizationstress}
shows plots of the stress component $\sigma_{xx}$ for both matched and mismatched
parameterizations for $p=2$. It can be seen that even for the coarse initial meshes shown in
Figures~\ref{fig:solidinfiniteplatepatch3schematic}a and b, for both matched
and mismatched parameterizations, the stress concentration in the circular cutout
is very close to the analytical solution $\sigma_{xx} = 30$, as shown in
Figure~\ref{fig:solid2Dinfiniteplatepatch3mismatchedparameterizationstress}a and
Figure~\ref{fig:solid2Dinfiniteplatepatch3mismatchedparameterizationstress}b.

\FloatBarrier
\subsection{Weakly continuous geometry: Large deformations in two-dimensions}
\label{subsec:weak_con_geom}
We now employ \Bezier dual mortaring to build weakly continuous geometry as described in Section~\ref{sec:weak_geom}. Since the weak continuity constraint is embedded into the geometric description, a standard finite element code can be employed to process the weakly continuous basis in exactly the same manner as a standard conforming basis. To demonstrate the effectiveness of the approach, we compare the displacements computed on a weakly continuous mesh to those computed on a similar continuous mesh for a large deformation, plane strain problem.

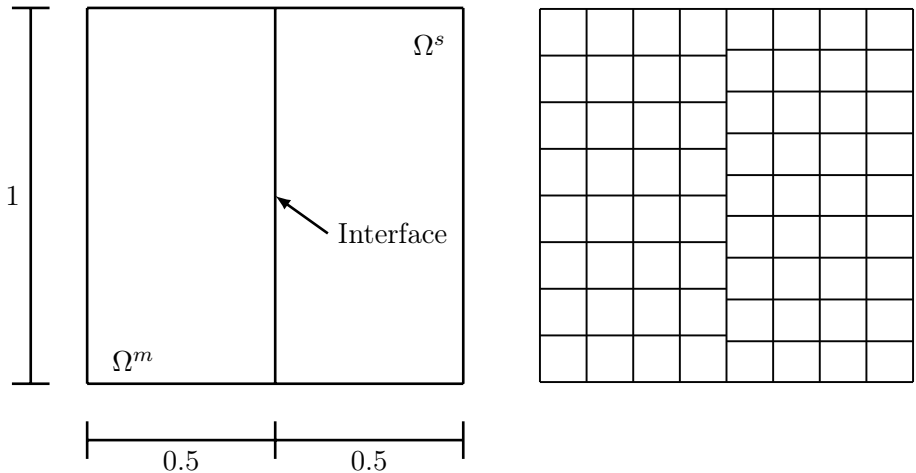
\begin{figure}[htb]
\centering
\begin{tikzpicture}
	\draw [line width = 1pt]
		(0,0) -- (5,0)
		(5,0) -- (5,5)
		(0,5) -- (5,5)
		(0,0) -- (0,5)
		(2.5,0) -- (2.5,5)
		;
		
	\draw (3.2,2) node[right] {Interface};
        \draw (0.2, 0.3 ) node[right]{$\Omega^m$};
        \draw (4.2, 4.5 ) node[right]{$\Omega^s$};
        
	\draw[-latex, line width = 1pt] (3.2,2) -- (2.5,2.5);
	
	\draw[line width = 1pt]
		(-0.5,0) -- (-1,0)
		(-0.5,5) -- (-1,5)
		(-0.75,0) -- (-0.75,5) node[pos=0.5,left] {1}
		
		(0,-0.5) -- (0,-1)
		(2.5,-0.5) -- (2.5,-1)
		(5,-0.5) -- (5,-1)
		
		(0,-0.75) -- (2.5,-0.75) node[pos=0.5,below]{0.5}
		(2.5,-0.75) -- (5,-0.75) node[pos=0.5,below]{0.5}
	;
	
	\node (myfirstpic) at (8.5,2.5) {\includegraphics[height=5cm]{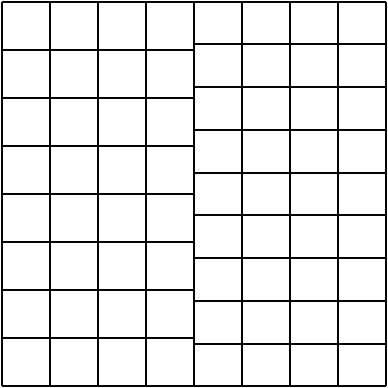}};
	
\end{tikzpicture}

\caption{The geometry and B\'ezier mesh.}
\label{fig:large-def-2d-geom-mesh}
\end{figure}

The initial geometry and the location of the interface are shown in Figure~\ref{fig:large-def-2d-geom-mesh}. For the weakly continuous mesh, the discretization does not match at the interface between the two patches and continuity is enforced weakly by building appropriately modified extraction operators. For the continuous mesh, the interface is treated as a $C^0$ interface. In both cases, the basis functions are quadratic maximally smooth B-splines. A B\'ezier element representation of the coarsest weakly continuous mesh is also shown in Figure~\ref{fig:large-def-2d-geom-mesh}. As can be seen, there is one additional element on the right side of the interface in the vertical direction. As the mesh is refined, the size of the elements on the left side of the interface are cut in half in each direction and the right side is refined so that there is always one additional element in the vertical direction. The continuous meshes are refined such that the element size is always the same as the element size on the left side of the weakly continuous meshes for a given refinement level.

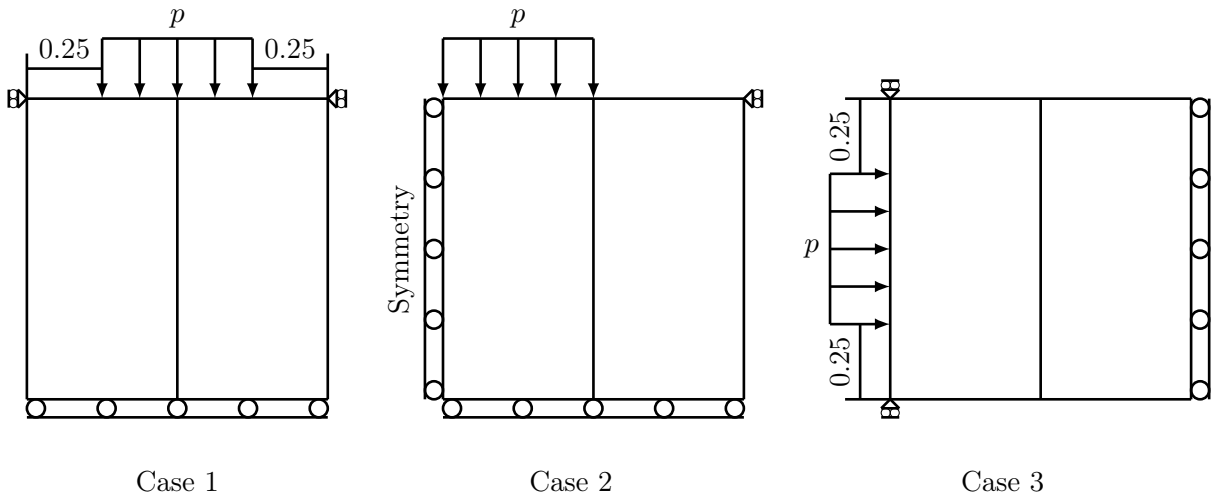
\begin{figure}[htb]
\centering
\begin{tabular}{ccc}
\begin{tikzpicture}[scale=0.8]
	\draw [line width = 1pt]
		(0,0) -- (5,0)
		(5,0) -- (5,5)
		(0,5) -- (5,5)
		(0,0) -- (0,5)
		(2.5,0) -- (2.5,5)
		;
		
	\draw[line width = 1pt] (1.25,6) -- (3.75,6) node[pos=0.5,above] {$p$};
	\draw[-latex, line width = 1pt] (1.25,6) -- (1.25,5);
	\draw[-latex, line width = 1pt] (1.875,6) -- (1.875,5);
	\draw[-latex, line width = 1pt] (2.5,6) -- (2.5,5);
	\draw[-latex, line width = 1pt] (3.125,6) -- (3.125,5);
	\draw[-latex, line width = 1pt] (3.75,6) -- (3.75,5);
	
	\draw[line width = 1pt] (0,5) -- (0,5.75);
	\draw[line width = 1pt] (0,5.5) -- (1.25,5.5) node[pos=0.5,above] {0.25};
	
	\draw[line width = 1pt] (5,5) -- (5,5.75);
	\draw[line width = 1pt] (3.75,5.5) -- (5,5.5) node[pos=0.5,above] {0.25};
	
	\draw[line width=1pt] (0.15,-0.15) circle (0.15cm);
	\draw[line width=1pt] (1.325,-0.15) circle (0.15cm);
	\draw[line width=1pt] (2.5,-0.15) circle (0.15cm);
	\draw[line width=1pt] (3.675,-0.15) circle (0.15cm);
	\draw[line width=1pt] (4.85,-0.15) circle (0.15cm);
	\draw[line width=1pt] (0,-0.30) -- (5,-0.30);
	
	\draw[line width=1pt] (-0.15,4.85) -- (0,5);
	\draw[line width=1pt] (-0.15,5.15) -- (0,5);
	\draw[line width=1pt] (-0.15,4.85) -- (-0.15,5.15);
	\draw (-0.225,4.925) circle (0.075cm);
	\draw (-0.225,5.075) circle (0.075cm);
	\draw[line width=1pt] (-0.30,4.85) -- (-0.30,5.15);

	\draw[line width=1pt] (5,5) -- (5.15,5.15);
	\draw[line width=1pt] (5,5) -- (5.15,4.85);
	\draw[line width=1pt] (5.15,5.15) -- (5.15,4.85);
	\draw (5.225, 4.925) circle (0.075cm);
	\draw (5.225, 5.075) circle (0.075cm);
	\draw[line width=1pt] (5.30,4.85) -- (5.30,5.15);

\end{tikzpicture}
& 
\begin{tikzpicture}[scale=0.8]
	\draw [line width = 1pt]
		(0,0) -- (5,0)
		(5,0) -- (5,5)
		(0,5) -- (5,5)
		(0,0) -- (0,5)
		(2.5,0) -- (2.5,5)
		;
		
	\draw[line width = 1pt] (0,6) -- (2.5,6) node[pos=0.5,above] {$p$};
	\draw[-latex, line width = 1pt] (0,6) -- (0,5);
	\draw[-latex, line width = 1pt] (0.625,6) -- (0.625,5);
	\draw[-latex, line width = 1pt] (1.25,6) -- (1.25,5);
	\draw[-latex, line width = 1pt] (1.875,6) -- (1.875,5);
	\draw[-latex, line width = 1pt] (2.5,6) -- (2.5,5);
	
	\draw[line width=1pt] (0.15,-0.15) circle (0.15cm);
	\draw[line width=1pt] (1.325,-0.15) circle (0.15cm);
	\draw[line width=1pt] (2.5,-0.15) circle (0.15cm);
	\draw[line width=1pt] (3.675,-0.15) circle (0.15cm);
	\draw[line width=1pt] (4.85,-0.15) circle (0.15cm);
	\draw[line width=1pt] (0,-0.30) -- (5,-0.30);
	

	\draw[line width=1pt] (5,5) -- (5.15,5.15);
	\draw[line width=1pt] (5,5) -- (5.15,4.85);
	\draw[line width=1pt] (5.15,5.15) -- (5.15,4.85);
	\draw (5.225, 4.925) circle (0.075cm);
	\draw (5.225, 5.075) circle (0.075cm);
	\draw[line width=1pt] (5.30,4.85) -- (5.30,5.15);

	\draw[line width=1pt] (-0.15,0.15) circle (0.15cm);
	\draw[line width=1pt] (-0.15,1.325) circle (0.15cm);
	\draw[line width=1pt] (-0.15,2.5) circle (0.15cm);
	\draw[line width=1pt] (-0.15,3.675) circle (0.15cm);
	\draw[line width=1pt] (-0.15,4.85) circle (0.15cm);
	\draw[line width=1pt] (-0.30,0) -- (-0.3,5) node[pos=0.5,above,rotate=90] {Symmetry};
\end{tikzpicture}
& 
\begin{tikzpicture}[scale=0.8]
	\draw [line width = 1pt]
		(0,0) -- (5,0)
		(5,0) -- (5,5)
		(0,5) -- (5,5)
		(0,0) -- (0,5)
		(2.5,0) -- (2.5,5)
		;
		
	\draw[line width = 1pt] (-1,1.25) -- (-1,3.75) node[pos=0.5,left] {$p$};
	\draw[-latex, line width = 1pt] (-1,1.25) -- (0,1.25);
	\draw[-latex, line width = 1pt] (-1,1.875) -- (0,1.8755);
	\draw[-latex, line width = 1pt] (-1,2.5) -- (0,2.5);
	\draw[-latex, line width = 1pt] (-1,3.125) -- (0,3.125);
	\draw[-latex, line width = 1pt] (-1,3.75) -- (0,3.75);
	
	\draw[line width = 1pt] (-0,0) -- (-0.75,0);
	\draw[line width = 1pt] (-0.5,0) -- (-0.5,1.25) node[pos=0.5,above,rotate=90] {0.25};
	
	\draw[line width = 1pt] (-0,5) -- (-0.75,5);
	\draw[line width = 1pt] (-0.5,3.75) -- (-0.5,5) node[pos=0.5,above,rotate=90] {0.25};

	\draw[line width=1pt] (5.15,0.15) circle (0.15cm);
	\draw[line width=1pt] (5.15,1.325) circle (0.15cm);
	\draw[line width=1pt] (5.15,2.5) circle (0.15cm);
	\draw[line width=1pt] (5.15,3.675) circle (0.15cm);
	\draw[line width=1pt] (5.15,4.85) circle (0.15cm);
	\draw[line width=1pt] (5.30,0) -- (5.3,5);
	
        
	\draw[line width=1pt] (0,5) -- (-0.15,5.15);
	\draw[line width=1pt] (0,5) -- (0.15,5.15);
	\draw[line width=1pt] (-0.15,5.15) -- (0.15,5.15);
	\draw (-0.075,5.225) circle (0.075cm);
	\draw (0.075,5.225) circle (0.075cm);
	\draw[line width=1pt] (-0.15,5.3) -- (0.15,5.3);

	\draw[line width=1pt] (0,0) -- (-0.15,-0.15);
	\draw[line width=1pt] (0,0) -- (0.15,-0.15);
	\draw[line width=1pt] (-0.15,-0.15) -- (0.15,-0.15);
	\draw (-0.075,-0.225) circle (0.075cm);
	\draw (0.075,-0.225) circle (0.075cm);
	\draw[line width=1pt] (-0.15,-0.3) -- (0.15,-0.3);
	
\end{tikzpicture}
\\ \\
Case 1 & Case 2 & Case 3
\end{tabular}
\caption{Load cases}
\label{fig:large-def-2d-loads}
\end{figure}

We compare the computed results for the three load cases with associated boundary conditions that are shown in Figure~\ref{fig:large-def-2d-loads}. The deformation is governed by the strain energy density functional that is given by
\begin{align}
	\psi = \lambda\left(\tfrac{1}{4}(J^2-1) - \tfrac{1}{2}\ln J\right) 
  + \tfrac{1}{2}\mu\left(\operatorname{tr}\left[\mathbf{b}\right]-3 - 2\ln J\right)
\end{align}
where $\lambda$ and $\mu$ are the typical Lam\'e parameters with
\begin{align}
	\lambda &= \frac{E\nu}{(1+\nu)(1-2\nu)}\\
	\mu &= \frac{E}{2(1+\nu)}
\end{align}
for Young's Modulus, $E$, and Poison's ratio, $\nu$. We use $E= 30 \times 10^9$
and $\nu=0.48$ for the results presented here. In addition, 
\begin{align}
  J = |\mathbf{F}| \quad \text{and} \quad \mathbf{b} = \mathbf{F} \mathbf{F}^{\text{T}}
\end{align}
where $\mathbf{F}$ is the deformation gradient, and $\mathbf{b}$ is the left
Cauchy-Green tensor.

The pressure boundary condition, $p$, is applied as a dead load in the reference configuration and is increased in twenty equal load increments to a maximum value of $100 \times 10^9$. At each load increment, the nonlinear problem is solved using a Newton-Raphson scheme with convergence satisfied when the residual is reduced by a factor of $10^8$.

\begin{figure}
\centering
	\begin{tabular}{cc}
		\includegraphics[scale=0.22]{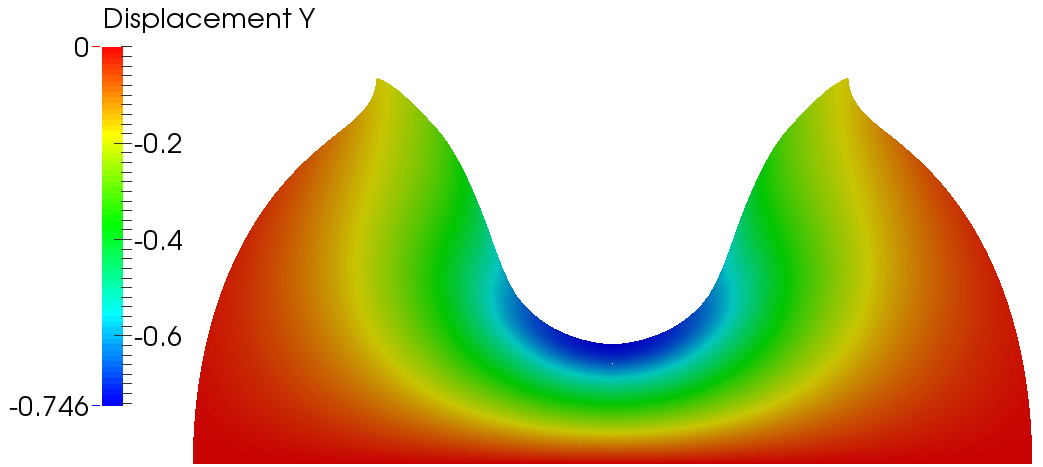} &
		\includegraphics[scale=0.22]{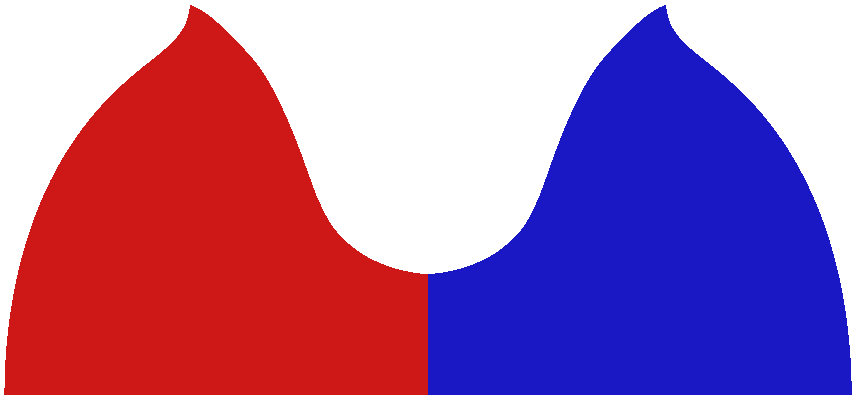}
	\end{tabular}
        \caption{Vertical displacement and deformed configuration - Case 1.}
\label{fig:large-def-2d-case1}
\end{figure}

\begin{figure}
\centering
	\begin{tabular}{cc}
		\includegraphics[scale=0.18]{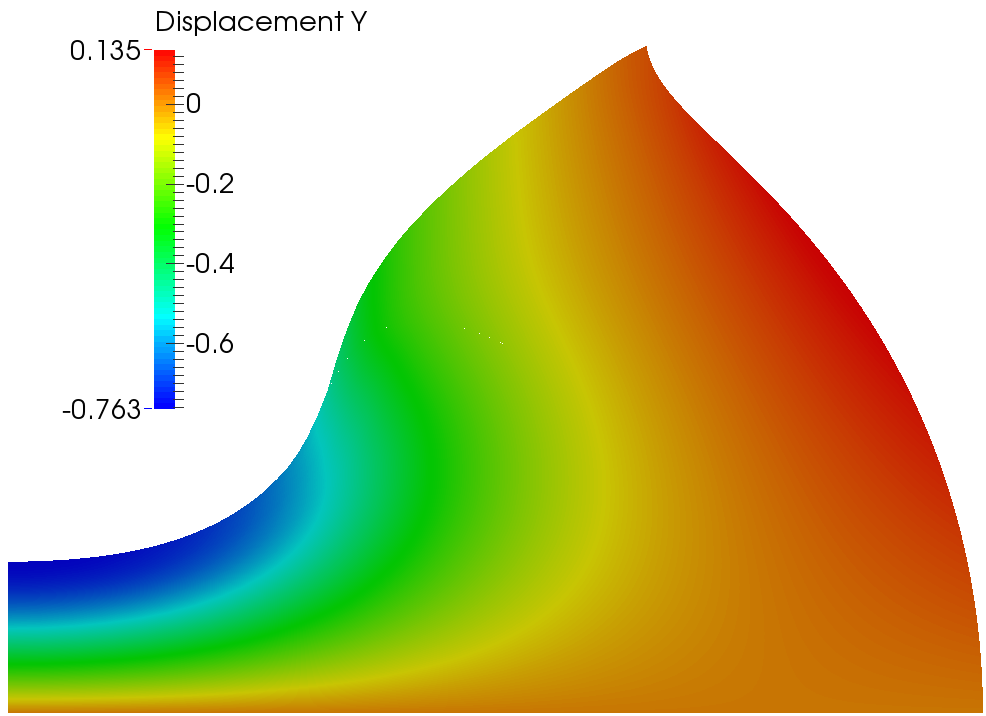} &
		\includegraphics[scale=0.18]{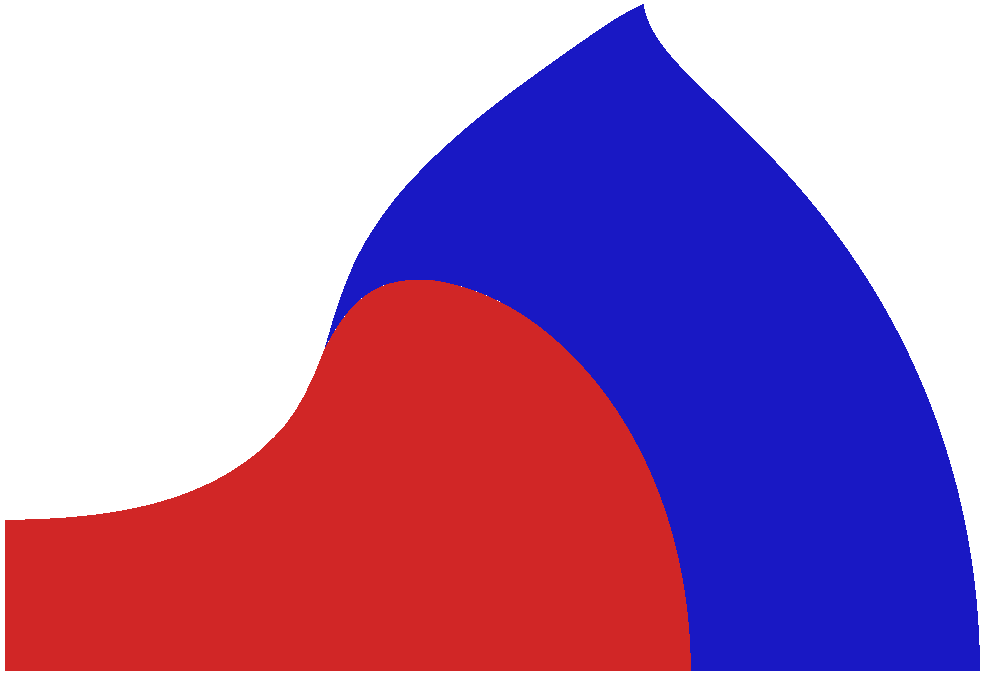}
	\end{tabular}
        \caption{Vertical displacement and deformed configuration - Case 2.}
\label{fig:large-def-2d-case2}
\end{figure}

\begin{figure}
\centering
	\begin{tabularx}{\textwidth}{>{\centering\arraybackslash}X>{\centering\arraybackslash}X}
		\includegraphics[scale=0.25]{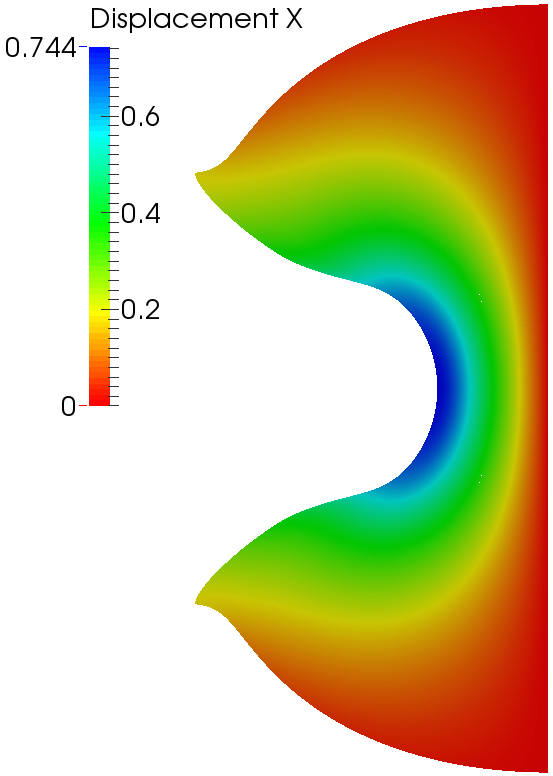} &
		\includegraphics[scale=0.25]{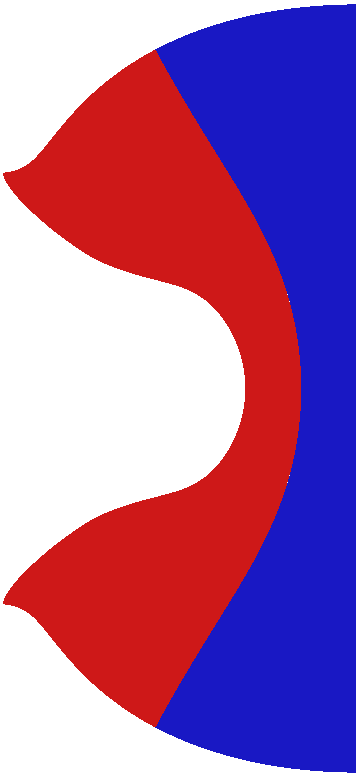}
	\end{tabularx}
        \caption{Horizontal displacement and deformed configuration - Case 3.}
\label{fig:large-def-2d-case3}
\end{figure}

The results of the computations are shown in Figures~\ref{fig:large-def-2d-case1} through~\ref{fig:large-def-2d-case3}. Each figure shows the unscaled deformation at the final load increment. On the left, the color scale indicates the magnitude of the displacement. On the right, each patch is shown as a distinct color so that the deformation of the interface between the two patches can clearly be seen. In all cases, the deformation of the interface is severe, but there is nothing in the displacement plot that indicates the presence of the weak interface.

\begin{figure}[!htb]
\centering
\includegraphics[width=0.5\textwidth]{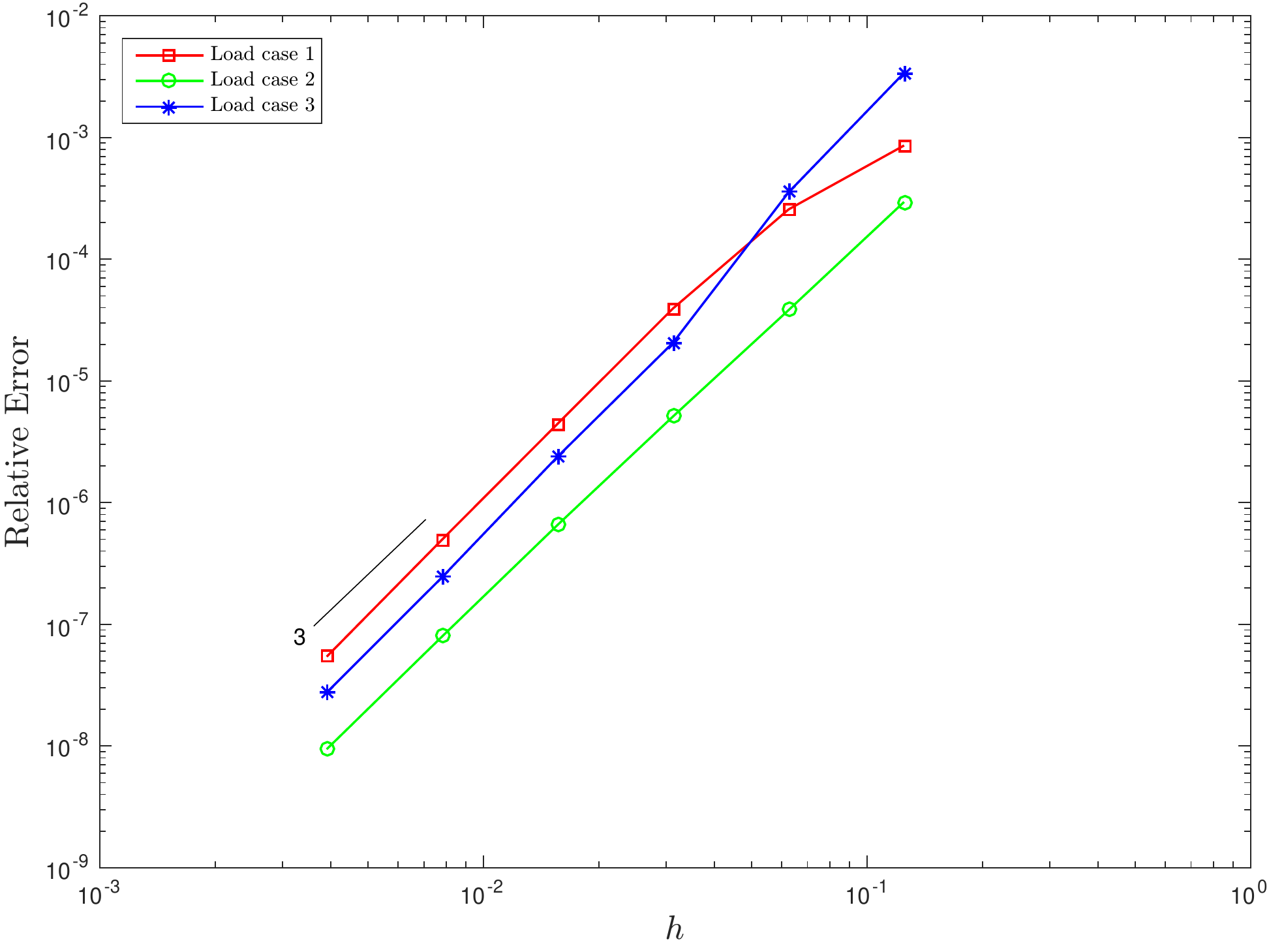}
\caption{Convergence rates of the $L^2$-relative error for the three large deformation load cases.}
\label{fig:large-def-2d-error}
\end{figure}

To quantify the accuracy of the weak geometry approach we compare the displacements, $u^{h,w}$, computed on the weakly continuous mesh to the displacements, $u^{h,c}$, computed on the continuous mesh. We define the relative error, $e_r$, to be the $L^2$-norm of the difference between the two considered displacements, i.e
\begin{align}
	e_r = \|u^{h,w} - u^{h,s}\|_{L^2}.
\end{align}
Using the triangle inequality,
\begin{align}
	\|u^{h,w}-u\|_{L^2} \leq \|u^{h,w} - u^{h,c}\|_{L^2} + \|u^{h,c}-u\|_{L^2}, \label{eq:triangle-inequality-large-def}
\end{align}
we see that the absolute error of the solution computed on the weakly continuous mesh case is bounded by the sum of the relative error and the absolute error of the solution computed on the continuous mesh case. Now, assuming that the solution computed on the continuous mesh case converges optimally, by \eqref{eq:triangle-inequality-large-def}, if the relative error converges optimally then we know that the absolute error of the solution computed on the weakly continuous mesh case also must converge optimally. The convergence rates of the relative error are plotted in Figure~\ref{fig:large-def-2d-error} for the three load cases. This figure clearly shows that the convergence rate of the relative error are cubic, which is the optimal rate for quadratic basis functions.

\section{Conclusions}
\label{sec:conclude}
We have introduced a new approach for the coupling of non-conforming higher-order smooth spline patches which
we call the isogeometric \Bezier dual mortar method. The construction of the underlying dual spline basis is based on \Bezier
extraction and projection and is applicable to any spline description which has
a \Bezier representation. The dual basis is refineable and the associated mortaring strategy preserves the sparsity of the stiffness matrix. The accuracy of the
coupling can be adaptively controlled by employing a dual basis refinement
scheme which can be used to recover optimal convergence rates without adding any
additional degrees-of-freedom to the global system. As a particular application of \Bezier dual
mortaring, we introduced weakly continuous geometry, where the weak continuity constraint is built into properly modified extraction operators. This allows
for the use of weakly coupled multi-patch geometry in design and as a basis for
standard finite element frameworks which do not employ any mortaring algorithms.

 We applied the isogemetric \Bezier dual mortar method to standard linear and nonlinear elastic test cases and B-spline and NURBS geometries. All tests show that the isogeometric
 \Bezier dual mortar method is robust and accurate, works for arbitrary
 master/slave pairings, and arbitrary parameterizations.
 \section*{Acknowledgment}
 M. A. Scott was supported through a grant from the Air Force Office of
 Scientific Research (FA9550-214-1-0113) and the Ford University Research
 Program, and W. Dornisch was partially supported by the German Research
 Foundation (DFG) through the Research Group FOR 1509 and the Collaborative
 Research Centre SFB 926. These supports are gratefully acknowledged.
 
\newpage
\section*{References}
\bibliographystyle{elsarticle-num-names}
\bibliography{bibliography1,bibliography2}
\newpage
\appendix
\section{A derivation of the weakly continuous extraction operator for element $e_{21}$ from Figure~\ref{fig:schematics_of_spliting_element}.}
\label{app:weak_extra_op}
For the example shown in Figure~\ref{fig:schematics_of_spliting_element}, the basis relation matrix 
  $\mathbf{G}_{\bar{N}^r, N^m}$ defined by (\ref{eq:global_basis_relation}) and the
  localized counterpart $\mathbf{G}^e_{\bar{N}^r, N^m}$ for the interface of element $e_{21}$ are 
   \begin{align}
     \mathbf{G}_{\bar{N}^r, N^m} =
     \begin{bmatrix}
       1& \tikzmark{De1top}{\frac{1}{3}}& 0          & 0                                & 0          & 0\\[0.5em]
       0& \frac{2}{3}                   & \frac{2}{3}& \frac{1}{3}                      & 0          & 0\\[0.5em]
       0& 0                             & \frac{1}{3}& \tikzmark{De1bottom}{\frac{2}{3}}& \frac{2}{3}& 0\\[0.5em]
       0& 0                             & 0          & 0                                & \frac{1}{3}& 1
  \end{bmatrix}^{\text{T}}
 \quad \text{and} \quad
 \mathbf{G}^e_{\bar{N}^r, N^m} = 
  \begin{bmatrix}
        \frac{1}{3}& \frac{2}{3}           & 0          \\[0.5em]
        0          & \frac{2}{3}& \frac{1}{3}\\[0.5em]
        0          & \frac{1}{3}& \frac{2}{3}\\
  \end{bmatrix}.
 \begin{tikzpicture}[overlay,remember picture] 
   \fill[ opacity=.4] (De1top.north west) rectangle (De1bottom.south east);
 \end{tikzpicture}
   \end{align}
  The standard \Bezier extraction operator $\mathbf{R}^{e,r}_{\xi_1}$ and the basis
  transformation matrix $\mathbf{M}$ are
   \begin{align}
 \mathbf{R}^{e,r}_{\xi_1} = 
   \begin{bmatrix}
     \frac{1}{3}& 0& 0\\[0.5em]
     \frac{2}{3}& 1& \frac{1}{2}\\[0.5em]
     0          & 0& \frac{1}{2}
   \end{bmatrix}
   \quad \text{and} \quad
    \mathbf{M} =
   \begin{bmatrix}
     1          & 0          & 0\\[0.5em]
     \frac{1}{2}& \frac{1}{2}& 0\\[0.5em]
     \frac{1}{4}& \frac{1}{2}& \frac{1}{4} 
     \end{bmatrix},
   \end{align}
   and the weakly continuous one-dimensional interface element extraction operator $\tilde{\mathbf{R}}^e_{\xi_1}$ for element $e_{21}$ is
  \begin{align}
  \tilde{\mathbf{R}}^e_{\xi_1} = (\mathbf{G}^e_{\bar{N}^r, N^m})^\text{T} \mathbf{R}^{e,r}_{\xi_1} \mathbf{M}^{-\text{T}} =
    \begin{bmatrix}
      \frac{1}{9} & - \frac{1}{9} & \frac{1}{9}\\[0.5em]
      \frac{2}{3} & \frac{2}{3} & 0\\[0.5em]
      \frac{2}{9} & \frac{4}{9} & \frac{8}{9} 
    \end{bmatrix}.
    \end{align}
The two standard one-dimensional \Bezier extraction
  operators for the original slave patch element $e_2$ are
  \begin{align}
    \mathbf{R}^e_{\xi_1} = 
  \begin{bmatrix}
    \frac{1}{2}& 0& 0\\[0.5em]
    \frac{1}{2}& 1& \frac{1}{2}\\[0.5em]
    0          & 0& \frac{1}{2}
  \end{bmatrix}
   \quad \text{and} \quad
    \mathbf{R}^e_{\xi_2} = 
  \begin{bmatrix}
    \frac{1}{2}& 0& 0\\[0.5em]
    \frac{1}{2}& 1& 0\\[0.5em]
    0          & 0& 1
  \end{bmatrix}.
\end{align}
  As shown in Figure~\ref{fig:schematics_of_spliting_element}, the interior
  basis functions of element $e_{21}$ are identical to those of element $e_2$,
  and only the interface basis functions are replaced by the refined interface
  basis. We decompose $\mathbf{R}^e_{\xi_2}$ into two submatrices
  $\mathbf{R}_1$ and $\mathbf{R}_2$ such that
  {\setstretch{1.35}
  \begin{align}
    \mathbf{R}^e_{\xi_2} = 
    \begin{lbmatrix}{1}
      \mathbf{R}_1 \\  \hdashline[2pt/2pt] 
      \mathbf{R}_2  
    \end{lbmatrix} =
  \begin{lbmatrix}{3}
    \frac{1}{2}& 0& 0\\
    \frac{1}{2}& 1& 0\\ \hdashline[2pt/2pt]
    0          & 0& 1
  \end{lbmatrix},
    \end{align}}
 where $\mathbf{R}_1$ is related to the interior basis functions and
 $\mathbf{R}_2$ is related to the interface basis functions. Then, the weakly continuous patch element extraction operator $\tilde{\mathbf{R}}^e$ for element $e_{21}$ can be computed as follows:
  \begin{align}
    \tilde{\mathbf{R}}^e =
    \begin{bmatrix}
      \mathbf{R}_1 \otimes \mathbf{R}^e_{\xi_1} \\[0.25em]
      \mathbf{R}_2 \otimes \tilde{\mathbf{R}}^e_{\xi_1}
    \end{bmatrix} =
    \begin{bmatrix}
      \frac{1}{2} \mathbf{R}^e_{\xi_1}& 0               & 0\\[0.5em]
      \frac{1}{2} \mathbf{R}^e_{\xi_1}& \mathbf{R}^e_{\xi_1}&0\\[0.5em]
    0                 & 0               &  \tilde{\mathbf{R}}^e_{\xi_1}
    \end{bmatrix} =
     \begin{bmatrix}
       \frac{1}{4}& 0          & 0          & 0          & 0& 0          & 0          & 0            & 0          \\[0.5em]
       \frac{1}{4}& \frac{1}{2}& \frac{1}{4}& 0          & 0& 0          & 0          & 0            & 0          \\[0.5em]
       0          & 0          & \frac{1}{4}& 0          & 0& 0          & 0          & 0            & 0          \\[0.5em]
       \frac{1}{4}& 0          & 0          & \frac{1}{2}& 0& 0          & 0          & 0            & 0          \\[0.5em]
       \frac{1}{4}& \frac{1}{2}& \frac{1}{4}& \frac{1}{2}& 1& \frac{1}{2}& 0          & 0            & 0          \\[0.5em]
       0          & 0          & \frac{1}{4}& 0          & 0& \frac{1}{2}& 0          & 0            & 0          \\[0.5em]
       0          & 0          & 0          & 0          & 0& 0          & \frac{1}{9}& - \frac{1}{9}& \frac{1}{9}\\[0.5em]
       0          & 0          & 0          & 0          & 0& 0          & \frac{2}{3}& \frac{2}{3}  & 0          \\[0.5em]
       0          & 0          & 0          & 0          & 0& 0          & \frac{2}{9}& \frac{4}{9}  & \frac{8}{9} 
      \end{bmatrix}.
    \end{align}
Note that the only difference between the weakly continuous element extraction operator
and the standard element extraction operator is that the last three rows are
modified. These rows correspond to interface basis functions. 
\end{document}